\documentclass[11pt]{amsart}

\usepackage{import}
\usepackage{cancel}

\newcommand{\N}{\mathbb{N}} 
\newcommand{\Z}{\mathbb{Z}} 
\newcommand{\R}{\mathbb{R}} 

\newcommand{\id}{\mathrm{id}} 
\newcommand{\1}{\mathbf{1}}
\DeclareMathOperator{\im}{im}

\newcommand{\del}{\partial} 
\newcommand{\Cinf}{C^\infty} 

\newcommand{\derpar}[2]{\frac{\partial #1}{\partial #2}}
\newcommand{\nabladerpar}[2]{\frac{\nabla #1}{\partial #2}}

\newcommand{\abs}[1]{\vert#1\vert}
\newcommand{\norm}[1]{\left\Vert#1\right\Vert}

\newcommand{\Ascr}{\mathscr{A}}
\newcommand{\Bscr}{\mathscr{B}}
\newcommand{\Mscrh}{\hspace{-2pt}\mathscr{M}}
\newcommand{\Mscr}{\mathscr{M}}
\newcommand{\Escr}{\mathscr{E}}
\newcommand{\Fscr}{\mathscr{F}}

\newcommand{\Nscr}{\mathscr{N}}
\newcommand{\Nscrh}{\hspace{-2pt}\mathscr{N}}

\newcommand{\Homscr}{\mathscr{H\hspace{-4pt}o\hspace{-1pt}m}}

\newcommand{\Hom}{\mathrm{Hom}} 
\newcommand{\B}{\mathrm{B}}
\newcommand{\Sh}{\mathrm{Sh}}
\newcommand{\Der}{\mathrm{Der}}

\newcommand{\ev}{\mathrm{ev}}
\newcommand{\ad}{\operatorname{ad}}
\newcommand{\Pair}{\mathrm{Pair}}
\newcommand{\Chat}{\hat{C}}
\newcommand{\Homhat}{\mathrm{\widehat{H}om}}
\newcommand{\Endhat}{\mathrm{\widehat{E}nd}}

\newcommand{\sk}{\mathrm{sk}}

\newcommand{\Hol}{\mathrm{Hol}}

\newcommand{\End}{\mathrm{End}}
\newcommand{\Curv}{\mathrm{Curv}}
\newcommand{\GL}{\mathrm{GL}}

\newcommand{\Gau}{\mathrm{G}}

\newcommand{\Ner}{\mathrm{N}}

\newcommand{\Endi}{\mathrm{End}_\infty}
\newcommand{\Curvi}{\mathrm{Curv}_\infty}
\newcommand{\GLi}{\mathrm{GL}_\infty}

\newcommand{\Gaui}{\mathrm{G}_\infty}
\newcommand{\qGaui}{\mathrm{qG}_\infty}
\newcommand{\Neri}{\mathrm{N}_\infty}

\newcommand{\HConn}{\mathrm{Conn}^\infty}
\newcommand{\HRep}{\mathrm{Rep}^\infty}
\newcommand{\Rep}{\mathrm{Rep}}
\newcommand{\Flat}{\mathrm{Flat}}
\newcommand{\RH}{\mathrm{RH}}

\newcommand{\sSet}{\mathrm{sSet}}

\newcommand{\sep}{:}

\calclayout

\makeatletter
\newtheorem*{rep@theorem}{\rep@title}\newcommand{\newreptheorem}[2]{
\newenvironment{rep#1}[1]{
    \def\rep@title{#2 \ref{##1}}
    \begin{rep@theorem}}
    {\end{rep@theorem}}
}
\makeatother

\newtheorem{theorem}{Theorem}[section]
\newtheorem*{theorem*}{Theorem}
\newreptheorem{theorem}{Theorem}
\newtheorem{corollary}[theorem]{Corollary}
\newtheorem{lemma}[theorem]{Lemma}
\newtheorem*{lemma*}{Lemma}
\newtheorem{proposition}[theorem]{Proposition}

\newtheorem*{question*}{Question}

\theoremstyle{definition}
\newtheorem{definition}[theorem]{Definition}
\newtheorem{example}[theorem]{Example}
\AtBeginEnvironment{example}{%
  \pushQED{\qed}%
}
\AtEndEnvironment{example}{\popQED}
\theoremstyle{remark}
\newtheorem{remark}[theorem]{Remark}
\AtBeginEnvironment{remark}{%
  \pushQED{\qed}%
}
\AtEndEnvironment{remark}{\popQED}
\newtheorem{notation}[theorem]{Notation}

\usepackage{mathrsfs}
\usepackage[
	scr=boondox,
	cal=euler 
]{mathalpha}
\usepackage{amssymb}
\usepackage{amsmath}
\usepackage{dsfont}
\usepackage{tensor}
\usepackage{esint}
\usepackage{mathtools}
\usepackage{svg}

\usepackage{xcolor}
\definecolor{color-citas}{RGB}{9,21,128}

\usepackage{hyperref}
\hypersetup{
    colorlinks=true,
    linkcolor=color-citas,
    citecolor=color-citas,
    urlcolor=color-citas
}
\usepackage{cleveref}
\usepackage{enumitem}
\usepackage{tikz-cd}
\svgsetup{inkscapearea=page}

\allowdisplaybreaks

\numberwithin{equation}{section}

\begin{document} 

\title[Higher Holonomy]{Higher Holonomy for Higher Lie Algebroid Connections on Graded Vector Bundles}
\author{Sven Holtrop}
\author{Jaime Pedregal}
\address{Department of Mathematics, Utrecht University, 3508 TA Utrecht, The Netherlands}
\email{s.t.holtrop@uu.nl, j.pedregalpastor@uu.nl}

\thanks{S.H. has been supported by NWO Open Competition grant OCENW.M.22.190 (\url{https://doi.org/10.61686/YTCMW72497}) and J.P. by the NWO Open Competition grant OCENW.M.22.264 (\url{https://doi.org/10.61686/XBDQH34231}) and by FEDER/AEI/MICINN through the grants PID2025-173178NB-I00 and PID2022-137667NA-I00.}

\begin{abstract}
We prove the existence of a Riemann--Hilbert $A_\infty$-functor from the curved dg-category of $\Ascr$-connections on graded vector bundles, where $\Ascr$ is a higher Lie algebroid, to the curved dg-category of $\Pi_\infty(\Ascr)$-connections on graded vector bundles, where $\Pi_\infty(\Ascr)$ is the homotopy groupoid of $\Ascr$. When restricted to the sub dg-categories of flat objects, it induces an $A_\infty$-functor from the dg-category of representations up to homotopy of $\Ascr$ to the dg-category of representations up to homotopy of $\Pi_\infty(\Ascr)$, which recovers the one from \cite{arias-abad-schatz,block-smith} when $\Ascr$ is a Lie algebroid. The key ingredient for defining such functors is a non-abelian Stokes's theorem for endomorphism-valued differential forms. Furthermore, for a graded vector bundle $E$, we introduce the notion of its higher general linear groupoid in such a way that $G$-connections on $E$, for a higher Lie groupoid $G$, become smooth functors from $G$ to the higher general linear groupoid of $E$. We combine this with the Riemann--Hilbert functor to define higher parallel transport functors for $\Ascr$-connections. Finally, we use these functors to define the higher holonomy groups of $\Ascr$-connections and endow them with the structure of augmented simplicial Lie groups.
\end{abstract}

\maketitle

\tableofcontents

\bigskip
\section{Introduction}

{\renewcommand{\thetheorem}{\Alph{theorem}}
\noindent It is well known that a flat connection on a vector bundle $E\to M$ induces a representation of $\pi_1(M,x)$ on $E_x$ for any $x\in M$, since parallel transport only depends on the homotopy class of paths. Actually, the celebrated Riemann--Hilbert correspondence states that there is an equivalence of categories
\[
\Flat(M) \cong \Rep(\pi_1(M,x)),
\]
where $\Flat(M)$ is the category of flat connections on vector bundles over $M$. Moreover, by fixing the vector bundle $E\to M$ and considering $\Flat(M,E)$, the category of flat connections on $E$, one gets an additional equivalence
\[
\Flat(M,E) \cong \Rep(\pi_1(M,x),E_x)\cong \Hom(\pi_1(M,x),\GL(E_x)),
\]
where $\Rep(\pi_1(M,x),E_x)$ refers to the category of $\pi_1(M,x)$-representations on the vector space $E_x$.

It is possible to get rid of the dependence on the base point by considering the relevant groupoids:
\[
\Flat(M) \cong \Rep(\Pi_1(M)),
\]
where $\Pi_1(M)\rightrightarrows M$ is the fundamental groupoid. Again, by fixing the vector bundle $E\to M$ one gets
\[
\Flat(M,E) \cong \Rep(\Pi_1(M),E)\cong\Hom(\Pi_1(M),\GL(E)),
\]
where $\GL(E)\rightrightarrows M$ is the general linear groupoid of $E$:
\[
\GL(E):=\{(y,A,x):x,y\in M,\, A:E_x\stackrel{\cong}{\to}E_y\}.
\]
For a flat connection on $E$, the corresponding morphism $\tau:\Pi_1(M)\to\GL(E)$ is its parallel transport functor, and from this perspective the holonomy of the connection at a point $x$ is
\[
\Hol_x:=\tau(\Pi_1(M)_x)\subseteq \GL(E_x),
\]
where $\Pi_1(M)_x$ denotes the isotropy at $x$.

In a nutshell, one can understand the holonomy of a flat connection in terms of a three-step process:
\begin{enumerate}
    \item a flat connection on $E$ defines a representation of $\Pi_1(M)$ on $E$, through the Riemann--Hilbert functor,
    \item such a representation defines a parallel transport functor $\tau:\Pi_1(M)\to\GL(E)$,
    \item the holonomy group is the image of the isotropy: $\Hol_x=\tau(\Pi_1(M)_x)$.
\end{enumerate}

Several questions arise: Does a similar process hold for non-flat connections? What about flat higher connections on graded vector bundles (also known as representations up to homotopy or $\infty$-local systems)? Even further, what about (possibly non-flat) higher connections of higher Lie algebroids on graded vector bundles?

One possible answer to the first question is provided by the classical notion of parallel transport, considered as a functor not from $\Pi_1(M)$, but from $\Pi_1^{\mathrm{thin}}(M)$, which is the groupoid of all paths in $M$ up to thin homotopy, an infinite-dimensional groupoid. As for the second question, the first step was taken by Block and Smith \cite{block-smith} by generalizing the Riemann--Hilbert functor. In this paper, we tackle the third question by generalizing this three-step process to higher Lie algebroid connections on graded vector bundles, which recovers the answers to the first two questions.

Our interest in studying higher Lie algebroid connections mainly arises from our interest in Lie 2-algebroid connections (in particular, Courant algebroid connections), whereas the reason for considering graded vector bundles is two-fold. Firstly, it has long been observed that certain Lie algebroid representations are unexpectedly missing. For instance, there is no notion of an adjoint representation for Lie algebroids. To find such an adjoint representation, one has to enlarge the framework to representations up to homotopy (also known as ruths) \cite{ariasabad-crainic-ruths-algebroids}, which are flat higher connections on graded vector bundles. Secondly, connections on graded vector bundles are better suited to seeing the full structure of a higher Lie algebroid, which is itself graded. If one has a higher algebroid $\mathscr{A}\to M$, then an $\Ascr$-connection on a vector bundle $E\to M$ is determined by a map 
\begin{equation*}
    \Gamma(E) \to C^1(\Ascr)\otimes_{\Cinf(M)}\Gamma(E),
\end{equation*}
where $C(\Ascr)$ is the dga of global graded functions on $\Ascr$. In particular, the connection only sees at most the degree 2 part of $C(\Ascr)$, through its curvature. However, if we allow for a graded vector bundle $E = \bigoplus_{i=0}^n E^i$, an $\Ascr$-connection on $E$ is determined by maps
\begin{equation*}
    \Gamma(E^k) \to \bigoplus_{i \geq 0} C^i(\Ascr)\otimes_{\Cinf(M)}\Gamma(E^{k+1-i}),
\end{equation*}
so one can probe the entirety of $\Ascr$ via such connections.

For a higher Lie algebroid $\Ascr\to M$, we implement the aforementioned three-step process as follows:
\begin{enumerate}
    \item we generalize the Riemann--Hilbert functor to non-flat higher Lie algebroid connections, so that an $\Ascr$-connection on a graded vector bundle $E$ defines a connection of the homotopy groupoid $\Pi_\infty(\Ascr)$ on $E$,
    \item we introduce the higher general linear groupoid $\GLi(E)$ of $E$ and prove that the $\Pi_\infty(\Ascr)$-connection on $E$ defines a higher parallel transport functor 
    \[
    \tau:\Pi_\infty(\Ascr)\to\GLi(E),
    \]
    \item we use this functor to define the higher holonomy groups of the connection.
\end{enumerate}
We will now explain in more detail what each one of these steps entails.

Our first main result, regarding the first step, is the extension of the Riemann--Hilbert functor to all connections.

\begin{theorem}[cf. Theorem \ref{thm:riemann-hilbert-higher-algebroids}] \label{thm:intro:riemann-hilbert}
    For each higher Lie algebroid $\Ascr\to M$, there is an $A_\infty$-functor
    \[
    \RH:\HConn(\Ascr)\to \HConn(\Pi_\infty(\Ascr)),
    \]
    which restricts to an $A_\infty$-functor
    \[
    \RH:\HRep(\Ascr)\to\HRep(\Pi_\infty(\Ascr))
    \]
    on the full subcategories of flat objects.
\end{theorem}

Here, $\HConn(\Ascr)$ and $\HConn(\Pi_\infty(\Ascr))$ are the curved dg-categories of $\Ascr$- and $\Pi_\infty(\Ascr)$-connections. An object in $\HConn(\Ascr)$ is a pair $(E,D)$, where $E\to M$ is a graded vector bundle and $D$ is a derivation of the $C(\Ascr)$-module
\[
C(\Ascr)\otimes_{\Cinf(M)}\Gamma(E)
\]
with symbol $Q$ (the cohomological vector field of $\Ascr$), i.e., such that
\[
D(\xi e)=(Q\xi)e+(-1)^{\abs{\xi}}\xi (De)
\]
for all $\xi\in C(\Ascr)$ homogeneous and $e\in C(\Ascr)\otimes_{\Cinf(M)}\Gamma(E)$. The objects of $\HRep(\Ascr)$ are such pairs $(E,D)$ with $D^2=0$. The definitions of $\HConn(\Pi_\infty(\Ascr))$ and $\HRep(\Pi_\infty(\Ascr))$ are similar in spirit, with $C(\Ascr)$ replaced by the dga of cochains $C(\Pi_\infty(\Ascr))$ (for the precise definitions, see Section \ref{sect:higher-groupoids}).

Theorem \ref{thm:intro:riemann-hilbert} can be placed in a tradition of results going back all the way to Chen in the 70s concerning iterated integrals \cite{chen,gugenheim,igusa,arias-abad-schatz,block-smith, antweiler, zeng2025higher}. The existence of this Riemann--Hilbert functor at the level of ruths for $\Ascr=TM$ was established in \cite{block-smith}, where it is furthermore proved that in this case it is actually an $A_\infty$-quasi-equivalence. This result was extended to projectively flat connections with a fixed curvature in \cite{antweiler}. In \cite{arias-abad-schatz}, the Riemann--Hilbert functor for ruths was generalized to the case when $\Ascr$ is an arbitrary Lie algebroid, and in \cite{zeng2025higher} it  was shown that it is an $A_\infty$-quasi-equivalence in the case that $\Ascr=F\subseteq TM$ is a regular foliation. One way to think of this functor at the level of ruths is as a way to integrate ruths of $\Ascr$ to ruths of $\Pi_\infty(\Ascr)$. The question of when such a ruth of $\Pi_\infty(\Ascr)$ descends to a ruth of $\Pi_1(\Ascr)$, when $\Ascr$ is a Lie algebroid, was addressed in \cite{brahic2014obstructions, brahic_ortiz}. Theorem \ref{thm:intro:riemann-hilbert} generalizes these known Riemann--Hilbert functors to non-flat connections of higher Lie algebroids.

The construction of the Riemann--Hilbert functor from Theorem \ref{thm:intro:riemann-hilbert} uses iterated integrals. We approach these integrals in a very concrete fashion via parallel transport. Our approach is a reformulation of the usual approach by Chen. Moreover, we extend the iterated integrals to the theory of curved $A_\infty$-algebras and $A_\infty$-categories, which allows for non-flat $\Ascr$-connections. The key result in this non-flat setting is the non-abelian Stokes's theorem (Theorem \ref{thm:non-abelian-stokes-for-pinched-cubes}). 

For the second step, we introduce the notion of the higher endomorphism groupoid $\Endi(E)$ of a graded vector bundle $E$ (Definition \ref{def:higher_end_groupoid}). This groupoid models the (graded) endomorphisms of $E$ and it has a canonical connection $D_E$ on $E$, which we call the universal connection on $E$ (Definition \ref{def:universal-connection-and-curvature-section}). We further show that there is a canonical simplicial vector bundle $\Curvi(E)\to\Endi(E)$ with a canonical section $K_E$, which we call the curvature section of $\Endi(E)$ (Definition \ref{def:universal-connection-and-curvature-section}), that represents the universal curvature $K_{D_E}$.

The main result in this section is that any $G$-connection on $E$, for any higher Lie groupoid $G$, arises as a pullback of the universal $\Endi(E)$-connection on $E$:

\begin{theorem}[cf. Theorems \ref{thm:simplicialization} and \ref{thm:groupoid-connections-as-functors}] \label{thm:intro:higher-general-linear-groupoid}
Let $E$ be an $n$-term vector bundle. Then, $\Endi(E)$ is a contractible Lie $(n+1)$-groupoid such that for any higher Lie groupoid $G$ over $M$, there is a one-to-one correspondence:
\begin{equation*}
    \left\{\begin{matrix}\text{$G$-connections}\\\text{on $E$}\end{matrix}\right\} \stackrel{1 - 1}{\longleftrightarrow} \left\{\begin{matrix}
        \text{smooth functors $F: G \to \Endi(E)$} \\
        \text{over $M$}
    \end{matrix}\right\}.
\end{equation*}
Furthermore, this correspondence induces an analogous 1-1 correspondence:
\begin{equation*}
    \left\{\begin{matrix}\text{ruths of $G$}\\\text{on $E$}\end{matrix}\right\} \stackrel{1 - 1}{\longleftrightarrow} \left\{\begin{matrix}
        \text{smooth functors $F: G \to \Endi(E)$} \\
        \text{over $M$ such that $K_E \circ F = 0$}
    \end{matrix}\right\}.
\end{equation*}

Explicitly, given a functor $F$, the $G$-connection on $E$ is the pullback along $F$ of the universal connection on $E$, with curvature $F^*K_{D_E}$. 
\end{theorem}

Theorem \ref{thm:intro:higher-general-linear-groupoid} provides an answer to the question posed in \cite{del2019general} of how to represent $G$-ruths on $n$-term vector bundles as a functor from $G$ to some adequate general linear $n$-groupoid of $E$. 

Inside $\Endi(E)$, there is a special subgroupoid that plays a central role in the theory of higher parallel transports: the general linear groupoid $\GLi(E)$ of $E$, defined as the subspace of $\Endi(E)$ where the 1-simplices are isomorphisms (Definition \ref{def:higher-general-linear-and-gauge-groupoid}). By parallel transport $\Pi_\infty(\Ascr)$-connections we mean the $\Pi_\infty(\Ascr)$-connections that arise from $\Ascr$-connections, see Theorem \ref{thm:intro:riemann-hilbert}.

\begin{theorem}[cf. Corollary \ref{cor:higher-parallel-transport} and Proposition \ref{prop:parallel-transport-functor-descends-to-n+1-homotopy-groupoid}] \label{thm:intro:higher-general-linear-groupoid-for-algebroids}
    Let $E$ be an $n$-term vector bundle. Then, $\GLi(E)$ is a contractible Lie $(n+1)$-groupoid such that for any higher Lie algebroid $\Ascr\to M$, there is a 1-1 correspondence:
    \begin{equation*}
        \left\{\begin{matrix}\text{parallel transport}\\\text{$\Pi_\infty(\Ascr)$-connections on $E$}\end{matrix}\right\} \stackrel{1 - 1}{\longleftrightarrow} \left\{\begin{matrix}
            \text{smooth functors $\tau: \Pi_{n+1}(\Ascr)\to\GLi(E)$} \\
            \text{over $M$} 
        \end{matrix}\right\}.
    \end{equation*}
    Furthermore, this correspondence induces an analogous 1-1 correspondence:
    \begin{equation*}
        \left\{\begin{matrix}\text{parallel transport ruths}\\\text{of $\Pi_\infty(\Ascr)$ on $E$}\end{matrix}\right\} \stackrel{1 - 1}{\longleftrightarrow} \left\{\begin{matrix}
            \text{smooth functors $\tau: \Pi_{n}(\Ascr)\to\GLi(E)$} \\
            \text{over $M$ such that $K_E \circ \tau = 0$} 
        \end{matrix}\right\},
    \end{equation*}
     where we consider $\Pi_n(\Ascr)$ and $\Pi_{n+1}(\Ascr)$ as diffeological higher groupoids.
\end{theorem}
Here, $\Pi_n(\Ascr)$ refers to the homotopy $n$-groupoid of $\Ascr$, defined as an appropriate quotient of $\Pi_\infty(\Ascr)$ (Definition \ref{def:homotopy-n-groupoid}).

In general, the functor induced by Theorem \ref{thm:intro:higher-general-linear-groupoid} for the ruth of a higher Lie groupoid $G$ on $E$ lands inside $K_E^{-1}(0)$. This is not a higher Lie groupoid itself. Set theoretically it corresponds to the dg-nerve \cite{kerodon} of a dg-category associated to $E$, and it has the structure of a quasi-smooth simplicial manifold \cite{behrend-liao-xu} (see Remark \ref{rm:quasi-smooth-groupoids}). However, we show that in the case of ruths of transitive higher Lie algebroids, the image of parallel transport naturally lies inside a higher Lie groupoid $\Gaui(E, \delta_x)$ inside $K_E^{-1}(0)$. 

As the final step, we want to use the higher parallel transport functors
\[
\tau:\Pi_{n+1}(\Ascr)\to\GLi(E)
\]
from Theorem \ref{thm:intro:higher-general-linear-groupoid-for-algebroids} to define higher holonomy groups of an $\Ascr$-connection at a point $x\in M$. One could naively hope that $\tau(\Pi_\infty(\Ascr)_x)$ is a valid candidate, but in general this fails to be a higher group. To solve this issue, we change focus and exploit the fact that $\GLi(E)$ is not only a higher groupoid, but also a simplicial groupoid.

\begin{theorem}[cf. Corollary \ref{cor:GLinf(E)-as-simplicial-grpd} and Proposition \ref{prop:canonical-rep-of-GLinf(E)[2]}] \label{thm:intro:GLinf(E)-is-simplicial}
If $E$ is an $n$-term vector bundle, then $\GLi(E)$ can be endowed with the structure of an augmented simplicial Lie groupoid $\GLi(E)[2]\rightrightarrows \GLi(E)[1]$ with a canonical representation.
\end{theorem}

We then use Theorem \ref{thm:intro:GLinf(E)-is-simplicial} to introduce the holonomy groups of an $\Ascr$-connection as augmented simplicial subgroups of the isotropy of $\GLi(E)[2]\rightrightarrows \GLi(E)[1]$, hence inheriting a canonical holonomy representation. Yamabe's theorem on arcwise connected subgroups of Lie groups \cite{yamabe} proves that these groups are smooth, finally completing our three-step program.

\begin{theorem}[cf. Theorem \ref{thm:holonomy-groups-are-smooth}]
The simplicial holonomy groups of an $\Ascr$-connection on $E$ are augmented simplicial Lie subgroups of the isotropy of $\GLi(E)[2]\rightrightarrows \GLi(E)[1]$.
\end{theorem}

In the case that $\Ascr$ is a Lie algebroid and $E$ is a non-graded vector bundle, this notion of holonomy recovers that of \cite{fernandes-poisson,fernandes-lie-algebroids}.

There are several other approaches to higher holonomy that do not fall inside the tradition of Chen and Igusa \cite{other-approaches-1, other-approaches-2, other-approaches-3, other-approaches-4, other-approaches-5, other-approaches-6, other-approaches-7, other-approaches-8, other-approaches-9, other-approaches-10, other-approaches-11, other-approaches-12, other-approaches-13, other-approaches-15, other-approaches-14, other-approaches-16, other-approaches-17}. Most of these develop a notion of higher holonomy in the setting of higher gauge theory. It is unclear to us how precisely our notion of higher holonomy on a graded vector bundle might connect to these notions of higher holonomy for higher gauge theories. This is an interesting direction for future research.

The structure of the paper is as follows. In Sections \ref{sect:higher-groupoids} and \ref{sect:higher-algebroids}, we recall basic notions of higher Lie groupoids and algebroids, and their connections and ruths. In Section \ref{sect:iterated-integrals}, we introduce iterated integrals for endomorphism-valued differential forms, first over pinched cubes and then over simplices. In Section \ref{sect:riemann-hilbert-functor}, we use these iterated integrals to define a Riemann--Hilbert functor for higher Lie algebroid connections. In Section \ref{sect:higher-general-linear-groupoid}, we introduce the higher endomophism groupoid of a graded vector bundle and some interesting subgroupoids (e.g. the higher general linear groupoid) and we prove that morphisms into the endomorphism groupoid model higher Lie groupoid connections and ruths. In Section \ref{sect:higher-holonomy}, we reinterpret the higher general linear groupoid as an augmented simplicial Lie groupoid, we define simplicial holonomy and we show that it is an augmented simplicial Lie subgroup of the isotropy of the higher general linear groupoid.
}
\bigskip

\noindent \textbf{Acknowledgements:} We are grateful to Miquel Cueca for useful discussions, to Ezra Getzler for pointing out relevant references, and to Patrick Antweiler for useful discussions and relevant references and for pointing out the fact that our cube-to-simplex family coincides with his. Finally, we would like to thank Alejandro Cabrera for useful discussions and suggesting the idea of a globular collection in order to simplify the exposition of Section \ref{sect:higher-general-linear-groupoid}. 

\bigskip
\noindent \textbf{AI statement:} No LLMs were used throughout this project.


\section{Higher Lie groupoids and connections} \label{sect:higher-groupoids}

\subsection{Higher Lie groupoids}

We first recall some basic notions from simplicial sets \cite[\href{https://kerodon.net/tag/0001}{Chapter 1}]{kerodon}.

\begin{definition}
The \textbf{simplex category} $\Delta$ has as objects the non-empty linearly ordered sets
\begin{equation*}
    [n] := \{0,\dots, n\}
\end{equation*}
with $n \in \N$ and as morphisms order-preserving maps.
\end{definition}

Recall that the morphisms of $\Delta$ are generated by the \textbf{face} and \textbf{degeneracy maps}. The face map $\delta^i : [n-1] \to [n]$ is the unique injective order preserving map that does not have $i$ in its image, whereas the degeneracy map $\sigma^i : [n + 1] \to [n]$ is the unique surjective order preserving map for which the preimage of $i$ contains 2 elements. Explicitly,
\[
\delta^i(j)=\left\{
\begin{array}{cc}
    j, & \text{if $j< i$,} \\
    j+1, & \text{if $j\geq i$,}
\end{array}
\right. \qquad\text{and}\qquad
\sigma^i(j)=\left\{
\begin{array}{cc}
    j, & \text{if $j\leq i$,} \\
    j-1, & \text{if $j>i$.}
\end{array}
\right.
\]

\begin{definition}
A \textbf{simplicial set} is a functor $X : \Delta^{\mathrm{op}} \to \mathrm{Set}$. Equivalently, it is a collection of sets $\{X_k\}_{k\in\N}$ together with \textbf{face} and \textbf{degeneracy} maps $d_i : X_k \to X_{k-1}$ and $s_i : X_k \to X_{k + 1}$, respectively, for $0\leq i\leq k$, satisfying the \textbf{simplicial relations}:
\begin{align}
    d_i d_j &= d_{j -1} d_i, \quad \text{if $i<j$,} \label{eq:simplicial-relations:faces} \\
    s_i s_j & = s_{j +1} s_i, \quad  \text{if $i \leq j$,} \label{eq:simplicial-relations:degeneracies} \\
    d_i s_j & = \begin{cases}
        s_{j-1} d_{i} & \text{if $i < j$,} \\
        \id, & \text{if $i = j$ or $i = j +1$,} \\
        s_j d_{i-1} &\text{if $i > j + 1$}.
    \end{cases} \label{eq:simplicial-relations:degeneracies-faces}
\end{align}

The \textbf{category $\sSet$ of simplicial sets} has as objects the simplicial sets and as morphisms the natural transformations.
\end{definition}

If $\alpha:[\ell]\to[k]$ is a morphism in $\Delta$, we will denote the induced map in $X$ by $\alpha^*:X_k\to X_\ell$. Whenever we write $\alpha:[\ell]\to[k]$, we will always mean a morphism in $\Delta$ i.e. order-preserving. Sometimes we will also use the relations for $\delta^i$ and $\sigma^i$, which are dual to \eqref{eq:simplicial-relations:faces}--\eqref{eq:simplicial-relations:degeneracies-faces}.

There are some special simplicial sets that play an important role in the theory.

\begin{definition}
The \textbf{$k$-simplex} $\Delta^k\in\sSet$ is defined as the Yoneda embedding of $[k]$
\[
(\Delta^k)_m:=\Hom_{\Delta}([m],[k]).
\]

If $X$ be a simplicial set and $A = (A_k)_{k \in \mathbb{N}}$ a collection of subsets with $A_k \subseteq X_k$, then the \textbf{simplicial set generated by $A$} is the smallest simplicial subset of $X$ containing $A$. Explicitly, the space of its $k$-simplices is
\begin{equation*}
    \{\alpha^*a \in X_k \sep \text{$\alpha:[\ell]\to[k]$ and $a\in A_\ell$, for $\ell\in\N$}\}.
\end{equation*}

The \textbf{boundary of $\Delta^k$}, or the \textbf{$k$-sphere}, is the simplicial subset $\del\Delta^k\subseteq \Delta^k$ of $\Delta^k$ generated by the elements $\delta^i : [k-1] \to [k]$ with $i = 0, \dots, k$. 

The \textbf{$(k,\ell)$-horn of $\Delta^k$}, for $0\leq \ell\leq k$, is the simplicial subset $\Lambda_{k,\ell}\subseteq \del\Delta^k$ generated by the elements $\delta^i$ with $i \neq \ell$.
\end{definition}

By the Yoneda lemma, for any simplicial set $X$, we have that $X_k\cong \Hom_{\sSet}(\Delta^k,X)$ is the set of $k$-simplices. This motivates the following definition.

\begin{definition}
For a simplicial set $X$, we define the set of \textbf{$k$-spheres} and of \textbf{$(k,\ell)$-horns} of $X$, respectively, as
\begin{align*}
\del_k X &:=\Hom_{\sSet}(\del\Delta^k,X), \\
X_{k,\ell} &:= \Hom_{\sSet}(\Lambda_{k,\ell},X).
\end{align*}
\end{definition}

We can naturally identify these sets with 
\begin{equation*} 
    \partial_kX = \left\{(x_0, \dots, x_k) \in \prod_{i = 0}^k X_{k-1} \sep d_ix_j = d_{j-1}x_i \text{ for all } i < j \right\}
\end{equation*}
and
\begin{equation*} 
    X_{k, \ell} = \left\{(x_0, \dots, \widehat{x_\ell}, \dots, x_k) \in \prod_{\substack{i=0 \\ i \neq \ell}}^k X_{k-1} \sep  \begin{matrix} d_ix_j = d_{j-1}x_i \\ \text{ for all } i < j \text{ with } i, j \neq \ell \end{matrix} \right\}.
\end{equation*}

By the simplicial relation \eqref{eq:simplicial-relations:faces}, we may think of the space of $k$-spheres as the natural codomain of the boundary map
\[
\partial_k:= (d_0, \dots, d_k) : X_k \to \prod_{i = 0}^kX_{k-1},
\] 
and the space of $(k,\ell)$-horns as the natural codomain of the map 
\[
d_{k, \ell}:=  (d_0, \dots, \widehat{d_\ell}, \dots, d_k) : X_k \to \prod_{\substack{i=0 \\ i \neq \ell} }^kX_{k-1}.
\]
Given a simplicial set $X$, we say that $x \in X_k$ is \textbf{non-degenerate} if it is not in the image of any degeneracy map. We write $\underline{X}$ for the non-degenerate elements of $X$. For instance, we have that 
\begin{equation*}
    \underline{\Delta^k} = \{\alpha : [i] \to [k]: \alpha \text{ is injective}\}.
\end{equation*}
Furthermore, if $A$ is a simplicial subset of $B$, then $a \in A_i$ is non-degenerate if and only if it is non-degenerate in $B$. Then, since $\alpha \in \Delta^k$ factors through $\delta^i$ if and only if $i \not\in \im \alpha$, we see that
\begin{align*}
    \underline{\partial\Delta^k} & = \{\alpha : [i] \to [k]: \alpha \text{ is injective and $\alpha \neq \id$}\}, \\
    \underline{\Lambda_{k,\ell}} &= \{\alpha : [i] \to [k] : \alpha \text{ is injective, and $\alpha \neq \delta^\ell,\id$}\}.
\end{align*}

\begin{definition}
A simplicial set $X$ is a \textbf{higher groupoid} (or \textbf{$\infty$-groupoid}) if all horns are fillable, i.e., if $d_{k,\ell}:X_k\to X_{k,\ell}$ is surjective for all $k>1$ and $0\leq \ell\leq k$. It is an \textbf{$n$-groupoid} if it is a higher groupoid such that $d_{k,\ell}$ is bijective for all $k>n$.

If the above conditions only hold for $0<\ell<k$, then $X$ is a \textbf{higher category} (or \textbf{$\infty$-category}) or an \textbf{$n$-category}, respectively.
\end{definition}

In the smooth (or diffeological) setting, we will only consider higher groupoids. For more details, see \cite{henriques2008integrating, zhu}.

\begin{definition}
A \textbf{simplicial manifold} is a functor $G : \Delta^{\mathrm{op}} \to \mathrm{Mfd}$. 
A simplicial manifold $G$ is a \textbf{higher Lie groupoid} (or \textbf{Lie $\infty$-groupoid}) if $G_{k,\ell}\subseteq \prod_{i\neq \ell}G_{k-1}$ is an embedded submanifold and $d_{k,\ell}:G_k\to G_{k,\ell}$ is a surjective submersion for all $k>1$ and $0\leq\ell\leq k$. It is a \textbf{Lie $n$-groupoid} if it is a higher Lie groupoid and $d_{k,\ell}$ is a diffeomorphism for all $k>n$.

A \textbf{simplicial diffeological space} is a functor $\Delta^{\mathrm{op}}\to \mathrm{Dflg}$, where $\mathrm{Dflg}$ is the category of diffeological spaces. A \textbf{diffeological higher groupoid} is a simplicial diffeological space whose underlying simplicial set is a higher groupoid.
\end{definition}

\begin{example}
    Given a Lie groupoid $G$, its nerve $\Ner(G)$ defines a Lie $1$-groupoid. Of particular interest to us will be the nerve of the pair groupoid $\Pair(M)$ of a manifold $M$, which we will write as $M^{\bullet + 1}$.
\end{example}

\begin{definition}
    Let $G$ be a higher Lie groupoid. Then, $G$ is called \textbf{contractible} if the boundary maps $\partial_k : G_k \to \partial_kG$ are surjective submersions for all $k \geq 0$.
\end{definition}

This definition makes sense because the sphere spaces $\partial_kG$ for a higher Lie groupoid $G$ are always manifolds \cite[Lem. 2.4]{zhu}.

\begin{example}
    The nerve of the pair groupoid $M^{\bullet + 1}$ is contractible. In this case $\del_k M^{\bullet+1}\cong M^{k+1}$ and $\del_k=\id$.
\end{example}

\subsection{Higher connections and representations up to homotopy}
For more details on representations up to homotopy of Lie groupoids, one can look at \cite{ariasabad-crainic-ruths-groupoids}. For a similar definition on a simplicial set or manifold, see section 2.1 of \cite{arias-abad-schatz}.

The following simplices play a special role in the theory of higher connections and representations up to homotopy.

\begin{definition} \label{def:front-and-back-faces}
The \textbf{front and back $\ell$-faces} on the $k$-simplex $\Delta^k$ are the $\ell$-simplices $\alpha^\ell,\beta^\ell:[\ell]\to[k]$ given, respectively, by
\[
\alpha^\ell(i):=i \quad\text{and}\quad \beta^\ell(i):=k-\ell+i.
\]
The \textbf{$\ell$th vertex} of $\Delta^k$ is the $0$-simplex $\upsilon^\ell:[0]\to[k]$ given by $\upsilon^\ell(0):=\ell$.

If $X$ is a simplicial set, we will write $a_\ell:=(\alpha^\ell)^*:X_k\to X_\ell$, $b_\ell:=(\beta^\ell)^*:X_k\to X_\ell$ and $v_\ell:=(\upsilon^\ell)^*:X_k\to X_0$.
\end{definition}

\begin{remark}
For a simplicial set $X$ we have  that, explicitly,
\begin{align*}
a_\ell &=d_{\ell+1}\dots d_k:X_k\to X_\ell, \\
b_\ell &=d_0\stackrel{(k-\ell)}{\dots}d_0:X_k\to X_\ell, \\
v_\ell &= a_0b_{k-\ell}:X_k\to X_0. \qedhere
\end{align*}
\end{remark}

\begin{definition}
If $G$ is a higher Lie groupoid, we define the \textbf{algebra of cochains of $G$} as the graded algebra $C(G)=\bigoplus_{k\geq 0} C^k(G)$ with $C^k(G):=\Cinf(G_k)$, with product
\[
fg:=(-1)^{k\ell}(a_k^*f)(b_\ell^*g), \quad\text{for $f\in C^k(G)$ and $g\in C^\ell(G)$,}
\]
and differential $d:C^{k-1}(G)\to C^k(G)$ given by $d := (-1)^k\sum_{i=0}^k(-1)^id_i^*$.

We also define the \textbf{subalgebra of normalized cochains of $G$} as $\Chat(G)=\bigoplus_{k\geq 0}\Chat^k(G)$ with
\[
\Chat^k(G):=\{f\in C^k(G): s_i^*f=0,\text{ for $i=0,\dots, k-1$}\}\subseteq C^k(G).
\]
\end{definition}

The sign in the product of $C(G)$ should be understood intuitively as follows: if we evaluate on $x\in G_{k+\ell}$, we have that
\[
(fg)(x)=(-1)^{k\ell}f(a_kx)g(b_\ell x),
\]
and the sign comes from commuting $a_kx$, which is a $k$-simplex, with $g$, which is a degree $\ell$ function. This intuition will guide all the sign conventions used throughout the paper.

It is well known that $C(G)$ is a dga and $\Chat(G)$ a sub-dga of $C(G)$.

\begin{definition}
Let $E=\bigoplus_{k\in\Z}E^k\to G_0$ be a graded vector bundle. We define the \textbf{module of $E$-valued cochains of $G$} as the bigraded vector space $C(G,E)=\bigoplus_{p,q\in\Z}C^p(G,E^q)$ with $C^p(G,E^q):=\Gamma(G_p,v_0^*E^q)$. Elements of $C^p(G,E^q)$ have bidegree $(p,q)$ and total degree $p+q$.

We define the \textbf{submodule of normalized $E$-valued cochains of $G$} as $\Chat(G,E)=\bigoplus_{p,q\in\Z}\Chat^p(G,E^q)\subseteq C(G,E)$ with
\[
\Chat^p(G,E^q):=\{e\in C^p(G,E^q): s_i^*e=0,\text{ for $i=0,\dots, p-1$}\}.
\]
\end{definition}

\begin{lemma}
$C(G,E)$ has a right module structure over the graded algebra structure of $C(G)$. From such structure, $\Chat(G,E)$ inherits a right module structure over the graded algebra structure of $\Chat(G)$.
\end{lemma}

Note that $C(G,E)\cong\Gamma(G_0,E)\otimes_{C^0(G)} C(G)$, so that the right module structure over $C(G,E)$ comes from the algebra structure of $C(G)$. Similarly, $\Chat(G,E)\cong \Gamma(G_0,E)\otimes_{C^0(G)}\Chat(G)$. The product is explicitly given by
\[
\begin{array}{ccc}
C^p(G,E^q)\times C^k(G) & \longrightarrow & C^{p+k}(G,E^q) \\
(e,f) & \longmapsto & (-1)^{pk}(a_p^*e)(b_k^*f).
\end{array}
\]

\begin{definition}
If $E_0,E_1\to G_0$ are two graded vector bundles, we denote by
\[
\Hom_{C(G)}(C(G,E_0),C(G,E_1))
\]
the space of $C(G)$-linear operators $C(G,E_0)\to C(G,E_1)$. It is naturally a bigraded vector space. The space of homogeneous elements of bidegree $(p,q)$ we denote by $\Hom_{C(G)}^{p,q}(C(G,E_0),C(G,E_1))$, whereas the space of homogeneous elements of total degree $k$ we denote by $\Hom_{C(G)}^k(C(G,E_0),C(G,E_1))$. We also use the notation
\[
\End_{C(G)}(C(G,E)):=\Hom_{C(G)}(C(G,E),C(G,E)).
\]

We say that a homomorphism $\Phi\in \Hom_{C(G)}(C(G,E_0),C(G,E_1))$ is \textbf{normalized} if $\Phi(\Chat(G,E_0))\subseteq \Chat(G,E_1)$. We write
\[
\Homhat_{C(G)}(C(G,E_0),C(G,E_1))
\]
for the space of normalized homomorphisms. We also write
\[
\Endhat_{C(G)}(C(G,E)):=\Homhat_{C(G)}(C(G,E),C(G,E)).
\]
\end{definition}

\begin{definition}
A \textbf{$G$-connection on $E$} is a degree $1$ derivation $D$ of $C(G,E)$ with symbol $d$ preserving the normalized cochains, i.e., it is an $\R$-linear degree 1 operator $D:C(G,E)\to C(G,E)$ such that
\begin{enumerate}[label=(\roman*)]
    \item it satisfies the Leibniz rule
    \[
    D(ef) = (De)f + (-1)^{|e|}e(df),
    \]
    for all $e\in C(G,E)$ homogeneous and $f\in C(G)$, and where $\abs{e}$ is the total degree of $e$,
    \item $D(\Chat(G,E))\subseteq \Chat(G,E)$.
\end{enumerate}

The degree $2$ endomorphism $D^2\in \Endhat_{C(G)}^2(C(G,E))$ is the \textbf{curvature} of $D$. We say that $D$ is \textbf{flat} or a \textbf{representation up to homotopy} (or \textbf{ruth} for short) if $D^2=0$.
\end{definition}

\begin{remark}
A ruth on $E$ makes $C(G,E)$ into a right dg-module over $C(G)$.
\end{remark}

Recall that if $E_0,E_1\to M$ are vector bundles, then the vector bundle $\Hom(E_0,E_1)\to M$ encodes the morphisms between $E_0$ and $E_1$, i.e.,
\[
\Hom_{\Cinf(M)}(\Gamma(E_0),\Gamma(E_1))\cong\Gamma(\Hom(E_0,E_1)).
\]
We have a similar phenomenon for vector-valued cochains on a higher Lie groupoid $G$, which allows to rewrite $G$-connections and their curvatures. This is well known, the case when $G$ is the nerve of a Lie groupoid having been worked out for the first time in \cite{ariasabad-crainic-ruths-groupoids}. We will not provide proofs, since they follow along the same lines as in \cite{ariasabad-crainic-ruths-groupoids}.

\begin{definition}
We define the bigraded vector space 
\[
C_G(\Hom(E_0,E_1))=\bigoplus_{p,q\in\Z} C_G^p(\Hom^q(E_0,E_1))
\]
with $C_G^p(\Hom^q(E_0,E_1)):=\Gamma(G_p, \Hom^q(v_p^*E_0,v_0^*E_1))$. We also define the subspace of normalized elements
\[
\Chat_G(\Hom(E_0,E_1))=\bigoplus_{p,q\in\Z} \Chat_G^p(\Hom^q(E_0,E_1))
\]
with
\[
\Chat^p_G(\Hom^q(E_0,E_1)):=\{\phi\in C^p_G(\Hom^q(E_0,E_1)): s_i^*\phi=0,\text{ for $i=0,\dots, p-1$}\}.
\]
We will also use the notation $C_G(\End(E)):=C_G(\Hom(E,E))$.
\end{definition}

\begin{lemma} \label{lem:hom-C(G,E)}
As bigraded vector spaces, 
\[
\Hom_{C(G)}(C(G,E_0),C(G,E_1))\cong C_G(\Hom(E_0,E_1)),
\]
and this isomorphism induces an isomorphism
\[
\Homhat_{C(G)}(C(G,E_0),C(G,E_1))\cong \Chat_G(\Hom(E_0,E_1)).
\]

Moreover, under this identification the evaluation map is given by
\[
\begin{array}{ccc}
C^p_G(\Hom^q(E_0,E_1))\times C^k(G,E_0^\ell) & \longrightarrow & C^{p+k}(G,E_1^{q+\ell}) \\
(\phi,e) & \longmapsto & (-1)^{p(k+\ell)}(a_p^*\phi)(b_k^*e)
\end{array}
\]
and the composition map by
\[
\begin{array}{ccc}
C^p_G(\Hom^q(E_1,E_2))\times C^k_G(\Hom^\ell(E_0,E_1)) & \longrightarrow & C^{p+k}_G(\Hom^{q+\ell}(E_0,E_2)) \\
(\phi,\psi) & \longmapsto & (-1)^{p(k+\ell)}(a_p^*\phi)(b_k^*\psi).
\end{array}
\]
\end{lemma}

\begin{proposition} \label{prop:connections-simplicial-description}
There is a canonical flat $G$-connection $D_0$ on $E$, given explicitly on $e\in C^p(G,E^q)$ by
\[
D_0e=(-1)^{p+q}\sum_{i=1}^{p+1}(-1)^id_i^*e.
\]

In particular, any $G$-connection on $E$ is equivalent to a degree 1 cochain $F\in C_G(\End(E))^1$ such that, if $F_i\in C^i_G(\End^{1-i}(E))$ denotes the components of $F$, then $F_i$ is normalized for all $i\neq 1$ and $F_1$ is such that $s_0^*F_1=\id$. 

Furthermore, the curvature of $F$ is a degree 2 normalized cochain $K\in \Chat_G(\End(E))^2$ whose components $K_i\in C^i_G(\End^{2-i}(E))$ are given by
\begin{equation} \label{eq:curvature_ruth}
    K_i=-\sum_{j=1}^{i-1}(-1)^jd_j^*F_{i-1} +\sum_{j= 0}^i (-1)^j(a_j^*F_j)(b_{i-j}^*F_{i-j}).
\end{equation}
\end{proposition}

\bigskip
\section{Higher Lie algebroids and connections} \label{sect:higher-algebroids}

\subsection{Higher Lie algebroids}

For the basic theory of graded manifolds and higher Lie algebroids presented here, see for instance \cite{fairon,jotz-mehta-papantonis}.

\begin{definition}
A \textbf{graded manifold} is a pair $\Mscr=(M,C_{\Mscrh})$, where $M$ is a smooth manifold and $C_{\Mscrh}$ is a sheaf of real graded commutative algebras which is locally isomorphic to $\Cinf_M\otimes_\R SV^*$, for some $\N$-graded vector space $V$ with $V_0=\R$, by which we mean that for every $x\in M$ there is an open neighborhood $U$ of $x$ such that
\[
C_{\Mscrh}|_U\cong (\Cinf_M\otimes_\R SV^*)|_U.
\]
Here, $SV^*$ is the graded symmetric algebra of $V^*$. The \textbf{algebra of functions} of $\Mscr$ we define to be $C(\Mscr):=C_{\Mscr}(M)$.

A \textbf{morphism} between graded manifolds $\Mscr=(M,C_{\Mscrh})$ and $\Nscr=(N,C_{\Nscrh})$ is a pair $\phi=(\phi^*,\phi_0)$, where $\phi_0:M\to N$ is a smooth map and $\phi^*:C_{\Nscrh}\to\phi_{0*}C_{\Mscrh}$ is a morphism of sheaves of graded commutative algebras. Recall that the pushforward sheaf $\phi_{0*}C_{\Mscrh}$ over $N$ is defined by
\[
\phi_{0*}C_{\Mscrh} (U):=C_{\Mscrh}(\phi_0^{-1}(U)) \qquad \text{for $U\subseteq N$ open.}
\]

A \textbf{degree $k$ derivation} of $\Mscr$, or \textbf{degree $k$ vector field}, is a degree $k$ morphism of sheaves of real graded vector spaces $X:C_{\Mscrh}\to C_{\Mscrh}$ satisfying the Leibniz rule
\[
X(\xi\eta)=(X\xi)\eta+(-1)^{k\abs{\xi}}\xi(X\eta),
\]
for $U\subseteq M$ open and $\xi,\eta\in C_{\Mscrh}(U)$. In such case we write $\abs{X}:=k$. The space of degree $k$ derivations we write as $\Der^k(\Mscr)$. The space of derivations is the graded vector space
\[
\Der(\Mscr):=\bigoplus_{k\in\Z}\Der^k(\Mscr).
\]
This space comes with a graded Lie bracket given by the graded commutator of derivations:
\[
[X,Y]:=XY-(-1)^{\abs{X}\abs{Y}}YX.
\]

A \textbf{higher Lie algebroid} or \textbf{Lie $\infty$-algebroid} (or differential graded manifold, or dg-manifold) is a pair $(\Ascr,Q)$, where $\Ascr$ is a graded manifold and $Q$ is a degree 1 vector field on $\Ascr$ such that $Q^2=0$. In such case we call $Q$ the \textbf{differential} (or cohomological vector field) of $(\Ascr,Q)$. If $C_{\Ascr}$ is generated by the degrees $0,\dots,n$ we call $\Ascr$ a \textbf{Lie $n$-algebroid}.

A \textbf{morphism of higher Lie algebroids} $\phi:(\Ascr, Q)\to(\Bscr,Q)$ is a morphism of graded manifolds $\phi:\Ascr\to\Bscr$ such that $\phi^*Q=Q\phi^*$.
\end{definition}

\begin{remark}
    What we called graded manifolds are sometimes called non-negatively graded manifolds. All our graded manifolds will be non-negatively graded.
\end{remark}

\begin{example}
Let $A\to M$ be a Lie algebroid. Then $(M,\Omega_A,d_A)$ is a Lie 1-algebroid, where $\Omega_A(U):=\Gamma(U,\wedge A^*)$. This provides an equivalence of categories between Lie algebroids and Lie 1-algebroids.
\end{example}

\begin{example}
Let $A=\bigoplus_{i=0}^\infty A_i\to M$ be a graded vector bundle. Then $(M,\Gamma_{SA^*})$ is a graded manifold, where $\Gamma_{SA^*}(U):=\Gamma(U,SA^*)$. Graded manifolds of this type are known as split graded manifolds. It is well known that any graded manifold is non-canonically isomorphic to a split one, a result sometimes known as Batchelor's theorem (in the case of $\Z_2$-graded manifolds) \cite{batchelor}. If $\Mscr$ is a graded manifold, an isomorphism $\Mscr\cong (M,\Gamma_{SA^*})$ is called a splitting of $\Mscr$.
\end{example}

\begin{definition}
Let $\Mscr$ be a graded manifold. A \textbf{graded vector bundle over $\Mscr$}, which we write as $\Escr\to\Mscr$, is a locally freely and finitely generated sheaf $\Escr$ of $\N$-graded (left) $C_{\Mscrh}$-modules. A graded vector bundle $\Escr$ is called \textbf{$n$-term} if $\Escr$ is generated by degrees $0,\dots, n-1$.
\end{definition}

\begin{example}
Let $\Mscr$ be a graded manifold, and consider the sheaf $T\Mscrh$ given by
\[
T\Mscrh (U):=\Der(\Mscr|_U),
\]
where $\Mscr|_U:=(U,C_{\Mscrh}|_U)$. Then $T\Mscrh\to\Mscr$ is a graded vector bundle, the \textbf{graded tangent bundle of $\Mscr$}, whose global sections are precisely vector fields on $\Mscr$.
\end{example}

We are also interested in vector bundle maps for graded vector bundles.

\begin{definition}
Let $\Escr,\Fscr\to\Mscr$ be graded vector bundles. A \textbf{degree $k$ morphism} $\Phi:\Escr\to\Fscr$ over the identity is a degree $k$ morphism of sheaves of $C_{\Mscrh}$-modules, i.e.,
\[
\Phi(\xi s)=(-1)^{k\abs{\xi}}\xi\Phi(s),
\]
for $U\subseteq M$ open, $\xi\in C_{\Mscrh}(U)$ and $s\in\Escr(U)$.  In such case we write $\abs{\Phi}:=k$. The space of degree $k$ morphisms between $\Escr$ and $\Fscr$ over the identity we denote by $\Hom^k(\Escr,\Fscr)$. The space of morphisms between $\Escr$ and $\Fscr$ is the graded vector space
\[
\Hom(\Escr,\Fscr):=\bigoplus_{k\in\Z} \Hom^k(\Escr,\Fscr).
\]
We also define the graded vector bundle of morphisms $\Homscr(\Escr,\Fscr)\to\Mscr$ to be
\[
\Homscr(\Escr,\Fscr)(U):=\Hom(\Escr|_U,\Fscr|_U), \qquad\text{for $U\subseteq M$ open.}
\]

We define the \textbf{dual} of $\Escr$ as $\Escr^*:=\Homscr(\Escr,\underline\R)$, where $\underline\R\to\Mscr$ is $C_{\Mscrh}$ regarded as a $C_{\Mscrh}$-module.

If $\Fscr$ has a different base, say $\Fscr\to\Nscr$, then a degree $k$ morphism $\Phi:\Escr\to\Fscr$ over $\phi:\Mscr\to\Nscr$ is a degree $k$ morphism of sheaves of modules $\Phi^*:\Fscr^*\to\phi_{0*}\Escr^*$ over $\phi^*:C_{\Nscrh}\to\phi_{0*}C_{\Mscrh}$, i.e.,
\[
\Phi^*(\xi\lambda)=(-1)^{k\abs{\xi}}(\phi^*\xi)(\Phi^*\lambda),
\]
for $U\subseteq N$ open, $\xi\in C_{\Nscrh}(U)$ and $\lambda\in\Fscr^*(U)$. In such case we write $\abs{\Phi}:=k$.
\end{definition}

\begin{example}
If $\Mscr$ is a graded manifold, then the \textbf{graded cotangent bundle of $\Mscr$} is the dual of $T\Mscrh$, which we write as $T^*\Mscr$, and whose global sections are referred to as differential 1-forms.
\end{example}

\subsection{Higher connections and representations up to homotopy}
The following material is standard. A good reference in the case of Lie algebroids is \cite{ariasabad-crainic-ruths-algebroids} and in the case of higher Lie algebroids, see \cite{jotz-mehta-papantonis}.

\begin{definition} \label{def:higher-algebroid-connections}
Let $(\Ascr\to M,Q)$ be a higher Lie algebroid and $E\to M$ a graded vector bundle. Let $\Escr:=C_{\Ascr}\otimes_{\Cinf_M}\Gamma_E\to\Ascr$ be the induced graded vector bundle over $\Ascr$, where $\Gamma_E$ is the sheaf of sections of $E$. An \textbf{$\Ascr$-connection} on $E$ is degree 1 derivation $D$ of $\Escr$ with symbol $Q$, i.e., a degree 1 morphism of sheaves of real graded vector spaces $D:\Escr\to\Escr$ satisfying the Leibniz rule
\[
D(\xi s)=(Q\xi) s+(-1)^{\abs{\xi}}\xi D s,
\]
for $U\subseteq M$ open, $\xi\in C_{\Ascr}(U)$ and $s\in\Escr(U)$. The degree 2 morphism of sheaves of real graded $C_{\Ascr}$-modules $D^2:\Escr\to\Escr$ we call the \textbf{curvature} of $D$. We say that $D$ is \textbf{flat} or a \textbf{representation up to homotopy} (or \textbf{ruth} for short) if $D^2=0$. 
\end{definition}

\begin{example}
If $A\to M$ is a Lie algebroid and $E\to M$ is a vector bundle, then an $A$-connection in the usual sense, meaning a real operator $\nabla:\Gamma(E)\to\Omega^1(A,E)$ satisfying the Leibniz rule
\[
\nabla(fs)=d_Af\otimes s+f\nabla s,
\]
for $f\in\Cinf(M)$ and $s\in\Gamma(E)$, is an $A$-connection in the sense of Definition \ref{def:higher-algebroid-connections} by considering the induced derivation $d^\nabla$ on $\Omega(A,E)$, whose curvature is $K^\nabla\wedge(-)$, where $K^\nabla\in\Omega^2(A,\End(E))$ is the usual curvature of $\nabla$. In the case of Lie algebroids, the term $A$-connection we will reserve for the usual meaning of the word, whereas connections in the sense of Definition \ref{def:higher-algebroid-connections} we will call higher $A$-connections.
\end{example}

Notice that if $E_i\to M$ are graded vector bundles, then
\[
C_{\Ascr}\otimes_{\Cinf_M}\Gamma_{\Hom(E_0,E_1)}\cong\Homscr(\Escr_0,\Escr_1).
\]
In particular, $C_{\Ascr}\otimes_{\Cinf_M}\Gamma_{E^*}\cong\Escr^*$.

\begin{lemma} \label{lem:higher-connection-on-homomorphisms}
Let $(E_i\to M,D_i)$ be graded vector bundles with $\Ascr$-connections for $i=0,1$, where $(\Ascr\to M,Q)$ is a higher Lie algebroid. Then there is an induced higher connection on $\Hom(E_0,E_1)\to M$ given by
\[
(D\Phi)s:=D_1(\Phi s)-(-1)^{\abs{\Phi}}\Phi(D_0s),
\]
for $U\subseteq M$ open, $\Phi\in\Homscr(\Escr_0,\Escr_1)(U)$ and $s\in \Escr_0(U)$.
\end{lemma}

\begin{corollary}
If $(E\to M,D)$ is a graded vector bundle with $\Ascr$-connection, where $(\Ascr\to M,Q)$ is a higher Lie algebroid, then there is an induced $\Ascr$-connection on $E^*$.
\end{corollary}

\begin{example}
Let $(\Ascr\to M,Q)$ be a higher Lie algebroid. The \textbf{adjoint representation of $\Ascr$} is $(T\Ascr,\ad(Q))$, where
\[
\ad(Q)X:=[Q,X]=QX-(-1)^{\abs{X}}XQ.
\]
Strictly speaking, in general it is not canonically a ruth, in the sense that $T\Ascr$ is not canonically of the form $C_{\Ascr}\otimes_{\Cinf_M}\Gamma_E$ for a graded vector bundle $E\to M$. The \textbf{coadjoint representation} of $\Ascr$ is the dual representation $(T^*\Ascr,\ad^*(Q))$.
\end{example}

It is well known that for Lie algebroids connections can be described in terms of geometric data. A proof of this fact can be found in the original paper \cite{ariasabad-crainic-ruths-algebroids}.

\begin{proposition} \label{prop:ruths-for-lie-algebroids}
Let $A\to M$ be a Lie algebroid and $E\to M$ a graded vector bundle. An $A$-connection $D$ on $E$ is equivalently given by a triple $(\delta,\nabla,\omega)$, where
\begin{enumerate}
	\item $\delta\in\Gamma(\End^1(E))$,
	\item $\nabla$ is a degree-preserving $A$-connection on $E$,
	\item $\omega=\sum_{i\geq 2}\omega_i$, with $\omega_i\in\Omega^i(A,\End^{1-i}(E))$.
\end{enumerate}
Its curvature is represented by $K=\sum_{i\geq 0}K_i$, with $K_i\in\Omega^i(A,\End^{2-i}(E))$ given by
\begin{enumerate}
	\item $K_0:=\delta^2\in\Gamma(\End^2(E))$,
	\item $K_1:=\nabla\delta\in\Omega^1(A,\End^1(E))$,
	\item $K_2:=K^\nabla+\omega_2\wedge\delta+\delta\wedge\omega_2\in\Omega^2(A,\End^0(E))$,
	\item $K_i:=d^\nabla\omega_{i-1}+\omega_i\wedge\delta+\delta\wedge\omega_i + \sum_j\omega_j\wedge\omega_{i-j}\in\Omega^i(A,\End^{2-i}(E))$, for $i>2$.
\end{enumerate}
\end{proposition}

Bear in mind that the wedge product in Proposition \ref{prop:ruths-for-lie-algebroids} is the one from Example \ref{ex:end-valued-forms-dga}.

\bigskip
\section{Iterated integrals} \label{sect:iterated-integrals}

\subsection{Curved differential graded algebras} \label{sect:curved-dg-algebras}
We will need the following standard material about curved differential graded algebras and $A_\infty$-maps between them. See for example section 1 of \cite{getzler-jones} or section 3 and 7 of \cite{positselski2011two}.
\begin{definition}
A \textbf{curved differential graded algebra} (curved dga for short) is a tuple $(A,\wedge,d,K)$, where $(A,\wedge)$ is a graded associative real algebra, $d$ is a degree 1 derivation, i.e., such that for all $a,b\in A$ homogeneous,
\[
d(a\wedge b)=da\wedge b+(-1)^{\abs{a}}a\wedge db,
\]
and $K\in A^2$ is such that $d^2=[K,\cdot]$ (where $[\cdot,\cdot]$ is the graded commutator) and $dK=0$. The element $K$ is called the \textbf{curvature}.
\end{definition}

\begin{example} \label{ex:end-valued-forms-dga}
Let $A\to M$ be a Lie algebroid and $(E\to M, \nabla)$ a graded vector bundle with a degree-preserving $A$-connection. Then $\Omega(A,\End(E))$ becomes a curved dga with its total grading, with wedge product given on decomposables by
\[
(\xi\otimes\Phi)\wedge (\eta\otimes\Psi):= (-1)^{\abs{\Phi}\abs{\eta}}(\xi\wedge\eta)\otimes(\Phi\Psi)
\]
and with differential $d^\nabla$ and curvature $K^\nabla$.
\end{example}

The theory of curved dgas is part of the wider theory of curved $A_\infty$-algebras (see, e.g., \cite{getzler-jones}).

Recall that if $A$ is a graded vector space, then $A[n]$ denotes the $n$-shifted vector space, given by $A[n]^k:=A^{n+k}$.

\begin{definition}
Let $A$ be a graded vector space. Consider the tensor coalgebra
\[
\B A:=T(A[1])=\bigoplus_{\ell\geq 0} A[1]^{\otimes \ell}
\]
with coproduct
\[
\Delta(a_1\otimes\dots\otimes a_\ell)=\sum_{i=0}^\ell (a_1\otimes\dots\otimes a_i)\otimes (a_{i+1}\otimes \dots\otimes a_\ell).
\]
An \textbf{$A_\infty$-structure on $A$} is a codifferential on $\B A$, i.e., a degree 1 linear map $D:\B A\to \B A$ such that
\[
\Delta D=(\id\otimes D+ D\otimes\id)\Delta
\]
and $D^2=0$.

If $(A_i,D_i)$ are $A_\infty$-algebras, for $i=0,1$, then a morphism of $A_\infty$-algebras $(A_0,D_0)\to (A_1,D_1)$ is a degree 0 map of coalgebras $\Psi:\B A_0\to \B A_1$ 
such that $\Psi D_0=D_1\Psi$.
\end{definition}

\begin{remark}
What we called here $A_\infty$-algebra is often called curved $A_\infty$-algebra.
\end{remark}

Codifferentials on $\B A$ have the following explicit description.

\begin{proposition}
An $A_\infty$-structure on $A$ is equivalently given by a collection of degree 1 maps $b_\ell:A[1]^{\otimes \ell}\to A[1]$, for $\ell\geq 0$, such that
\[
\sum_{i+j+k=\ell} b_{i+k+1}\circ (\id^{\otimes i}\otimes b_j\otimes \id^{\otimes k})=0,
\]
and a morphism $(A_0,D_0)\to (A_1,D_1)$ is equivalently given by a collection of degree 0 maps $\psi_\ell:A_0[1]^{\otimes \ell}\to A_1[1]$, for $\ell >0$, such that
\begin{equation} \label{eq:broken_down_eq}
\sum_{k\geq 0}\sum_{i_1+\dots+i_k=\ell} b^1_k \circ (\psi_{i_1}\otimes\dots\otimes \psi_{i_k}) =\sum_{i+j+k=\ell} \psi_{i+k+1}\circ (\id^{\otimes i}\otimes b^0_j\otimes \id^{\otimes k}).
\end{equation}
The identification is provided by setting
\begin{align*}
D &=\sum_{\ell\geq 0} \sum_{i+j+k=\ell} \id^{\otimes i}\otimes b_j\otimes \id^{\otimes k}, \\
\Psi &= \sum_{\ell> 0} \sum_{i_1+\dots+i_k=\ell}\psi_{i_1}\otimes \dots \otimes \psi_{i_k}.
\end{align*}

Under such identification, the composition of two morphisms $\Psi_0$ and $\Psi_1$ is given by the maps
\[
(\psi^1\psi^0)_\ell=\sum_{i_1+\dots +i_k=\ell} \psi^1_k\circ(\psi^0_{i_1}\otimes\dots\otimes \psi^0_{i_k}).
\]
\end{proposition}

\begin{remark}
Notice that a priori no $\psi_0$ term is allowed in a morphism of $A_\infty$-algebras, as this would lead to convergence issues in the above equation for $\Psi$ \cite{getzler-jones}. However, in case the $A_\infty$-algebras involved only have a finite number of non-zero $b_\ell$ operators, then equation \eqref{eq:broken_down_eq} makes sense even with a $\psi_0$-term. In the present paper, we will always only have a $b_0$, $b_1$ and $b_2$. So, we always allow for a $\psi_0$ term. 
\end{remark}

For an $A_\infty$-algebra $(A,b)$, the element $b_0(1)\in A[1]^1$ is called the \textbf{curvature} of $A$.

It is also usual to present $A_\infty$-algebras using $A$ itself instead of $A[1]$, the difference resulting in a change of signs. If $s:A\to A[1]$ is the \textbf{suspension map}, defined by $sa:=a$, which is of degree $-1$, then we can transform each $b_\ell$ into a degree $2-\ell$ map $m_\ell:A^{\otimes \ell}\to A$ by asking for the following diagram to commute:
\[
\begin{tikzcd}
	{A[1]^{\otimes \ell}} & {A[1]} \\
	{A^{\otimes \ell}} & A.
	\arrow["{b_\ell}", from=1-1, to=1-2]
	\arrow["{(s^{-1})^{\otimes \ell}}"', from=1-1, to=2-1]
	\arrow["{m_\ell}"', from=2-1, to=2-2]
	\arrow["s"', from=2-2, to=1-2]
\end{tikzcd}
\]
The compatibility equations for $m$ are the same as for $b$ with additional signs.

The result that we are interested in is the following simple application:

\begin{corollary} \label{cor:bar-complex}
A curved dga $(A,\wedge,d,K)$ is the same as an $A_\infty$-algebra $(A,m)$ with $m_0(1)=K$, $m_1=d$, $m_2=\wedge$ and $m_n=0$ for $n\geq 3$. In particular, every curved dga carries a codifferential on its tensor coalgebra.
\end{corollary}

We will need the explicit description of the codifferential from Corollary \ref{cor:bar-complex} for our later computations. If $(A,\wedge, d,K)$ is a curved dga, then the codifferential $D$ acts as follows:
\begin{equation} \label{eq:D-bar-complex}
\begin{split}
D & (sa_1\otimes \dots\otimes sa_\ell)= \\
&= \sum_{i=1}^\ell (-1)^{ \sum_{k=1}^{i-1}(\abs{a_k}-1)} sa_1\otimes\dots\otimes sa_{i-1}\otimes sda_i\otimes sa_{i+1}\otimes \dots\otimes sa_\ell \\
&\qquad +\sum_{i=1}^{\ell-1} (-1)^{\sum_{k=1}^i (\abs{a_k}-1)} sa_1\otimes\dots\otimes sa_{i-1}\otimes s(a_i\wedge a_{i+1})\otimes sa_{i+2}\otimes\dots\otimes sa_\ell \\
&\qquad +\sum_{i=1}^{\ell+1}(-1)^{\sum_{k=1}^{i-1}(\abs{a_k}-1)} sa_1\otimes\dots\otimes sa_{i-1}\otimes sK \otimes sa_i\otimes\dots\otimes sa_\ell. 
\end{split}
\end{equation}

Regarding morphisms, we will use throughout the following notation: if $\psi:A_0\to A_1$ is a map of $A_\infty$-algebras, then we will write
\[
\psi_\ell(sa_1\otimes\dots\otimes sa_\ell)=s\hat\psi_\ell(a_1\otimes\dots\otimes a_\ell),
\]
for $\hat\psi_\ell:A_0^{\otimes\ell}\to A_1$ of degree $1-\ell$. Using such notation, we get the following explicit description for $A_\infty$-maps between curved dgas, which will be useful in the following sections.

\begin{corollary} \label{cor:A-infty-map-for-curved-dgas}
An $A_\infty$-map between two curved dgas $(A_i,\wedge, d_i, K_i)$, for $i=0,1$, is the same as a family of maps $\hat\psi_\ell:A_0^{\otimes\ell}\to A_1$ of degree $1-\ell$, for $\ell\geq 0$ such that
\begin{align*}
d_1(\hat\psi_\ell & (a_1\otimes\dots\otimes a_\ell))+\sum_{i+j=\ell}(-1)^{\sum_{k=1}^i(\abs{a_k}-1)} \hat\psi_i(a_1\otimes\dots\otimes a_i)\wedge \hat\psi_j(a_{i+1}\otimes\dots\otimes a_\ell)= \\
&=\sum_{i=1}^\ell (-1)^{ \sum_{k=1}^{i-1}(\abs{a_k}-1)} \hat\psi_\ell(a_1\otimes\dots\otimes a_{i-1}\otimes d_0a_i\otimes a_{i+1}\otimes\dots\otimes a_\ell) \\
&\qquad +\sum_{i=1}^{\ell-1} (-1)^{\sum_{k=1}^i (\abs{a_k}-1)} \hat\psi_{\ell-1}(a_1\otimes\dots\otimes a_{i-1}\otimes(a_i\wedge a_{i+1})\otimes a_{i+2}\otimes\dots\otimes a_\ell) \\
&\qquad +\sum_{i=1}^{\ell+1}(-1)^{\sum_{k=1}^{i-1}(\abs{a_k}-1)} \hat\psi_{\ell+1}(a_1\otimes\dots\otimes a_{i-1}\otimes K_0\otimes a_i\otimes \dots\otimes a_\ell).
\end{align*}
The equation for $\ell=0$ reads
\[
K_1+d_1 \left(\hat\psi_0(1)\right)-\hat\psi_0(1)\wedge\hat\psi_0(1)=\hat\psi_1(K_0).
\]
\end{corollary}

\begin{proof}
The equations for $\psi:A_0\to A_1$ being an $A_\infty$-map when $b^i_\ell=0$ for $\ell\geq 3$ are
\begin{align*}
b_1^1\circ\psi_\ell & +\sum_{i+j=\ell} b_2^1\circ(\psi_i\otimes\psi_j)=\sum_{i+j=\ell}\psi_{\ell+1}\circ (\id^{\otimes i}\otimes b_0^0\otimes\id^{\otimes j})\\
&\qquad +\sum_{i+j=\ell-1}\psi_\ell\circ(\id^{\otimes i}\otimes b_1^0\otimes\id^{\otimes j}) +\sum_{i+j=\ell-2} \psi_{\ell-1}\circ (\id^{\otimes i}\otimes b_2^0\otimes\id^{\otimes j})
\end{align*}
for $\ell>0$, and
\[
b_0^1(1)+b_1^1\circ\psi_0+b_2^1\circ(\psi_0\otimes\psi_0)=\psi_1\circ b_0^0. \qedhere
\]
\end{proof}

Lastly, we recall the so called Maurer--Cartan elements of $A_\infty$-algebras.

\begin{definition}
Let $(A,b)$ be an $A_\infty$-algebra and $x\in A[1]$. The \textbf{curvature of $x$} we define to be the element
\[
K_x:=\sum_{\ell\geq 0} b_\ell(x^{\otimes \ell})\in A[1].
\]

Let $\psi:A_0\to A_1$ be a map of $A_\infty$-algebras and $x\in A_0[1]$. The \textbf{image of $x$ by $\psi$} we define to be the element
\[
\psi(x):=\sum_{\ell\geq 0} \psi_\ell(x^{\otimes \ell})\in A_1[1].
\]
\end{definition}

Typically an element $x\in A[1]^0$ such that $K_x=0$ is called a Maurer--Cartan element, although we will not use this concept.

To see how $K$ interacts with $\psi$, we need to use the \textbf{sign map} $\epsilon:A\to A$, defined on homogeneous elements by $\epsilon(x):=(-1)^{\abs{x}}x$.

\begin{proposition} \label{prop:curvature-in-A-infty-algebras}
Let $\psi:A_0\to A_1$ be a map of $A_\infty$-algebras and $x\in A_0[1]$. Then
\[
K_{\psi(x)}=\sum_{i,j\geq 0}\psi_{i+j+1}(\epsilon(x)^{\otimes i}\otimes K_x\otimes x^{\otimes j}).
\]
\end{proposition}

\begin{proof}
Direct computation:
\begin{align*}
K_{\psi(x)} &=\sum_{\ell\geq 0} b_\ell^1\left(\left(\sum_{k\geq 0} \psi_k(x^{\otimes k})\right)^{\otimes \ell}\right) \\
&=\sum_{\ell\geq 0} \sum_{k\geq 0}\sum_{k_1+\dots +k_\ell=k} b_\ell^1 \circ (\psi_{k_1}\otimes\dots \otimes \psi_{k_\ell}) (x^{\otimes k}) \\
&=\sum_{k\geq 0} \sum_{i+p+j=k} \psi_{i+j+1}\circ(\id^{\otimes i}\otimes b_p^0\otimes\id^{\otimes j}) (x^{\otimes k}) \\
&=\sum_{i,j\geq 0} \psi_{i+j+1}\left(\epsilon(x)^{\otimes i} \otimes \sum_{p\geq 0} b_p^0(x^{\otimes p}) \otimes x^{\otimes j} \right) \\
&=\sum_{i,j\geq 0} \psi_{i+j+1}(\epsilon(x)^{\otimes i}\otimes K_x\otimes x^{\otimes j}). \qedhere
\end{align*}
\end{proof}

We can also translate these concepts to the $m$ and $\hat\psi$ notation. For $a\in A$, we set
\[
\hat\psi (a):=\sum_{\ell\geq 0}\hat\psi_\ell(a^{\otimes\ell}),
\]
so that $\psi(sa)=s\hat\psi(a)$. In line with this, we also set
\[
\hat K_a:=\sum_{\ell\geq 0} (-1)^{\ell(\ell-1)/2}m_\ell\left(\bigotimes_{i=1}^\ell \epsilon^{\ell-i}(a)\right),
\]
so that $K_{sa}=s\hat K_a$. We can then rewrite Proposition \ref{prop:curvature-in-A-infty-algebras} as
\begin{equation} \label{eq:curvature-in-A-infty-algebras-hat-notation}
\hat K_{\hat\psi(a)}=\sum_{i,j\geq 0}(-1)^i\hat\psi_{i+j+1}(\epsilon(a)^{\otimes i}\otimes\hat K_a\otimes a^{\otimes j}).
\end{equation}

\subsection{Iterated integrals for cubes} \label{sect:iterated-integrals-cubes}

Throughout this section $A\to M$ will always be a Lie algebroid and $(E\to M,\nabla)$ a graded vector bundle with a degree-preserving $A$-connection.

Let $I=[0,1]$ be the unit interval. For the cube $I^k$ we will typically use coordinates $(x_0,x)$, with $x=(x_1,\dots, x_{k-1})\in I^{k-1}$. We will denote the coordinate vector fields on $I^k$ by $\del_i$, for $i=0,\dots, k$.

\begin{definition}
A \textbf{pinched $k$-cube} in $A$ is a Lie algebroid map $\phi:TI^k\to A$ covering $\phi_0:I^k\to M$ such that
\[
\phi|_{(m,x)}=0, \qquad x\in I^{k-1},\, m=0,1.
\]
In particular, this implies that
\[
\phi_0(m,x)=\phi_0(m,0), \qquad x\in I^{k-1},\, m=0,1.
\]
The set of pinched cubes in $A$ we denote by $I_\infty(A)$, inside which the set of pinched $k$-cubes we denote by $I_\infty(A)_k$.

The \textbf{face maps} $\delta_i^m:I^{k-1}\to I^k$, for $i=0,\dots, k-1$ and $m=0,1$, we define by
\[
\delta_i^m(x_0,x):=(x_0,\dots, x_{i-1},m,x_i,\dots, x_{k-2}),
\]
which induce face maps $\delta_i^{m*}:I_\infty(A)_k\to I_\infty(A)_{k-1}$. We also define the \textbf{source} (for $m=0$) \textbf{and target} (for $m=1$) \textbf{maps} $w_m:I_\infty(A)_k\to M$ by
\[
w_m(\phi)=\phi_0(m,\dots, m).
\]
\end{definition}

If $\phi\in I_\infty(A)_k$, we will denote $(\phi^*\nabla)_{\del_i}$ interchangeably by $\nabla_{\del_i}$ or $\nabladerpar{}{x_i}$.
 
We can now define integrals along such pinched cubes.

\begin{definition} \label{def:integral-over-pinched-cubes}
Let $\phi\in I_\infty(A)_k$ and $\eta\in\Omega^k(I^k,\End(\phi^*E))$. Then we define the \textbf{integral of $\eta$ over $\phi$} as
\[
\int_\phi\eta:=\int_{I^{k-1}}\int_0^1 \tau_{0\leftarrow x_0,x} \eta(\del_0,\dots,\del_{k-1})(x_0,x)\tau_{x_0\leftarrow 1,x}dx_0dx,
\]
which is an element of $\Hom(E_{w_1(\phi)}, E_{w_0(\phi)})$, where $\tau_{y_0\leftarrow x_0,x}$ is parallel transport with respect to $\phi^*\nabla$ along the curve $t\mapsto ((1-t)x_0+ty_0,x)$. If we want to emphasize the dependence on the connection $\nabla$ we will write $\tau^\nabla$.

If $k>1$ and $\eta\in\Omega^{k-1}(I^k,\End(\phi^*E))$, we also define the \textbf{integral of $\eta$ along the boundary $\del\phi$} as
\[
\int_{\del\phi}\eta:=\sum_{i=1}^{k-1}\sum_{m=0}^1 (-1)^{i+m+1} \int_{(\delta_i^m)^*\phi} (\delta_i^m)^*\eta\in\Hom(E_{w_1(\phi)},E_{w_0(\phi)}).
\]
If $k=1$, we define the integral of $\eta$ along the boundary $\del\phi$ as
\[
\int_{\del\phi}\eta:=\tau_{0\leftarrow 1}\eta(w_1(\phi))-\eta(w_0(\phi))\tau_{0\leftarrow 1}.
\]
\end{definition}

\begin{figure}
    \centering
    \includesvg[width=0.9\textwidth]{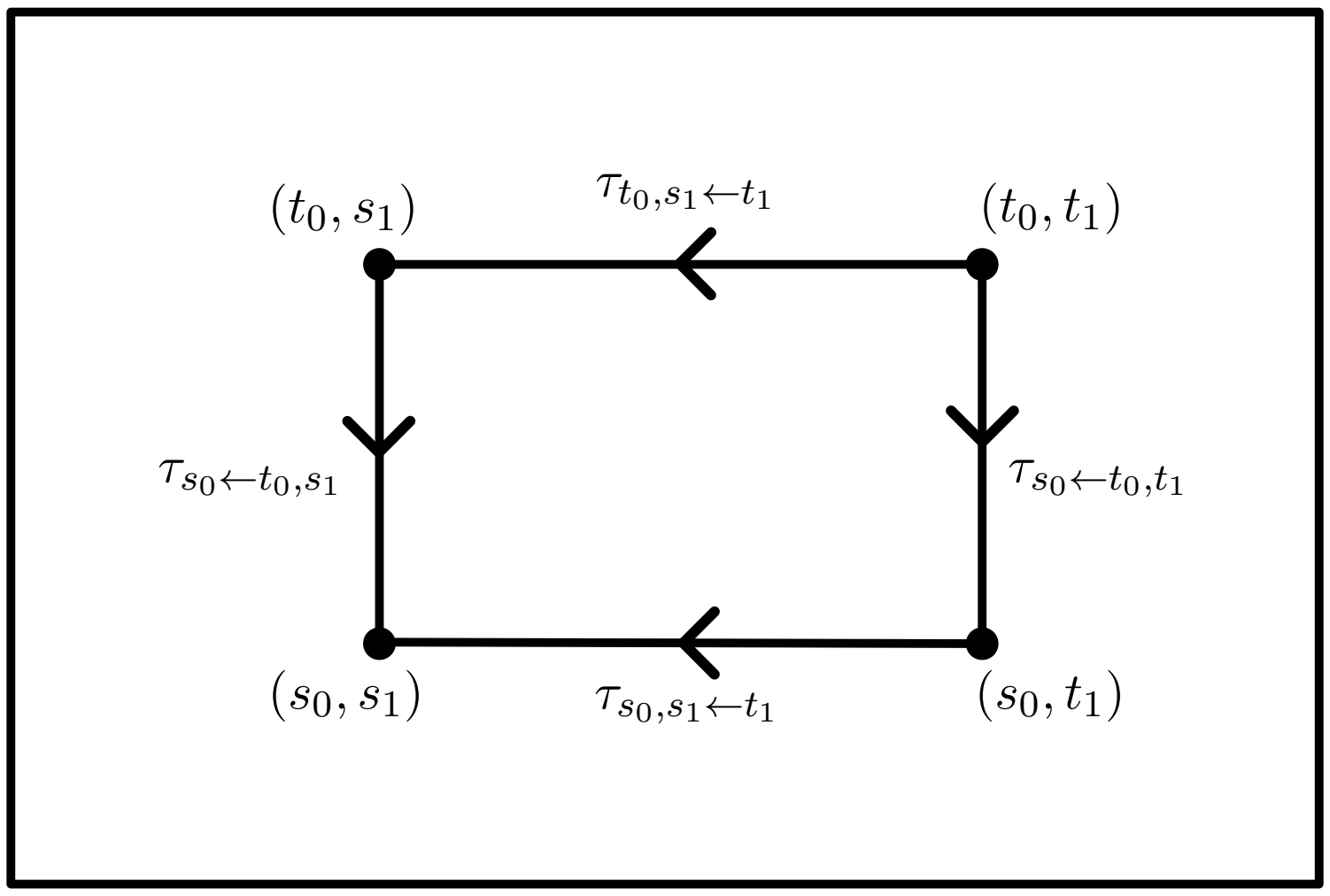}
    \caption{An illustration of the paths in Lemma \ref{lem:difference-parallel-transport}.}
    \label{fig:difference-parallel-transport}
\end{figure}

Definition \ref{def:integral-over-pinched-cubes} will be an appropriate definition of integral if it satisfies some version of Stokes's theorem. To check if Stokes's theorem holds, we first need to control how parallel transport differs when we follow different paths in the cube. This is a well known formula, but we give a proof here for the sake of completeness. We illustrate the situation in Figure \ref{fig:difference-parallel-transport}.

\begin{lemma} \label{lem:difference-parallel-transport}
Let $(V\to I^2,\nabla)$ be a vector bundle with connection and let $(t_0,t_1)$ and $(s_0,s_1)$ be points in $I^2$. If $\tau_{x_0\leftarrow v,x_1}$ and $\tau_{x_0, v\leftarrow x_1}$ denote, respectively, parallel transport along the curves $t\mapsto ((1-t)v+tx_0,x_1)$ and $t\mapsto (x_0, (1-t)x_1+tv)$, then
\begin{multline*}
\tau_{s_0\leftarrow t_0,s_1}\tau_{t_0,s_1\leftarrow t_1} -\tau_{s_0,s_1\leftarrow t_1} \tau_{s_0\leftarrow t_0, t_1} = \\
=\int_{t_1}^{s_1} \int_{t_0}^{s_0} \tau_{s_0,s_1\leftarrow x_1}\tau_{s_0\leftarrow x_0,x_1}K^\nabla(\del_0,\del_1)(x_0,x_1) \tau_{x_0\leftarrow t_0,x_1}\tau_{t_0,x_1\leftarrow t_1}dx_0dx_1.
\end{multline*}
\end{lemma}

\begin{proof}
If $e\in\Gamma(E)$ is such that $\nabla_{\del_0}e=0$ and $\nabla_{\del_1}e(t_0,x_1)=0$, then
\[
\tau_{x_0\leftarrow t_0,x_1} \derpar{}{x_0}\left( \tau_{t_0\leftarrow x_0,x_1}\nabla_{\del_1}e(x_0,x_1)\right)=\nabla_{\del_0}\nabla_{\del_1}e(x_0,x_1)=K^\nabla(\del_0,\del_1)e(x_0,x_1),
\]
which gives
\begin{align*}
\nabla_{\del_1}e(s_0,x_1) &=\tau_{s_0\leftarrow t_0,x_1}\int_{t_0}^{s_0}\derpar{}{x_0}\left(\tau_{t_0\leftarrow x_0,x_1}\nabla_{\del_1}e(x_0,x_1)\right)dx_0 \\
&=\int_{t_0}^{s_0}\tau_{s_0\leftarrow x_0,x_1}K^\nabla(\del_0,\del_1)e(x_0,x_1) dx_0.
\end{align*}
Using that $e(x_0,x_1)=\tau_{x_0\leftarrow t_0,x_1}e(t_0,x_1)=\tau_{x_0\leftarrow t_0,x_1}\tau_{t_0,x_1\leftarrow t_1}e(t_0,t_1)$, we can write
\begin{align*}
\tau_{s_0,t_1\leftarrow x_1}\nabla_{\del_1}e(s_0,x_1) &=\derpar{}{x_1}(\tau_{s_0,t_1\leftarrow x_1}e(s_0,x_1)) \\
&=\derpar{}{x_1}(\tau_{s_0,t_1\leftarrow x_1}\tau_{s_0\leftarrow t_0,x_1}\tau_{t_0,x_1\leftarrow t_1})e(t_0,t_1),
\end{align*}
and
\begin{multline*}
\int_{t_0}^{s_0}\tau_{s_0\leftarrow x_0,x_1}K^\nabla(\del_0,\del_1)e(x_0,x_1) dx_0= \\
=\left(\int_{t_0}^{s_0}\tau_{s_0\leftarrow x_0,x_1}K^\nabla(\del_0,\del_1)(x_0,x_1)\tau_{x_0\leftarrow t_0,x_1}dx_0\right) \tau_{t_0,x_1\leftarrow t_1}e(t_0,t_1).
\end{multline*}
This finally leads to
\begin{multline*}
\tau_{s_0,t_1\leftarrow s_1}\tau_{s_0\leftarrow t_0,s_1}\tau_{t_0,s_1\leftarrow t_1}-\tau_{s_0\leftarrow t_0,t_1} =\int_{t_1}^{s_1}\derpar{}{x_1}(\tau_{s_0,t_1\leftarrow x_1}\tau_{s_0\leftarrow t_0,x_1}\tau_{t_0,x_1\leftarrow t_1})dx_1 \\
=\int_{t_1}^{s_1}\int_{t_0}^{s_0}\tau_{s_0,t_1\leftarrow x_1}\tau_{s_0\leftarrow x_0,x_1}K^\nabla(\del_0,\del_1)(x_0,x_1)\tau_{x_0\leftarrow t_0,x_1}\tau_{t_0,x_1\leftarrow t_1}dx_0dx_1. \qedhere
\end{multline*} 
\end{proof}

\begin{corollary} \label{cor:difference-parallel-transport}
Let $(V\to I^2,\nabla)$ be a vector bundle with connection and $B\in\Gamma(\End(E))$. Then
\begin{align*}
\nabladerpar{}{t_1}(\tau_{t_0\leftarrow x_0,t_1} & B(x_0,t_1)\tau_{x_0\leftarrow s_0,t_1}) = \tau_{t_0\leftarrow x_0,t_1}\nabladerpar{}{t_1}(B(x_0,t_1))\tau_{x_0\leftarrow s_0,t_1} \\
&\qquad -\left( \int_{t_0}^{x_0}\tau_{t_0\leftarrow u,t_1}K^\nabla(\del_0,\del_1)(u,t_1)\tau_{u\leftarrow x_0,t_1} du\right) B(x_0,t_1)\tau_{x_0\leftarrow s_0,t_1} \\
&\qquad -\tau_{t_0\leftarrow x_0,t_1}B(x_0,t_1)\left( \int_{x_0}^{s_0}\tau_{x_0\leftarrow u,t_1}K^\nabla(\del_0,\del_1)(u,t_1)\tau_{u\leftarrow s_0,t_1}du\right).
\end{align*}
\end{corollary}

\begin{proof}
It follows immediately from re-expressing Lemma \ref{lem:difference-parallel-transport} as
\[
\tau_{s_0\leftarrow t_0, t_1}=\tau_{s_0,t_1\leftarrow s_1} \mathfrak{R}\tau_{t_0,s_1\leftarrow t_1},
\]
for
\[
\mathfrak{R}=\tau_{s_0\leftarrow t_0,s_1}-\int_{t_1}^{s_1} \int_{t_0}^{s_0} \tau_{s_0,s_1\leftarrow v}\tau_{s_0\leftarrow u,v}K^\nabla(\del_0,\del_1)(u,v) \tau_{u\leftarrow t_0,v}\tau_{t_0,v\leftarrow s_1}dudv. \qedhere
\]
\end{proof}

We can now check whether Stokes's theorem holds in the simplest case: if $\phi$ is a 2-cube and $\eta\in\phi^*\Omega^1(A,\End(E))$, then by Corollary \ref{cor:difference-parallel-transport},
\begin{align*}
\int_\phi d^\nabla\eta&= \int_0^1\int_0^1 \tau_{0\leftarrow x_0,x_1}d^\nabla\eta(\del_0,\del_1)(x_0,x_1)\tau_{x_0\leftarrow 1,x_1}dx_0dx_1 \\
&=\int_0^1\int_0^1 \tau_{0\leftarrow x_0,x_1}\left( \nabla_{\del_0}(\eta(\del_1))-\nabla_{\del_1}(\eta(\del_0))
\right) (x_0,x_1)\tau_{x_0\leftarrow 1,x_1}dx_0dx_1 \\
&=\int_0^1\int_0^1 \derpar{}{x_0}\left( \tau_{0\leftarrow x_0,x_1}\eta(\del_1)(x_0,x_1)\tau_{x_0\leftarrow 1,x_1} \right)dx_0dx_1 \\
&\qquad -\int_0^1\int_0^1  \tau_{0\leftarrow x_0,x_1}\nabla_{\del_1}(\eta(\del_0))(x_0,x_1)\tau_{x_0\leftarrow 1,x_1} dx_0dx_1 \\
&= -\int_0^1\int_0^1  \derpar{}{x_1}\left( \tau_{0\leftarrow x_0,x_1}\eta(\del_0)(x_0,x_1)\tau_{x_0\leftarrow 1,x_1}\right) dx_0dx_1 \\
&\qquad -\int_0^1\int_0^1\int_0^{x_0}  \tau_{0\leftarrow u,x_1} K^\nabla(\del_0,\del_1)(u,x_1) \tau_{u\leftarrow x_0,x_1}\eta(\del_0)(x_0,x_1)\tau_{x_0\leftarrow 1,x_1} dudx_0dx_1 \\
&\qquad -\int_0^1\int_0^1\int_{x_0}^1  \tau_{0\leftarrow x_0,x_1} \eta(\del_0)(x_0,x_1) \tau_{x_0\leftarrow u,x_1}K^\nabla(\del_0,\del_1)(u,x_1)\tau_{u\leftarrow 1,x_1} dudx_0dx_1 \\
&=\int_{\del\phi} \eta + \text{curvature terms}.
\end{align*}

As expected, the presence of curvature seems to spoil Stokes's theorem. In the face of this, one option would be to give up and conclude that there is no such theorem for endomorphism-valued forms. Alternatively, we can notice that the curvature terms actually have some structure to them. If we write them in the suggestive form
\[
\int_{\del\phi} \eta=\int_{\phi} \left( \eta+K^\nabla\otimes\eta-\eta\otimes K^\nabla \right),
\]
then we recognize the right-hand side as the integral of $D\eta$, for $D$ the differential on the bar complex \eqref{eq:D-bar-complex}. The question then becomes: can we define integrals on elements of $\B (\Omega(A,\End(E))$? The answer to this question is yes, via iterated integrals.

\begin{definition} \label{def:iterated-integral-over-pinched-cubes}
Let $\phi\in I_\infty(A)_k$. We define the \textbf{horizontal integral} over $\phi$ as the bidegree $(-1,0)$ linear operator 
\[
\fint_\phi:\Omega(I^k,\End(\phi^*E))\to\Omega(I^k,\End(\phi^*E))
\]
defined as follows: if $\eta\in\Omega^{\ell+1}(I^k,\End(\phi^*E))$ and $(x_0,x)\in I^k$, then
\[
\left(\fint_\phi\eta\right)(\del_{i_1},\dots, \del_{i_\ell})(x_0,x):=\int_{x_0}^1\tau_{x_0\leftarrow t,x}\eta(\del_0,\del_{i_1},\dots, \del_{i_\ell})(t,x)\tau_{t\leftarrow x_0,x}dt.
\]

If $\eta_i\in\Omega^{p_i}(I^k,\End^{q_i}(\phi^*E))$, for $i=1,\dots,\ell$, are such that $\sum_{i=1}^\ell p_i=k+\ell-1$, then we define the \textbf{integral of $\eta_1\otimes\dots\otimes\eta_\ell$} as the iterated integral
\begin{align*}
\int_\phi(\eta_1& \otimes\dots\otimes\eta_\ell):=\int_\phi\left( \eta_1\wedge \fint_\phi\left( \eta_2\wedge\fint_\phi\left(\eta_3\wedge \dots \fint_\phi\left(\eta_{\ell-1}\wedge\fint_\phi\eta_\ell\right)\dots\right)\right)\right),
\end{align*}
which is an element in $\Hom^{q_1+\dots+q_\ell}(E_{w_1(\phi)},E_{w_0(\phi)})$.

If $\sum_{i=1}^\ell p_i=k+\ell-2$, we also define the \textbf{integral of $\eta_1\otimes\dots\otimes\eta_\ell$ along the boundary $\del\phi$} as
\[
\int_{\del\phi}(\eta_1\otimes\dots\otimes\eta_\ell):=\sum_{i=1}^{k-1}\sum_{m=0}^1 (-1)^{i+m+1} \int_{\delta_m^{i*}\phi} \delta_m^{i*}(\eta_1\otimes\dots\otimes\eta_\ell),
\]
where $\delta_m^{i*}(\eta_1\otimes\dots\otimes\eta_\ell):=\delta_m^{i*}\eta_1\otimes\dots\otimes \delta_m^{i*}\eta_\ell$.
\end{definition}

Iterated integrals were originally introduced by Chen \cite{chen} as solutions to a series of differential equations. He also gave an expression for them, which he used to compute the cohomology of loop spaces. Since then, iterated integrals have been used by several authors in the setting of de Rham's theorem and higher Riemann--Hilbert correspondences \cite{gugenheim,igusa,block-smith,arias-abad-schatz}. The definition we have provided above gives an alternative description without reference to the path space of a manifold. The definition in the spirit of Chen can be roughly described as follows:
\begin{enumerate}
    \item Denote by $\pi_i:(I^k)^\ell\to I^k$ the $i$th projection and by $\ev_\ell: \Delta^\ell\times I^{k-1}\to(I^k)^\ell$ the multi-evaluation map
    \[
    \ev_\ell(t,x)=((t_1,x),\dots, (t_\ell,x)),
    \]
    where
    \[
    \Delta^\ell:=\{(t_1,\dots, t_\ell)\in\R^\ell: 1\geq t_1\geq\dots\geq t_\ell\geq 0\}\subseteq\R^\ell.
    \]

    \item For each $i$ consider the anchored vector bundle map $\phi_i:=\phi\circ\pi_{i*}\circ\ev_{\ell*}:T(\Delta^\ell\times I^{k-1})\to A$, whose base map $\phi_{i0}$ is explicitly given by
    \[
    \phi_{i0}(t,x)=\phi_0(t_i,x),
    \]
    and consider the vector bundle map
        \[
	   \End(\phi_1^*E)\otimes\dots\otimes\End(\phi_\ell^*E)\longrightarrow \End(\phi_\ell^*E)
	   \]
	   defined on $A_1\otimes\dots\otimes A_\ell\in\End(E_{\phi_0(t_1,x)})\otimes\dots\otimes \End(E_{\phi_0(t_\ell,x)})$ by
	   \[
	   \tau_{t_\ell\leftarrow t_1,x}A_1\tau_{t_1\leftarrow t_2,x}A_2\dots A_{\ell-1}\tau_{t_{\ell-1}\leftarrow t_\ell,x}A_\ell\in\End(E_{\phi_0(t_\ell,x)}).
	   \]
	   This induces a wedge product
       \[
       \Omega(\Delta^\ell\times I^k,\End(\phi_1^*E))\otimes\dots \otimes \Omega(\Delta^\ell\times I^k,\End(\phi_\ell^*E))\longrightarrow\Omega(\Delta^\ell\times I^k,\End(\phi_\ell^*E)).
       \]

       \item define the integral of $\eta_1\otimes\dots\otimes \eta_\ell$ by
       \[
        \int_\phi (\eta_1\otimes\dots\otimes \eta_\ell) :=\int_{\Delta^\ell}\int_{I^k}\ev_\ell^*\pi_1^*\eta_1\wedge\dots\wedge \ev_\ell^*\pi_\ell^*\eta_\ell,
       \]
       where the reader can make precise what this last integral should mean, using the parallel transport of $\nabla$.    
\end{enumerate}

Chen's and our definition coincide, they are both given explicitly by the following formula: for $\eta_i$ of bidegree $(p_i,q_i)$, if we write $P_i:=\sum_{j=1}^i(p_j-1)$:
\begin{equation} \label{eq:iterated-integral-explicit}
\begin{aligned}
\int_\phi (\eta_1\otimes \dots\otimes \eta_\ell) &=(-1)^{\sum_{i=1}^{\ell-1}q_i\sum_{j=i+1}^\ell(p_j-1)}\sum_{\sigma\in\Sh(\eta)}(-1)^\sigma \int_{I^{k-1}}\int_0^1\int_{t_\ell}^1\dots \int_{t_2}^1 \\
&\qquad \tau_{0\leftarrow t_\ell,x}\eta_1(\del_0,\del_{\sigma(1)},\dots, \del_{\sigma(P_1)})(t_\ell,x)\tau_{t_\ell\leftarrow t_{\ell-1},x}\dots \\
&\qquad \dots \tau_{t_{i+1}\leftarrow t_i,x}\eta_i(\del_0,\del_{\sigma(P_{i-1}+1)},\dots, \del_{\sigma(P_i}))(t_i,x)\tau_{t_i\leftarrow t_{i-1},x}\dots \\
&\qquad \dots \tau_{t_2\leftarrow t_1,x}\eta_\ell(\del_0,\del_{\sigma(P_{\ell-1}+1)},\dots, \del_{\sigma(P_\ell)})(t_1,x)\tau_{t_1\leftarrow 1,x}, \\
&\qquad dt_1 \dots dt_{\ell-1}dt_\ell dx,
\end{aligned}
\end{equation}
where
\[
\Sh(\eta):=\Sh(p_1-1, p_2-1,\dots, p_\ell-1)\subseteq\mathfrak{S}_{k+1}.
\]

It was known for a long time that such expressions appear even in the case of classical (non-higher) parallel transport, as the following example shows. This case is known among physicists as the path-ordered exponential or Dyson's series.

\begin{example} \label{ex:parallel-transport-as-iterated-integral}
Let $\alpha\in\Omega^1(A,\End^0(E))$. Then, the parallel transports of $\nabla+\alpha$ and $\nabla$ along an $A$-path $a$ are related by
\[
\tau^{\nabla+\alpha}_{a^{-1}}=\tau^\nabla_{a^{-1}}+\sum_{\ell\geq 1} \int_a (-1)^\ell\alpha^{\otimes \ell}.
\]
Indeed, let $s$ be a section of $E$ that is $\nabla$-parallel along $a^{-1}$ and let
\begin{multline*}
\tilde s(t):=s(t)+\sum_{\ell\geq 1}\int_t^1\int_{t_\ell}^1\dots \int_{t_2}^1 \tau^\nabla_{t\leftarrow t_\ell}\alpha(a(t_\ell))\tau^\nabla_{t_\ell\leftarrow t_{\ell-1}}\dots \\
\dots\tau^\nabla_{t_2\leftarrow t_1}\alpha(a(t_1))\tau^\nabla_{t_1\leftarrow 1}s(1) dt_1\dots dt_{\ell-1}dt_\ell.
\end{multline*}

To see that $\tilde s$ converges, let $\alpha_t:=\tau_{1\leftarrow t}\alpha(a(t))\tau_{t\leftarrow 1}\in\End(a_0(1))$ and
\[
P_t:=1+\sum_{\ell\geq 1} \int_t^1\int_{t_\ell}^1\dots \int_{t_2}^1\alpha_{t_\ell}\dots \alpha_{t_1}dt_1\dots dt_{\ell-1}dt_\ell,
\]
so that $\tilde s(t)=\tau_{t\leftarrow 1}P_ts(1)$. Fix a norm $\norm{\cdot}$ on $\End(E_{a_0(1)})$ and let $C:=\max_{t\in I} \norm{\alpha_t}>0$ (if $\alpha\neq 0$). Then
\begin{align*}
\norm{P_t}&\leq 1+\sum_{\ell\geq 1}\int_t^1\int_{t_\ell}^1\dots\int_{t_2}^1 C^\ell dt_1\dots dt_{\ell-1}dt_\ell  \\
&=1+\sum_{\ell\geq 1} C^\ell \int_0^{1-t}\int_0^{t_2}\dots \int_0^{t_\ell}dt_\ell\dots dt_2dt_1 \\
&=1+\sum_{\ell\geq 1} \frac{((1-t)C)^\ell}{\ell!}=e^{(1-t)C},
\end{align*}
so $P_t$ converges uniformly. Moreover, a straightforward computation gives $\frac{d}{dt}P_t =-\alpha_tP_t$, so
\begin{align*}
(\nabla+\alpha)_a\tilde s(t) &=(\nabla_a+\tau_{t\leftarrow 1}\alpha_t\tau_{1\leftarrow t})(\tau_{t\leftarrow 1}P_ts(1)) \\
&= \tau_{t\leftarrow 1}\left(\frac{d}{dt}P_t+\alpha_tP_t\right)s(1)=0.
\end{align*}
\end{example}

We now want to check if some form of Stokes's theorem holds for iterated integrals with respect to the differential $D$ from \eqref{eq:D-bar-complex}. First of all, since we are dealing with pinched cubes, we are interested not in the whole algebra $\Omega(I^k,\End(\phi^*E))$, but only those elements that come from pinched cubes. In particular, we might restrict our attention to the dg-ideal $\phi^*\Omega(A,\End(E))$, whose elements are guaranteed to come from a pinched cube. The reason for this is that elements of $\phi^*\Omega(A,\End(E))$ vanish at the endpoints, and this will be crucial for our arguments.

In order to compute the Stokes theorem, we first need to understand how the horizontal integral $\fint_\phi$ interacts with the differential $d^\nabla$. Concretely, we want to compute the graded commutator
\[
\left[ d^\nabla, \fint_\phi\right]=d^\nabla\circ\fint_\phi+\fint_\phi\circ\  d^\nabla.
\]

\begin{proposition} \label{prop:commutator-d-horint}
Let $\phi\in I_\infty(A)_k$ and $\eta\in\phi^*\Omega(A,\End(E))$ homogeneous. Then
\begin{align*}
d^\nabla\left(\fint_\phi\eta\right)+\fint_\phi(d^\nabla\eta) &=-\eta-\fint_\phi\left(K^\nabla\wedge\fint_\phi\eta\right)+(-1)^{\abs{\eta}}\fint_\phi\left(\eta\wedge\fint_\phi K^\nabla\right) \\
&\qquad -(-1)^{\abs{\eta}}\fint_\phi\eta\wedge\fint_\phi K^\nabla,
\end{align*}
where $\abs{\eta}$ is the total degree of $\eta$.
\end{proposition}

\begin{proof}
Let $\eta\in\phi^*\Omega^p(A,\End^q(E))$. Let $\alpha=(\alpha_0,\dots, \alpha_{p-1})$ be a multiindex, and write $\hat\alpha^j$ for the multiindex $(\alpha_0,\dots, \widehat{\alpha_j},\dots, \alpha_{p-1})$ and $\del_\alpha:=\del_{\alpha_0}\wedge\dots\wedge\del_{\alpha_{p-1}}$.

Assume first that $\alpha_j\neq 0$ for all $j$. Then, by Corollary \ref{cor:difference-parallel-transport},
\begin{align*}
d^\nabla\left(\fint_\phi \eta\right) &(\del_\alpha)(x_0,x)=\sum_{j=0}^{p-1} (-1)^j\nabla_{\del_{\alpha_j}}\left(\left(\fint_\phi\eta\right)(\del_{\hat\alpha_j})\right) (x_0,x) \\
&=\sum_{j=0}^{p-1} (-1)^j\int_{x_0}^1 \nabla_{\del_{\alpha_j}} \left(\tau_{x_0\leftarrow t,x}\eta(\del_0\wedge\del_{\hat\alpha_j})(t,x)\tau_{t\leftarrow x_0,x)}\right)dt \\
&=\sum_{j=0}^{p-1} (-1)^j \left( \int_{x_0}^1\tau_{x_0\leftarrow t,x}\nabla_{\del_{\alpha_j}}\big(\eta(\del_0\wedge\del_{\hat\alpha_j}\big)(t,x)\tau_{t\leftarrow x_0,x)} dt \right. \\
&\qquad + \int_{x_0}^1\int_t^{x_0} \tau_{x_0\leftarrow u, x}K^\nabla(\del_0\wedge\del_{\alpha_j})(u,x)\tau_{u\leftarrow t,x}\eta(\del_0\wedge\del_{\hat\alpha_j})(t,x)\tau_{t\leftarrow x_0,x}dudt \\
&\qquad \left. -\int_{x_0}^1 \int_t^{x_0} \tau_{x_0\leftarrow t,x}\eta(\del_0\wedge\del_{\hat\alpha_j})(t,x)\tau_{t\leftarrow u,x}K^\nabla(\del_0\wedge\del_{\alpha_j})(u,x)\tau_{u\leftarrow x_0,x} dudt\right) \\
&=\sum_{j=0}^{p-1} (-1)^j \left( \int_{x_0}^1\tau_{x_0\leftarrow t,x}\nabla_{\del_{\alpha_j}}\big(\eta(\del_0\wedge\del_{\hat\alpha_j})\big)(t,x)\tau_{t\leftarrow x_0,x)} dt \right. \\
&\qquad + \int_{x_0}^1\int_u^1 \tau_{x_0\leftarrow u, x}K^\nabla(\del_0\wedge\del_{\alpha_j})(u,x)\tau_{u\leftarrow t,x}\eta(\del_0\wedge\del_{\hat\alpha_j})(t,x)\tau_{t\leftarrow x_0,x}dtdu \\
&\qquad -\int_{x_0}^1 \int_t^1 \tau_{x_0\leftarrow t,x}\eta(\del_0\wedge\del_{\hat\alpha_j})(t,x)\tau_{t\leftarrow u,x}K^\nabla(\del_0\wedge\del_{\alpha_j})(u,x)\tau_{u\leftarrow x_0,x} dudt \\
&\qquad \left. +\int_{x_0}^1 \int_{x_0}^1 \tau_{x_0\leftarrow t,x}\eta(\del_0\wedge\del_{\hat\alpha_j})(t,x)\tau_{t\leftarrow u,x}K^\nabla(\del_0\wedge\del_{\alpha_j})(u,x)\tau_{u\leftarrow x_0,x} dudt \right) \\
&=\left(-\fint_\phi (d^\nabla\eta)-\fint_\phi\left(K^\nabla\wedge\fint_\phi\eta\right)+(-1)^{p+q}\fint_\phi\left(\eta\wedge\fint_\phi K^\nabla\right) \right. \\
&\qquad \left. -(-1)^{p+q}\fint_\phi\eta\wedge\fint_\phi K^\nabla \right)(\del_\alpha)(x_0,x) \\
&\qquad +\int_{x_0}^1\derpar{}{t} \left(\tau_{x_0\leftarrow t,x}\eta(\del_\alpha)(t,x)\tau_{t\leftarrow x_0,x)}\right)dt,
\end{align*}
and this gives the result, since $\eta(\del_\alpha)(1,x)=0$ because $\eta\in\phi^*\Omega(A,\End(E))$.

If $\alpha_0=0$, then, since $\del_0\wedge\del_{\hat\alpha_j}\neq 0$ only for $j=0$,
\begin{align*}
d^\nabla\left(\fint_\phi \eta\right)& (\del_\alpha)(x_0,x) = \sum_{j=0}^{p-1} (-1)^j\nabla_{\del_{\alpha_j}}\left(\left(\fint_\phi\eta\right)(\del_{\hat\alpha_j})\right) (x_0,x)\\
&=\sum_{j=0}^{p-1} (-1)^j\int_{x_0}^1 \nabla_{\del_{\alpha_j}} \left(\tau_{x_0\leftarrow t,x}\eta(\del_0\wedge\del_{\hat\alpha_j})(t,x)\tau_{t\leftarrow x_0,x)}\right)dt \\
&= \int_{x_0}^1 \nabla_{\del_{\alpha_0}} \left(\tau_{x_0\leftarrow t,x}\eta(\del_\alpha)(t,x)\tau_{t\leftarrow x_0,x)}\right)dt \\
&= \int_{x_0}^1 \derpar{}{x_0}\left( \tau_{x_0\leftarrow t,x}\eta(\del_0\wedge\del_{\hat\alpha_0})(t,x)\tau_{t\leftarrow x_0,x)}\right)dt \\
&=-\eta(\del_\alpha)(x_0,x). \qedhere
\end{align*}
\end{proof}

We can now finally state and prove our version of the non-abelian Stokes's theorem.

\begin{theorem}[Non-abelian Stokes's theorem] \label{thm:non-abelian-stokes-for-pinched-cubes}
Let $\phi$ be a pinched cube and $\eta\in\B(\phi^*\Omega(A,\End(E)))$. Then
\[
\int_{\del\phi}\eta=\int_\phi D\eta.
\]
\end{theorem}

\begin{proof}
Assume $\phi$ is a $(k+1)$-cube. We will prove it for $\eta=\eta_1\otimes\dots\otimes\eta_\ell$ by induction on $\ell$.

Let $\ell=1$, so that $\eta\in\phi^*\Omega^k(A,\End^q(E))$. Then, by Corollary \ref{cor:difference-parallel-transport}, if we use the same notation as in the proof of Proposition \ref{prop:commutator-d-horint} for the multiindex $\alpha=(\alpha_0,\dots, \alpha_k)$ with $\alpha_i=i$,
\begin{align*}
\int_\phi d^\nabla\eta &= \int_{I^k}\int_0^1 \tau_{0\leftarrow x_0,x}d^\nabla\eta(\del_\alpha)(x_0,x)\tau_{x_0\leftarrow 1,x}dx_0dx \\
&=\sum_{i=0}^k(-1)^i\int_{I^k}\int_0^1\tau_{0\leftarrow x_0,x}\nabla_{\del_{\alpha_i}}\big( \eta(\del_{\hat\alpha_i})\big) (x_0,x)\tau_{x_0\leftarrow 1,x}dx_0dx \\
&=\int_{I^k}\int_0^1\derpar{}{x_0}\big( \tau_{0\leftarrow x_0,x} \eta(\del_{\hat\alpha_0})(x_0,x)\tau_{x_0\leftarrow 1,x} \big)dx_0dx \\
&\qquad +\sum_{i=1}^k(-1)^i\left( \int_{I^k}\int_0^1\derpar{}{x_i}\big( \tau_{0\leftarrow x_0,x} \eta(\del_{\hat\alpha_i})(x_0,x)\tau_{x_0\leftarrow 1,x} \big)dx_0dx \right.\\
&\qquad -\int_{I^k}\int_0^1\int_0^{x_0} \tau_{0\leftarrow u,x}K^\nabla(\del_0\wedge\del_{\alpha_i})(u,x)\tau_{u\leftarrow x_0,x} \eta(\del_{\hat\alpha_i})(x_0,x)\tau_{x_0\leftarrow 1,x}dudx_0dx \\
&\qquad \left. -\int_{I^k}\int_0^1\int_{x_0}^1\tau_{0\leftarrow x_0,x}\eta(\del_{\hat\alpha_i})(x_0,x)\tau_{x_0\leftarrow u,x} K^\nabla(\del_0\wedge\del_{\alpha_i})(u,x)\tau_{u\leftarrow 1,x}dudx_0dx\right).
\end{align*}
Since $\eta(\del_{\hat\alpha_0})(0,x)=0$ and $\eta(\del_{\hat\alpha_0})(1,x)=0$, the first term vanishes. For the second term, if we write $x_i^m:=(x_1,\dots, x_{i-1},m,x_{i+1},\dots, x_k)$, we have that
\begin{align*}
\int_{I^k}\int_0^1\derpar{}{x_i} & \big( \tau_{0\leftarrow x_0,x} \eta(\del_{\hat\alpha_i})(x_0,x)\tau_{x_0\leftarrow 1,x} \big)dx_0dx = \\
&=\sum_{m=0}^1 (-1)^{m+1}\int_0^1\dots\int_0^1 \tau_{0\leftarrow x_0,x_i^m}\eta(\del_{\hat\alpha_i})(x_0,x_i^m)\tau_{x_0\leftarrow 1,x_i^m}dx_0\dots\widehat{dx_i}\dots dx_k,
\end{align*}
so that the second terms assembles into $\int_{\del\phi}\eta$. The third term gives
\[
-\int_\phi\left(K^\nabla\wedge\fint_\phi\eta\right),
\]
whereas the fourth one gives
\[
(-1)^{k+q}\int_\phi\left(\eta\wedge\fint_\phi K^\nabla\right).
\]
This ends the proof for the case $\ell=1$.

Assume now that the result holds up to a certain $\ell-1$. Then
\begin{align*}
\int_{\del \phi} &(\eta_1\otimes\dots\otimes\eta_\ell) = \int_{\del \phi}\left( \eta_1\otimes\dots \otimes\eta_{\ell-2}\otimes\left( \eta_{\ell-1}\wedge\fint_\phi\eta_\ell\right) \right)\\
&=\int_\phi D\left( \eta_1\otimes\dots \otimes\eta_{\ell-2}\otimes\left( \eta_{\ell-1}\wedge\fint_\phi\eta_\ell\right) \right) \\
&=\int_\phi D(\eta_1\otimes\dots\otimes\eta_{\ell-2})\otimes\left( \eta_{\ell-1}\wedge\fint_\phi\eta_\ell\right) + \\
&\qquad + (-1)^{[\eta]_{\ell-2}}\int_\phi\left( \eta_1\otimes\dots\otimes\eta_{\ell-2}\otimes d^\nabla\left( \eta_{\ell-1}\wedge\fint_\phi\eta_\ell\right)\right) \\
&\qquad +(-1)^{[\eta]_{\ell-2}}\int_\phi\left(\eta_1\otimes\dots\otimes\eta_{\ell-3}\otimes \left( \eta_{\ell-2}\wedge \eta_{\ell-1}\wedge\fint_\phi\eta_\ell\right) \right) \\
&\qquad +(-1)^{[\eta]} \int_\phi\left( \eta_1\otimes\dots\otimes\eta_{\ell-2}\otimes \left( \eta_{\ell-1}\wedge\fint_\phi\eta_\ell\right) \otimes K^\nabla \right).
\end{align*}
By Proposition \ref{prop:commutator-d-horint},
\begin{align*}
d^\nabla &\, \left( \eta_{\ell-1}\wedge\fint_\phi\eta_\ell\right) +(-1)^{[\eta_{\ell-1}]+[\eta_\ell]}\eta_{\ell-1}\wedge\fint_\phi\eta_\ell\wedge\fint_\phi K^\nabla \\
&=d^\nabla\eta_{\ell-1}\wedge \fint_\phi\eta_\ell-(-1)^{\abs{\eta_{\ell-1}}}\eta_{\ell-1}\wedge\fint_\phi d^\nabla\eta_\ell \\
&\qquad -(-1)^{\abs{\eta_{\ell-1}}}\eta_{\ell-1}\wedge\eta_\ell-(-1)^{\abs{\eta_{\ell-1}}} \eta_{\ell-1}\wedge \fint_\phi\left( K^\nabla\wedge \fint_\phi \eta_\ell\right) \\
&\qquad+(-1)^{\abs{\eta_{\ell-1}}+\abs{\eta_\ell}}\eta_{\ell-1}\wedge \fint_\phi\left(\eta_\ell\wedge\fint_\phi K^\nabla\right) \\
&\qquad -(-1)^{\abs{\eta_{\ell-1}}+\abs{\eta_\ell}} \eta_{\ell-1}\wedge \fint_\phi\eta_\ell\wedge \fint_\phi K^\nabla \\
&\qquad +(-1)^{[\eta_{\ell-1}]+[\eta_\ell]}\eta_{\ell-1}\wedge\fint_\phi\eta_\ell\wedge\fint_\phi K^\nabla,
\end{align*}
and this finally proves the result.
\end{proof}

The use of the codifferential $D$ from \eqref{eq:D-bar-complex} in Theorem \ref{thm:non-abelian-stokes-for-pinched-cubes} indicates that this theorem is related to an $A_\infty$-map, as we will now explain. Before that, we have to endow $I_\infty(A)$ with a smooth structure. For our purposes it is enough to work with smooth structures in a weak sense: diffeology (see Appendix \ref{app:diffeology} for the basics on diffeologies).

\begin{lemma}
$I_\infty(A)$ is a diffeological space, where the set of plots for $I_\infty(A)_k$ is the set of maps $p:U\to I_\infty(A)_k$ for some open $U\subseteq \R^n$ such that the induced map
\[
\begin{array}{cccc}
\hat p_0: & U\times I^k & \longrightarrow & M \\
& (y,x) & \longmapsto & p(y)_0(x)
\end{array}
\]
is smooth and for every $\alpha\in\Omega^\ell(A)$ the map $p_\alpha^*:U\to \Omega^\ell(I^k)$ given by $y\mapsto p(y)^*\alpha$ is smooth.
\end{lemma}

\begin{proof}
Firstly, if $p:\R^n\to I_\infty(A)_k$ is constant, say $p(x)=\phi$ for some $\phi\in I_\infty(A)_k$, then $\hat p_0(y,x)=\phi_0(x)$ and $p_\alpha^*(x)=\phi^*\alpha$, both of which are smooth, so $p$ is a plot. Secondly, if $p:U\to I_\infty(A)_k$ is a map such that for every $x\in U$ there is a neighborhood $V$ containing $x$ such that $p|_V$ is a plot, then by the glueing property for smooth maps $U\times I^k\to M$ and for smooth functions on $U\times I^k$ we have that $p$ is a plot. Lastly, if $p:U\to I_\infty(A)_k$ is a plot and $F:V\to U$ is smooth, for some open $V\subseteq \R^m$, then
\[
\widehat{(p\circ F)}_0=\hat p_0\circ (F\times\id)
\]
and
\[
(p\circ F)^*_\alpha=p_\alpha^*\circ F,
\]
both of which are smooth, so $p\circ F$ is a plot.
\end{proof}

\begin{definition}
We define the bigraded vector space
\[
C_{I_\infty(A)}(\End(E))=\bigoplus_{p,q\geq 0} C^p_{I_\infty(A)}(\End^q(E))
\]
with
\[
C^p_{I_\infty(A)}(\End^q(E)):=\Gamma(I_\infty(A)_p, \Hom^q(w_1^*E,w_0^*E)).
\]
We define the bidegree $(1,0)$ map $d$, which we call the \textbf{differential}, defined on $C^p_{I_\infty(A)}(\End^q(E))$ for $p>0$ by
\[
d:=\sum_{i=1}^p\sum_{m=0}^1 (-1)^{i+m+1}d_i^{m*},
\]
where $d_i^m=\delta_i^{m*}:I_\infty(A)_{p+1}\to I_\infty(A)_p$, and for $p=0$ by
\[
dB(\phi):=\tau_{0\leftarrow 1}B(w_1(\phi))-B(w_0(\phi))\tau_{0\leftarrow 1},
\]
for $B\in\Gamma(\End(E))$.
\end{definition}

\begin{lemma}
If $i\leq j$, we have that $\delta_i^m\delta_j^\ell=\delta_{j+1}^\ell\delta_i^m$. As a consequence, the pair $(C_{I_\infty(A)}(\End(E)),d)$ is a cochain complex. In particular, it is an $A_\infty$-algebra with $m_n=0$ for $n\geq 2$.
\end{lemma}

\begin{proof}
Let $(x_0,x)\in I^p$. Then
\[
\delta_i^m\delta_j^\ell(x_0,x)=(x_0,\dots, x_{i-1},m,x_i,\dots,x_{j-1},\ell,x_j,\dots, x_{p-1})=\delta_{j+1}^\ell\delta_i^m(x_0,x).
\]
This implies that $d^\ell_jd^m_i=d^m_id^\ell_{j+1}$. Then $d^2$ on $C^p_{I_\infty(A)}(\End^q(E))$ for $p>0$ is
\begin{align*}
d^2&= \sum_{i=1}^{p+1}\sum_{j=1}^p\sum_{m,\ell=0}^1 (-1)^{i+j+m+\ell} (d^m_jd^\ell_i)^* \\
&=\sum_{j=1}^p\sum_{i=1}^j \sum_{m,\ell=0}^1 (-1)^{i+j+m+\ell} (d^m_jd^\ell_i)^* + \sum_{j=1}^p\sum_{i=j+1}^{p+1} \sum_{m,\ell=0}^1 (-1)^{i+j+m+\ell} (d^m_jd^\ell_i)^* \\
&=\sum_{j=1}^p\sum_{i=1}^j \sum_{m,\ell=0}^1 (-1)^{i+j+m+\ell} (d^\ell_id^m_{j+1})^* + \sum_{j=1}^p\sum_{i=j+1}^{p+1} \sum_{m,\ell=0}^1 (-1)^{i+j+m+\ell} (d^m_jd^\ell_i)^* \\
&=\sum_{j=1}^p\sum_{i=j+1}^{p+1} \sum_{m,\ell=0}^1 (-1)^{i+j+m+\ell}\left(  (d^m_jd^\ell_i)^* - (d^\ell_j d^m_i)^*\right)=0.
\end{align*}
If $p=0$ and $B\in\Gamma(\End^q(E))$ and $\phi\in I_\infty(A)_2$, then
\begin{align*}
d^2B(\phi) &=dB(\phi\circ\delta_1^0)-dB(\phi\circ\delta_1^1) \\
&=\tau_{0\leftarrow 1}B(w_1(\phi\circ\delta_1^0))-B(w_0(\phi\circ\delta_1^0))\tau_{0\leftarrow 1} \\
&\qquad -\tau_{0\leftarrow 1}B(w_1(\phi\circ\delta_1^1))+B(w_0(\phi\circ\delta_1^1))\tau_{0\leftarrow 1},
\end{align*}
and the result follows from
\begin{align*}
w_1(\phi\circ\delta_1^0) &=\phi_0(1,0)=\phi_0(1,1)=w_1(\phi\circ\delta_1^1) \\
w_0(\phi\circ\delta_1^0) &=\phi_0(0,0)=\phi_0(0,1)=w_0(\phi\circ\delta_1^1). \qedhere
\end{align*}
\end{proof}

\begin{definition}
For $\ell\geq 1$ we define the map
\[
S^\nabla_\ell:\Omega(A,\End(E))[1]^{\otimes\ell}\longrightarrow C_{I_\infty(A)}(\End(E))[1]
\]
by setting
\[
\hat S^\nabla_\ell(\eta_1\otimes\dots\otimes \eta_\ell)(\phi):=\int_\phi\phi^*(\eta_1\otimes\dots\otimes\eta_\ell), \qquad\text{if $\ell>0$}
\]
and
\[
\hat S^\nabla_0(1)(\phi):=\tau^\nabla_{\phi^{-1}},
\]
where in this case $\phi:TI\to A$ is an $A$-path and $\phi^{-1}$ denotes the reverse path.
\end{definition}

Then the non-abelian Stokes theorem can be rephrased as follows.

\begin{proposition}
The collection $\{S^\nabla_\ell\}_{\ell\geq 0}$ defines an $A_\infty$-map
\[
S^\nabla:\Omega(A,\End(E))\longrightarrow C_{I_\infty(A)}(\End(E)).
\]
\end{proposition}

\begin{proof}
Using Corollary \ref{cor:A-infty-map-for-curved-dgas}, since $C_{I_\infty(A)}(\End(E))$ is a cochain complex, we see that $S$ is an $A_\infty$-map if and only if
\[
d(\hat S^\nabla_\ell(\eta_1\otimes\dots\otimes\eta_\ell)) =\hat S^\nabla(D(\eta_1\otimes\dots\otimes\eta_\ell))
\]
for $\ell>0$, and $d(\hat S^\nabla_0(1))=\hat S^\nabla_1(K^\nabla)$. The first equation is equivalent to Theorem \ref{thm:non-abelian-stokes-for-pinched-cubes}, by noticing that
\begin{align*}
d(\hat S^\nabla_\ell(\eta_1\otimes\dots\otimes\eta_\ell))(\phi)&=\sum_{i=1}^k\sum_{m=0}^1(-1)^{i+m+1}\int_{\delta_m^{i*}\phi} \delta_m^{i*}\phi^*(\eta_1\otimes\dots\otimes\eta_\ell) \\
&=\int_{\del\phi}\phi^*(\eta_1\otimes\dots\otimes\eta_\ell), \\
\hat S^\nabla(D(\eta_1\otimes\dots\otimes\eta_\ell)) &= \int_\phi \phi^*(D(\eta_1\otimes\dots\otimes\eta_\ell)).
\end{align*}
The second condition is equivalent to the fact that for any $\phi\in I_\infty(A)_2$ we have
\[
\tau^\nabla_{\phi(\cdot,0)^{-1}}-\tau^\nabla_{\phi(\cdot,1)^{-1}}=\int_\phi \phi^*K^\nabla,
\]
which holds by Lemma \ref{lem:difference-parallel-transport}.
\end{proof}

The following remark will be very important in Section \ref{sect:riemann-hilbert-functor}.

\begin{remark} \label{rem:integral-different-bundles}
Everything that we did here can be generalized by letting the $\eta_i$ be defined on different vector bundles. Indeed, if $(E_i\to M,\nabla^i)$ are graded vector bundles with degree preserving $A$-connections, for $i=0,\dots, \ell$, and
\[
\eta_i\in \Omega^{p_i}(I^k,\Hom^{q_i}(\phi^*E_i,\phi^*E_{i-1}))
\]
are such that $\sum_{i=1}^\ell( p_i-1)+1=k$, then one can adapt Definition \ref{def:iterated-integral-over-pinched-cubes} to
\[
\int_\phi(\eta_1\otimes\dots\otimes\eta_\ell)\in\Hom^{q_1+\dots+q_\ell}(E_{\ell,w_1(\phi)},E_{0,w_0(\phi)}).
\]
Explicitly, this integral is given by \eqref{eq:iterated-integral-explicit} but replacing the parallel transport $\tau$ by the appropriate parallel transports $\tau^i$. This new version also satisfies a non-abelian Stokes's theorem, the analog of Theorem \ref{thm:non-abelian-stokes-for-pinched-cubes}.
\end{remark}

\subsection{Iterated integrals for simplices}
The general theory in this section closely follows section 4.2 of \cite{igusa} and section 3.2 of \cite{arias-abad-schatz}. However, instead of working with the path spaces of the standard simplices, as done there, we work with the simplices themselves. 
Throughout this section $A\to M$ will again be a Lie algebroid and $(E\to M,\nabla)$ a graded vector bundle with a degree-preserving $A$-connection.

First of all, let us set basic definitions and notation on simplices. We think of the standard simplex as the space
\[
\Delta^k:=\{(t_1,\dots, t_k)\in\R^k: 1\geq t_1\geq\dots\geq t_k\geq 0\}\subseteq\R^k.
\]
The collection $\{\Delta^k\}_k$ is a cosimplicial set with face maps $\delta_i:\Delta^k\to\Delta^{k+1}$ given by
\[
\delta_i(t_1,\dots, t_k):=\left\{
\begin{array}{ll}
(1,t_1,\dots,t_k), & i=0 \\
(t_1,\dots, t_{i-1},t_i,t_i, t_{i+1},\dots, t_k), & i=1,\dots, k \\
(t_1, \dots, t_k,0), & i=k+1
\end{array}
\right.
\]
and degeneracy maps $\sigma_i:\Delta^k\to\Delta^{k-1}$ given by
\[
\sigma_i(t_1,\dots, t_k) := (t_1,\dots, \widehat{t_i},\dots, t_k).
\]
The front-face, back-face and vertex maps
\begin{align*}
\alpha_i &:= \delta_k\dots\delta_{i+1}:\Delta^i\to \Delta^k, \\
\beta_i &:= \delta_0\stackrel{(k-i)}{\dots}\delta_0:\Delta^i\to \Delta^k, \\
\upsilon_i&:= \beta_{k-i}\alpha_0:\Delta^0\to\Delta^k
\end{align*}
are given explicitly by
\begin{align*}
\alpha_i(t_1,\dots, t_i) &=(t_1,\dots, t_i,0,\dots, 0), \\
\beta_i(t_1,\dots, t_i) &=(1,\dots, 1, t_1,\dots, t_i), \\
\upsilon_i &=(1,\stackrel{(i)}{\dots}, 1,0,\stackrel{(k-i)}{\dots},0).
\end{align*}

\begin{definition}
The \textbf{homotopy groupoid $\Pi_\infty(A)$ of $A$} is defined by setting $\Pi_\infty(A)_k$ to be the space of Lie algebroid maps $\sigma:T\Delta^k\to A$ and face and degeneracy maps given, respectively, by
\begin{align*}
d_i\sigma &:= \sigma\circ\delta_{i*}, \\
s_i\sigma &:= \sigma\circ \sigma_{i*}.
\end{align*}
\end{definition}

\begin{lemma}
$\Pi_\infty(A)$ is a diffeological higher groupoid.
\end{lemma}

\begin{proof}
See Lemma \ref{lem:pi-infty-higher-lie-algebroids} for the case of a higher Lie algebroid.
\end{proof}

In Section \ref{sect:iterated-integrals-cubes} we showed how to integrate over pinched cubes, which eventually resulted in an $A_\infty$-map $ S^\nabla:\Omega(A,\End(E))\to C_{I_\infty(A)}(\End(E))$. We would like now to be able to integrate over simplices, instead of cubes. The strategy is easy: find a way to identify pinched cubes with simplices and use $S^\nabla$ to integrate. Said otherwise, we would like to use $S^\nabla$ to define an $A_\infty$-map
\[
\Sigma^\nabla:\Omega(A,\End(E))\longrightarrow C_{\Pi_\infty(A)}(\End(E)).
\]
Here the we regard $C_{\Pi_\infty(A)}(\End(E))$ as a dga, with the product as in Lemma \ref{lem:hom-C(G,E)} and differential $d=[D_0,\cdot]$, where $D_0$ is the canonical flat $G$-connection on $E$ from Proposition \ref{prop:connections-simplicial-description}. Explicitly, $d$ acts on $C^p_{\Pi_\infty(A)}(\End^q(E))$ by
\[
d:=(-1)^{p+q}\sum_{i=1}^p(-1)^id_i^*.
\]

\begin{definition} \label{def:Sigma}
We define a \textbf{cube-to-simplex family} to be a family of maps $\theta=\{\theta_k\}_k$, where $\theta_k:I^k\to\Delta^k$ is piecewise smooth with $\theta_k(0,x)=v_0$ and $\theta_k(1,x)=v_k$ for all $x\in I^{k-1}$, and such that $\theta_1=\id$. By piecewise smooth we mean that there are $a_0=0<a_1<\dots<a_{n-1}<1=a_n$ such that $\theta|_{(a_{i-1},a_i)\times I^{k-1}}$ is smooth for all $i=1,\dots, n$.

We define $\Sigma^\nabla$ subject to $\theta$ as follows: for $\ell>0$, $\eta_i\in\Omega^{p_i}(A,\End^{q_i}(E))$ and $\sigma\in\Pi_\infty(A)_k$ we set
\begin{equation} \label{eq:def-of-sigma}
\hat\Sigma^\nabla_\ell(\eta_1\otimes\dots\otimes\eta_\ell)(\sigma):=(-1)^{\nu(\eta_1,\dots,\eta_\ell)} \hat S^\nabla_\ell(\eta_1\otimes\dots\otimes\eta_\ell)(\theta_k^*\sigma),
\end{equation}
where
\[
\nu(\eta_1,\dots,\eta_\ell):=\sum_{j=1}^\ell \frac{1}{2}p_j(p_j-1)+\sum_{j=1}^{\ell}(p_j-1)\sum_{s=1}^{j-1}(p_s-1)+\left( \sum_{j=1}^\ell q_j\right)\sum_{j=1}^\ell(p_j-1),
\]
and for $\ell=0$ and $\sigma\in\Pi_\infty(A)_1$ we set
\begin{equation} \label{eq:def-of-sigma-0}
\hat\Sigma^\nabla_0(1)(\sigma):=\tau^\nabla_{\sigma^{-1}}.
\end{equation}
\end{definition}

The question now becomes which cube-to-simplex families make $\Sigma^\nabla$ into an $A_\infty$-map.

\begin{theorem} \label{thm:good-cube-to-simplex-family-first-version}
For the map $\Sigma^\nabla$ as in Definition \ref{def:Sigma} to be an $A_\infty$-map it is enough if the following holds: for all $\eta_s\in \Omega^{p_s}(A,\End^{q_s}(E))$, with $s=1,\dots, \ell$, if we denote $P_i^j:=\sum_{k=i}^j(p_k-1)$ and $Q_i^j:=\sum_{k=i}^jq_k$, and $k:=1+P_1^\ell$, then for all $\sigma\in\Pi_\infty(A)_{k+1}$ and $i=1,\dots,k$,
\begin{equation} \label{eq:good-cube-to-simplex:easy}
\hat S^\nabla_\ell(\eta_1\otimes\dots\otimes\eta_\ell)(\delta_i^{1*}\theta_{k+1}^*\sigma)=\hat S^\nabla_\ell(\eta_1\otimes\dots\otimes\eta_\ell)(\theta_k^*\delta_i^*\sigma),
\end{equation}
and
\begin{align}
& \hat S^\nabla_\ell (\eta_1 \otimes\dots\otimes\eta_\ell)(\delta_i^{0*}\theta_{k+1}^*\sigma) = \label{eq:good-cube-to-simplex:hard}\\
&\quad =\sum_{j=0}^\ell (-1)^{Q_1^j P_{j+1}^\ell}
\hat S^\nabla_j(\eta_1\otimes\dots\otimes\eta_j)(\theta_i^*\alpha_i^*\sigma) \nonumber\hat S^\nabla_{\ell-j}(\eta_{j+1}\otimes\dots\otimes\eta_\ell)(\theta_{k-i}^*\beta_{k-i}^*\sigma). \nonumber
\end{align}
\end{theorem}

\begin{proof}
We want to deduce \eqref{eq:good-cube-to-simplex:easy} and \eqref{eq:good-cube-to-simplex:hard} from imposing Corollary \ref{cor:A-infty-map-for-curved-dgas} on $\Sigma^\nabla$. The right-hand side of the first equation in Corollary \ref{cor:A-infty-map-for-curved-dgas} is, for $\sigma\in\Pi_\infty(A)_{k+1}$,
\begin{align*}
& \sum_{i=1}^\ell (-1)^{P_1^{i-1}+Q_1^{i-1}} \hat \Sigma^\nabla_\ell (\eta_1\otimes\dots\otimes d^\nabla\eta_i\otimes\dots\otimes\eta_\ell)(\sigma) \\
&\qquad +\sum_{i=1}^{\ell-1} (-1)^{P_1^i+Q_1^i} \hat\Sigma^\nabla_{\ell-1}(\eta_1\otimes\dots\otimes(\eta_i\wedge\eta_{i+1})\otimes\dots\otimes\eta_\ell)(\sigma) \\
&\qquad +\sum_{i=1}^{\ell+1}(-1)^{P_1^{i-1}+Q_1^{i-1}} \hat\Sigma^\nabla_{\ell+1}(\eta_1\otimes\dots\otimes\eta_{i-1}\otimes K^\nabla\otimes\eta_i\otimes\dots\otimes\eta_\ell)(\sigma) \\
&=\sum_{i=1}^\ell (-1)^{P_1^{i-1}+Q_1^{i-1}+\nu(\eta_1,\dots, d^\nabla\eta_i,\dots,\eta_\ell)} \hat S^\nabla_\ell(\eta_1\otimes\dots\otimes d^\nabla\eta_i\otimes\dots\otimes\eta_\ell)(\theta_{k+1}^*\sigma) \\
&\qquad +\sum_{i=1}^{\ell-1} (-1)^{P_1^i+Q_1^i+\nu(\eta_1,\dots,\eta_i\wedge\eta_{i+1},\dots,\eta_\ell)} \hat S^\nabla_{\ell-1}(\eta_1\otimes\dots\otimes(\eta_i\wedge\eta_{i+1})\otimes\dots\otimes\eta_\ell)(\theta_{k+1}^*\sigma) \\
&\qquad +\sum_{i=1}^{\ell+1}(-1)^{P_1^{i-1}+Q_1^{i-1}+\nu(\eta_1,\dots,\eta_{i-1},K^\nabla,\eta_i,\dots,\eta_\ell)} \hat S^\nabla_{\ell+1}(\eta_1\otimes\dots\otimes\eta_{i-1}\otimes K^\nabla\otimes\eta_i\otimes\dots\otimes\eta_\ell)(\theta_{k+1}^*\sigma).
\end{align*}
But notice that
\begin{align*}
\nu(\eta_1,\dots,d^\nabla\eta_i,\dots,\eta_\ell) -\nu(\eta_1,\dots,\eta_\ell) & = p_i + P_1^{i-1}+P_{i+1}^\ell+Q_1^\ell\\
&=1+P_1^\ell+Q_1^\ell \\
\nu(\eta_1,\dots,\eta_i\wedge \eta_{i+1},\dots,\eta_\ell)-\nu(\eta_1,\dots,\eta_\ell) &= p_ip_{i+1}-p_{i+1}(p_i-1)+P_1^i+P_{i+2}^\ell+Q_1^\ell \\
&=1+P_1^\ell+Q_1^\ell \\
\nu(\eta_1,\dots,\eta_{i-1},K^\nabla,\eta_i,\dots,\eta_\ell)-\nu(\eta_1,\dots,\eta_\ell) &=1+P_1^\ell+Q_1^\ell.
\end{align*}
Hence, using as well that $S^\nabla$ is an $A_\infty$-map, we can rewrite the previous expression as
\begin{align} \label{eq:proof-Sigma-A-infty:rhs-corollary} 
& (-1)^{1+P_1^\ell+Q_1^\ell+\nu(\eta_1,\dots,\eta_\ell)} \left( \sum_{i=1}^\ell (-1)^{P_1^{i-1}+Q_1^{i-1}} \hat S^\nabla_\ell (\eta_1\otimes\dots\otimes d^\nabla\eta_i\otimes\dots\otimes\eta_\ell)(\theta_{k+1}^*\sigma) \right. \nonumber \\
&\qquad\qquad +\sum_{i=1}^{\ell-1} (-1)^{P_1^i+Q_1^i} \hat S^\nabla_{\ell-1}(\eta_1\otimes\dots\otimes(\eta_i\wedge\eta_{i+1})\otimes\dots\otimes\eta_\ell)(\theta_{k+1}^*\sigma) \\
&\qquad\qquad \left. +\sum_{i=1}^{\ell+1}(-1)^{P_1^{i-1}+Q_1^{i-1}} \hat S^\nabla_{\ell+1}(\eta_1\otimes\dots\otimes\eta_{i-1}\otimes K^\nabla\otimes\eta_i\otimes\dots\otimes\eta_\ell)(\theta_{k+1}^*\sigma) \right) \nonumber \\
&\qquad =(-1)^{1+P_1^\ell+Q_1^\ell+\nu(\eta_1,\dots,\eta_\ell)} d(\hat S^\nabla_\ell(\eta_1\otimes\dots\otimes\eta_\ell))(\theta_{k+1}^*\sigma) \nonumber \\
&\qquad =(-1)^{1+P_1^\ell+Q_1^\ell+\nu(\eta_1,\dots,\eta_\ell)} \sum_{i=1}^k (-1)^i\left( \hat S^\nabla_\ell (\eta_1\otimes\dots\otimes\eta_\ell)(\delta_i^{1*}\theta_{k+1}^*\sigma) \right. \nonumber \\
&\qquad\qquad\qquad\qquad\qquad\qquad\qquad\qquad \left. -\hat S^\nabla_\ell (\eta_1\otimes\dots\otimes\eta_\ell)(\delta_i^{0*}\theta_{k+1}^*\sigma) \right). \nonumber
\end{align}

On the other hand, the left-hand side of the first equation in Corollary \ref{cor:A-infty-map-for-curved-dgas} is
\begin{align*}
 d( & \hat\Sigma^\nabla_\ell(\eta_1\otimes\dots\otimes\eta_\ell))(\sigma)+\sum_{i=0}^\ell (-1)^{P_1^i+Q_1^i}\left( \hat\Sigma^\nabla_i(\eta_1\otimes\dots\otimes\eta_i)\hat\Sigma^\nabla_{\ell-j}(\eta_{i+1}\otimes\dots\otimes\eta_\ell)\right)(\sigma) \\
&=(-1)^{1+P_1^\ell+Q_1^\ell+\nu(\eta_1,\dots,\eta_\ell)}\sum_{i=1}^k(-1)^i \hat S^\nabla_\ell(\eta_1\otimes\dots\otimes\eta_\ell)(\theta_{k-1}^*\delta_i^*\sigma) \\
&\qquad +\sum_{i=0}^\ell (-1)^{P_1^i+Q_1^i+(P_1^i+1)(P_{i+1}^\ell+Q_{i+1}^\ell+1)+\nu(\eta_1,\dots,\eta_i)+\nu(\eta_{i+1},\dots,\eta_\ell)} \\
&\qquad \hat S^\nabla_i(\eta_1\otimes\dots\otimes\eta_i)(\theta_{P_1^i+1}^*(\alpha^{P_1^i+1})^*\sigma) \hat S^\nabla_{\ell-i}(\eta_{i+1}\otimes\dots\otimes\eta_\ell)(\theta_{P_{i+1}^\ell+1}^*(\beta^{P_{i+1}^\ell+1})^*\sigma).
\end{align*}
But notice that
\begin{align*}
\nu(\eta_1,\dots,\eta_\ell) & -\nu(\eta_1,\dots,\eta_i)-\nu(\eta_{i+1},\dots,\eta_\ell) =  \\
&= \sum_{j=1}^\ell (p_j-1)P_1^{j-1} -Q_1^\ell P_1^\ell-\sum_{j=1}^i(p_j-1)P_1^{j-1}-Q_1^iP_1^i \\
&\qquad -\sum_{j=i+1}^\ell (p_j-1)P_{i+1}^{j-1} +Q_{i+1}^\ell P_{i+1}^\ell \\
&= P_{i+1}^\ell P_1^i+Q_1^iP_{i+1}^\ell+Q_{i+1}^\ell P_1^i,
\end{align*}
so we can rewrite the previous expression as
\begin{align*}
(-1)^{1+P_1^\ell+Q_1^\ell +\nu(\eta_1\otimes\dots\otimes\eta_\ell)} & \left( 
\sum_{i=1}^k(-1)^i \hat S^\nabla_\ell(\eta_1\otimes\dots\otimes\eta_\ell)(\theta_{k-1}^*\delta_i^*\sigma)
\right. \\
&\qquad +\sum_{i=0}^\ell (-1)^{P_1^i+Q_1^iP_{i+1}^\ell} \hat S^\nabla_i(\eta_1\otimes\dots\otimes\eta_i)(\theta_{P_1^i+1}^*(\alpha^{P_1^i+1})^*\sigma) \\
&\qquad\qquad\qquad\qquad  \left.\hat S^\nabla_{\ell-i}(\eta_{i+1}\otimes\dots\otimes\eta_\ell)(\theta_{P_{i+1}^\ell+1}^*(\beta^{P_{i+1}^\ell+1})^*\sigma) \vphantom{\sum_{i=0}^\ell} \right).
\end{align*}
Comparing this to \eqref{eq:proof-Sigma-A-infty:rhs-corollary}, we get that $\Sigma^\nabla$ is an $A_\infty$-map if and only if
\begin{align*}
\sum_{i=1}^k (-1)^i &\left( \hat S^\nabla_\ell (\eta_1,\dots,\eta_\ell)(\delta_i^{1*}\theta_{k+1}^*\sigma) -\hat S^\nabla_\ell (\eta_1\otimes\dots\otimes\eta_\ell)(\delta_i^{0*}\theta_{k+1}^*\sigma) \right) \\
&=\sum_{i=1}^k(-1)^i \hat S^\nabla_\ell(\eta_1\otimes\dots\otimes\eta_\ell)(\theta_{k-1}^*\delta_i^*\sigma) \\
&\qquad +\sum_{i=0}^\ell (-1)^{P_1^i+Q_1^iP_{i+1}^\ell} \hat S^\nabla_i(\eta_1\otimes\dots\otimes\eta_i)(\theta_{P_1^i+1}^*(\alpha^{P_1^i+1})^*\sigma) \\
&\qquad\qquad\qquad\qquad \hat S^\nabla_{\ell-i}(\eta_{i+1}\otimes\dots\otimes\eta_\ell)(\theta_{P_{i+1}^\ell+1}^*(\beta^{P_{i+1}^\ell+1})^*\sigma).
\end{align*}

On the one hand, we see that \eqref{eq:good-cube-to-simplex:easy} already gives an equality for the first half of the sum, the one involving $\delta_1^i$. For the second part, notice that
\begin{align*}
\sum_{i=0}^\ell (-1)^{P_1^i+Q_1^iP_{i+1}^\ell} & \hat S^\nabla_i(\eta_1\otimes\dots\otimes\eta_i)(\theta_{P_1^i+1}^*(\alpha^{P_1^i+1})^*\sigma) \\
&\qquad\qquad\qquad \hat S^\nabla_{\ell-i}(\eta_{i+1}\otimes\dots\otimes\eta_\ell)(\theta_{P_{i+1}^\ell+1}^*(\beta^{P_{i+1}^\ell+1})^*\sigma)= \\
&=\sum_{i=1}^k \sum_{j=0}^\ell (-1)^{P_1^j+Q_1^jP_{j+1}^\ell} \hat S^\nabla_j(\eta_1\otimes\dots\otimes\eta_j)(\theta_i^*\alpha_i^*\sigma) \\
&\qquad\qquad\qquad\qquad \hat S^\nabla_{\ell-j}(\eta_{j+1}\otimes\dots\otimes\eta_\ell)(\theta_{k-i}^*\beta_{k-i}^*\sigma) \\
&=-\sum_{i=1}^k (-1)^i\sum_{P_1^j+1=i} (-1)^{Q_1^jP_{j+1}^\ell} \hat S^\nabla_j(\eta_1\otimes\dots\otimes\eta_j)(\theta_i^*\alpha_i^*\sigma) \\
&\qquad\qquad\qquad\qquad \hat S^\nabla_{\ell-j}(\eta_{j+1}\otimes\dots\otimes\eta_\ell)(\theta_{k-i}^*\beta_{k-i}^*\sigma).
\end{align*}
It is clear now that \eqref{eq:good-cube-to-simplex:hard} gives an equality for the summands involving $\delta_i^0$.
\end{proof}

We would like now to characterize more precisely the cube-to-simplex families which make \eqref{eq:good-cube-to-simplex:easy} and \eqref{eq:good-cube-to-simplex:hard} work. Condition \eqref{eq:good-cube-to-simplex:easy} is easy: we can simply ask $\delta_i\theta_k=\theta_{k+1}\delta_i^1$. Condition \eqref{eq:good-cube-to-simplex:hard} is much harder, as it somehow transforms a single integral into a composition of two integrals. To study \eqref{eq:good-cube-to-simplex:hard} we need to understand better the structure of the set of pinched cubes.

\begin{definition} \label{def:cubical-face-composition}
We define the \textbf{front-} and \textbf{back-face maps}, respectively, as
\begin{align*}
\alpha_i &:=(\delta_i^0)^{k-i}: I^i\longrightarrow I^k, \\
\beta_i &:= (\delta_0^1)^{k-i}: I^i\longrightarrow I^k.
\end{align*}
Explicitly, they are given by
\begin{align*}
\alpha_i(x) & =(x,0,\stackrel{(k-i)}{\dots},0), \\
\beta_i(x) &= (1,\stackrel{(k-i)}{\dots},1,x).
\end{align*}

We define the \textbf{cubical face-composition} map $\mu_i:I^{k-1}\to I^k$ as follows: consider $I^{k-1}=I\times I^{i-1}\times I^{k-i-1}$ and use coordinates $(t,x,y)$, with $t\in I$, $x\in I^{i-1}$ and $y\in I^{k-i-1}$. Then
\begin{equation}
\mu_i(t,x,y):= \begin{cases}
    \alpha_i(2t,x), & \text{if $0\leq t\leq \frac{1}{2}$}, \\
    \beta_{k-i}(2t-1,y), & \text{if $\frac{1}{2}< t\leq 1$}.
\end{cases} \label{eq:cubical_face_comp_maps}   
\end{equation}
\end{definition}

\begin{figure}
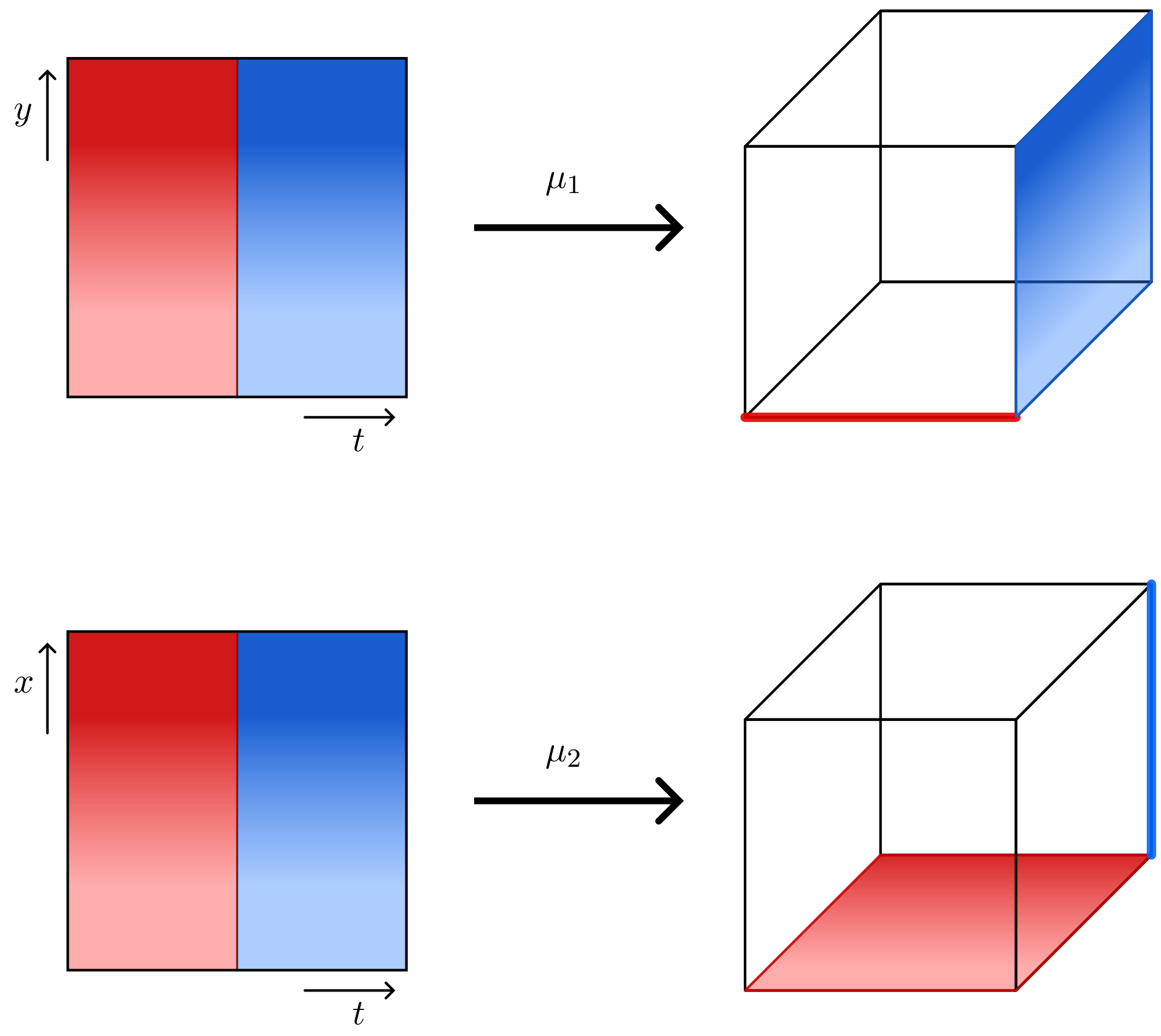

    \centering
    \includesvg[width=0.9\textwidth]{images/mu1_map.svg} \\ \vspace{3mm}
    \includesvg[width=0.9\textwidth]{images/mu2_map.svg}
    \caption{A graphical representation of the maps $\mu_1 : I \times I^0 \times I \times \to I^3$ and $\mu_2 : I \times I \times I^0 \to I^3$ of equation \eqref{eq:cubical_face_comp_maps}.}
    \label{fig:placeholder}
\end{figure}

Intuitively, $\mu_i$ gives us a way of composing a pinched $i$-cube and a pinched $(k-i)$-cube that are composable. More concretely, $\mu_i$ will give us a way of making more precise the conditions for $\theta$ such that \eqref{eq:good-cube-to-simplex:hard} holds. Figure \ref{fig:placeholder} shows a graphical representation of $\mu_i$ in low dimensions.

We only need one final ingredient.

\begin{definition}
A \textbf{simplicial face-composition family} is a family $\lambda=\{\lambda_i\}_i$ of smooth maps $\lambda_i:\Delta^i\times\Delta^{k-i}\to\Delta^k$ such that
\begin{align*}
\lambda_i(t,v_0) &=\alpha_i(t), \\
\lambda_i(v_i,s) &=\beta_{k-i}(s).
\end{align*}
\end{definition}

\begin{definition}
A cube-to-simplex family $\theta$ will be called \textbf{good} with respect to a simplicial face-composition family $\lambda$ if the diagrams
\begin{equation} \label{diag:good-cube-to-simplex:easy}
\begin{tikzcd}
	{I^{k-1}} & {I^k} \\
	{\Delta^{k-1}} & {\Delta^k}
	\arrow["{\delta_i^1}", from=1-1, to=1-2]
	\arrow["{\theta_{k-1}}"', from=1-1, to=2-1]
	\arrow["{\theta_k}", from=1-2, to=2-2]
	\arrow["{\delta_i}"', from=2-1, to=2-2]
\end{tikzcd}
\end{equation}
and
\begin{equation} \label{diag:good-cube-to-simplex:hard}
\begin{tikzcd}
	{I^{k-1}} & {I^k} & {\Delta^k} \\
	{I^i\times I^{k-i}} & &{\Delta^i\times\Delta^{k-i}} 
	\arrow["{\delta_i^0}", from=1-1, to=1-2]
	\arrow["{\mu_i}"', from=1-1, to=2-1]
	\arrow["{\theta_k}", from=1-2, to=1-3]
	\arrow["{\theta_i\times\theta_{k-i}}"', from=2-1, to=2-3]
    \arrow["{\lambda_i}"', from=2-3, to=1-3]
\end{tikzcd}
\end{equation}
commute for $0<i<k$, possibly up to reparametrization, by which we mean that $\theta_k\delta_i^1=\delta_i\theta_{k-i}\varphi_i^1$ and $\theta_k\delta_i^0=\lambda_i(\theta_i\times\theta_{k-i})\mu_i\varphi_i^0$, where $\varphi_i^m:I^{k-1}\to I^{k-1}$ are reparametrizations of the kind $\varphi_i^m=\chi_i^m\times\id$, for $\chi_i^m:I\to I$ piecewise smooth with $\chi_i^m(0)=0$ and $\chi_i^m(1)=1$.
\end{definition}

\begin{theorem}
Good cube-to-simplex families satisfy equations \eqref{eq:good-cube-to-simplex:easy} and \eqref{eq:good-cube-to-simplex:hard}.
\end{theorem}

\begin{proof}
Clearly \eqref{eq:good-cube-to-simplex:easy} holds, by \eqref{diag:good-cube-to-simplex:easy}. Regarding \eqref{eq:good-cube-to-simplex:hard}, we will prove it for $\ell=2$ and the general case will easily be seen to follow along the same lines.

We use coordinates $(t,x,y)$ in $I^{k-1}$ as in Definition \ref{def:cubical-face-composition}. First of all, if $\xi$ is any homogeneous form on $I^k$, then
\[
\mu_i^*\xi\left(\derpar{}{t},\derpar{}{x_{s_1}},\dots, \derpar{}{x_{s_q}},\derpar{}{y_{s_{q+1}}},\dots, \derpar{}{y_{s_{\abs{\xi}-1}}}\right)(t,x,y)
\]
can only be non-zero if $q=\abs{\xi}-1$, in which case it equals
\[
\left\{
\begin{array}{ll}
2\alpha_i^*\xi(\derpar{}{t},\derpar{}{x_{s_1}},\dots, \derpar{}{x_{s_{\abs{\xi}-1}}})(2t,x), & 0\leq t\leq \frac{1}{2} \\
0, & \frac{1}{2}< t\leq 1,
\end{array}
\right. 
\]
or $q=0$, in which case it equals
\[
\left\{
\begin{array}{ll}
0, & 0\leq t\leq \frac{1}{2} \\
2\beta_{k-i}^*\xi(\derpar{}{t},\derpar{}{y_{s_1}},\dots, \derpar{}{y_{s_{\abs{\xi}-1}}})(2t-1,y), & \frac{1}{2}< t\leq 1.
\end{array}
\right.
\]

Let $k=p_1+p_2$ and $\sigma\in\Pi_\infty(A)_k$, and let $\phi:=\delta_i^{0*}\theta_k^*\sigma \in I_\infty(A)_{k-1}$. Consider
\[
\xi_j :=(\theta_i\times\theta_{k-i})^*\lambda_i^*\sigma^*\eta_j 
\]
Then, by the previous computation for $\mu_i^*\xi_j$ and using the explicit formula \eqref{eq:iterated-integral-explicit}, we get that
\[
(-1)^{q_1(p_2-1)}\int_\phi (\mu_i^*\xi_1\otimes \mu_i^*\xi_2)
\]
can only be non-zero when $i=p_1$. In such case, if $p_1\neq p_2$, the only surviving term is
\begin{align*}
\int_{I^{k-i-1}}\int_{I^{i-1}} & \int_0^1\int_u^1 \tau_{0\leftarrow u,x,y} \mu_i^* \xi_1\left(\derpar{}{t},\derpar{}{x_1},\dots, \derpar{}{x_{i-1}}\right)(u,x,y)\tau_{u\leftarrow t,x,y} \\
&\qquad \mu_i^*\xi_2\left(\derpar{}{t},\derpar{}{y_1},\dots, \derpar{}{y_{k-i-1}}\right)(t,x,y) \tau_{t\leftarrow 1,x,y} dtdudxdy.
\end{align*}
By the previous computation, the integrand is only non-zero when $0\leq u\leq \frac{1}{2}$ and when $\frac{1}{2}\leq t\leq 1$, so this equals
\begin{align*}
\int_{I^{k-i-1}}&\int_{I^{i-1}}\int_0^{1/2}\int_{1/2}^1 \tau_{0\leftarrow 2u,x,0}2\alpha_i^* \xi_1\left(\derpar{}{t},\derpar{}{x_1},\dots, \derpar{}{x_{i-1}}\right)(2u,x)\tau_{2u\leftarrow 1,x,0} \\
&\qquad \tau_{0\leftarrow 2t-1,1,y}2\beta_{k-i}^*\xi_2\left(\derpar{}{t},\derpar{}{y_1},\dots, \derpar{}{y_{k-i-1}}\right)(2t-1,y)\tau_{2t-1\leftarrow 1,1,y} dtdudxdy = \\
&=\int_{I^{i-1}}\int_0^1\tau_{0\leftarrow u,x}\alpha_i^*\xi_1\left(\derpar{}{t},\derpar{}{x_1},\dots, \derpar{}{x_{i-1}}\right)(u,x)\tau_{u\leftarrow 1,x}dudx \\
&\qquad \int_{I^{k-i-1}}\int_0^1 \tau_{0\leftarrow t,y}\beta_{k-i}^*\xi_2\left(\derpar{}{t}, \derpar{}{y_1},\dots, \derpar{}{y_{k-i-1}}\right) (t,y)\tau_{t\leftarrow 1,y}dtdy \\
&=\left( \int_{\phi_i}\alpha_i^*\xi_1\right)\left(\int_{\phi_{k-i}}\beta_{k-i}^*\xi_2\right),
\end{align*}
where $\phi_i:=\alpha_i^*(\theta_i\times\theta_{k-i})^*\lambda_i^*\sigma\in I_\infty(A)_i$ and $\phi_{k-i}:=\beta_{k-i}^*(\theta_i\times\theta_{k-i})^*\lambda_i^*\sigma\in I_\infty(A)_{k-i}$. This last computation made sense because
\[
\sigma\circ\lambda_i\circ (\theta_i\times\theta_{k-i})(\alpha_i(1,x))=\sigma(v_i,v_0)=\sigma\circ\lambda_i\circ (\theta_i\times\theta_{k-i})(\beta_{k-i}(0,y)).
\]
The only case where we could have an additional term is when $p_1=p_2$. In such case, a priori we would also have the terms
\begin{align*}
\int_{I^{k-i-1}}\int_{I^{i-1}} & \int_0^1\int_u^1\tau_{0\leftarrow u,x,y}\mu_i^* \xi_1\left(\derpar{}{t},\derpar{}{y_1},\dots, \derpar{}{y_{i-1}}\right)(u,x,y)\tau_{u\leftarrow t,x,y} \\
&\qquad \mu_i^*\xi_2\left(\derpar{}{t},\derpar{}{x_1},\dots, \derpar{}{x_{k-i-1}}\right)(t,x,y)\tau_{t\leftarrow 1,x,y} dtdudxdy.
\end{align*}
But, since the integrand is only non-zero when $0\leq t\leq \frac{1}{2}$ and $\frac{1}{2}\leq u\leq 1$, an immediate inspection of the limits of integration gives that this contribution vanishes.

We finally conclude that, using \eqref{diag:good-cube-to-simplex:hard},
\begin{align*}
\hat S^\nabla_2(\eta_1\otimes\eta_2)(\delta_i^{0*}\theta_k^*\sigma) &=(-1)^{q_1(p_2-1)}\int_\phi(\mu_i^*\xi_1\otimes\mu_i^*\xi_2) \\
&=\left(\int_{\phi_i}\alpha_i^*\xi_1\right)\left(\int_{\phi_{k-i}}\beta_{k-i}^*\xi_2\right) \\
&=(-1)^{Q_1^1P_2^2}\hat S^\nabla_1(\eta_1)(\theta_i^*\alpha_i^*\sigma)\hat S^\nabla_1(\eta_2)(\theta_{k-i}^*\beta_{k-i}^*\sigma),
\end{align*}
where we used that
\[
\phi_i(t,x)=\sigma\circ\lambda_i\circ(\theta_i\times\theta_{k-i})(\alpha_i(t,x))=\sigma(\alpha_i\theta_i(t,x)),
\]
and similarly for $\phi_{k-i}$.
\end{proof}

\begin{corollary}
With respect to a good cube-to-simplex family, $\Sigma^\nabla$ defines an $A_\infty$-map.
\end{corollary}

The only thing that remains is to see is that this theory is non-empty, i.e., that we can find at least one good cube-to-simplex family. The families we use coincide with the ones from \cite{antweiler}. These are different from the ones used in \cite{igusa} and \cite{arias-abad-schatz}.

\begin{definition} \label{def:face-composition-family-example}
We define the map $\lambda_i:\Delta^i\times\Delta^{k-i}\to\Delta^k$ by
\[
\lambda_i(t,s):=((1-s_1)t+s_1v_i,s).
\]
\end{definition}

\begin{lemma}
The family $\lambda=\{\lambda_i\}_i$ from Definition \ref{def:face-composition-family-example} is a simplicial face-composition family.
\end{lemma}

\begin{proof}
First of all we check that indeed $\lambda_i$ lands in $\Delta^k$. Explicitly we have that
\[
\lambda_i(t,s)=((1-s_1)t_1+s_1,\dots,(1-s_1)t_i+s_1,s_1,\dots,s_{k-i}),
\]
which clearly lies in $\Delta^k$. Checking the condition for it being a simplicial face-composition map is a straightforward computation:
\begin{align*}
\lambda_i(t,v_0) &=(t,v_0)=\alpha_i(t), \\
\lambda_i(v_i,s) &=(v_i,s)=\beta_{k-i}(s). \qedhere
\end{align*}
\end{proof}

\begin{definition} \label{def:cube-to-simplex-family-example}
We define the map $\theta_k:I^k\to\Delta^k$ as the composition
\[
I^k\stackrel{\tilde\theta_k}{\longrightarrow} I^k\stackrel{p_k}{\longrightarrow} \Delta^k,
\]
where $\tilde\theta_k(t,x_1,\dots,x_{k-1})$ is given by the evaluation at $t$ of the piecewise affine path in $I^k$ that, in $k$ steps, passes through the points
\begin{align*}
(0,x_1,\dots,x_{k-1})&\to (1,x_1,\dots,x_{k-1})\to (1,1,x_2,\dots,x_{k-1})\to\dots \\
\dots &\to (1,\dots,1,x_{k-1})\to (1,\dots,1),
\end{align*}
and $p_k$ is given by
\[
p_k(x_1,\dots,x_k):=(t_1,\dots, t_k),\qquad t_i:=x_1\dots x_i.
\]
\end{definition}

Explicitly, $\tilde\theta_k$ is given by
\[
\tilde\theta_k\left(\frac{i}{k},x_1,\dots,x_{k-1}\right) =\left\{
\begin{array}{ll}
    (0,x_1,\dots,x_{k-1}), & i=0 \\
    (1,\dots, 1, x_i,\dots, x_{k-1}), & i=1,\dots, k
\end{array}
\right.
\]
and interpolating affinely between these.

\begin{figure}
    \centering
    \includesvg[width=0.9\textwidth]{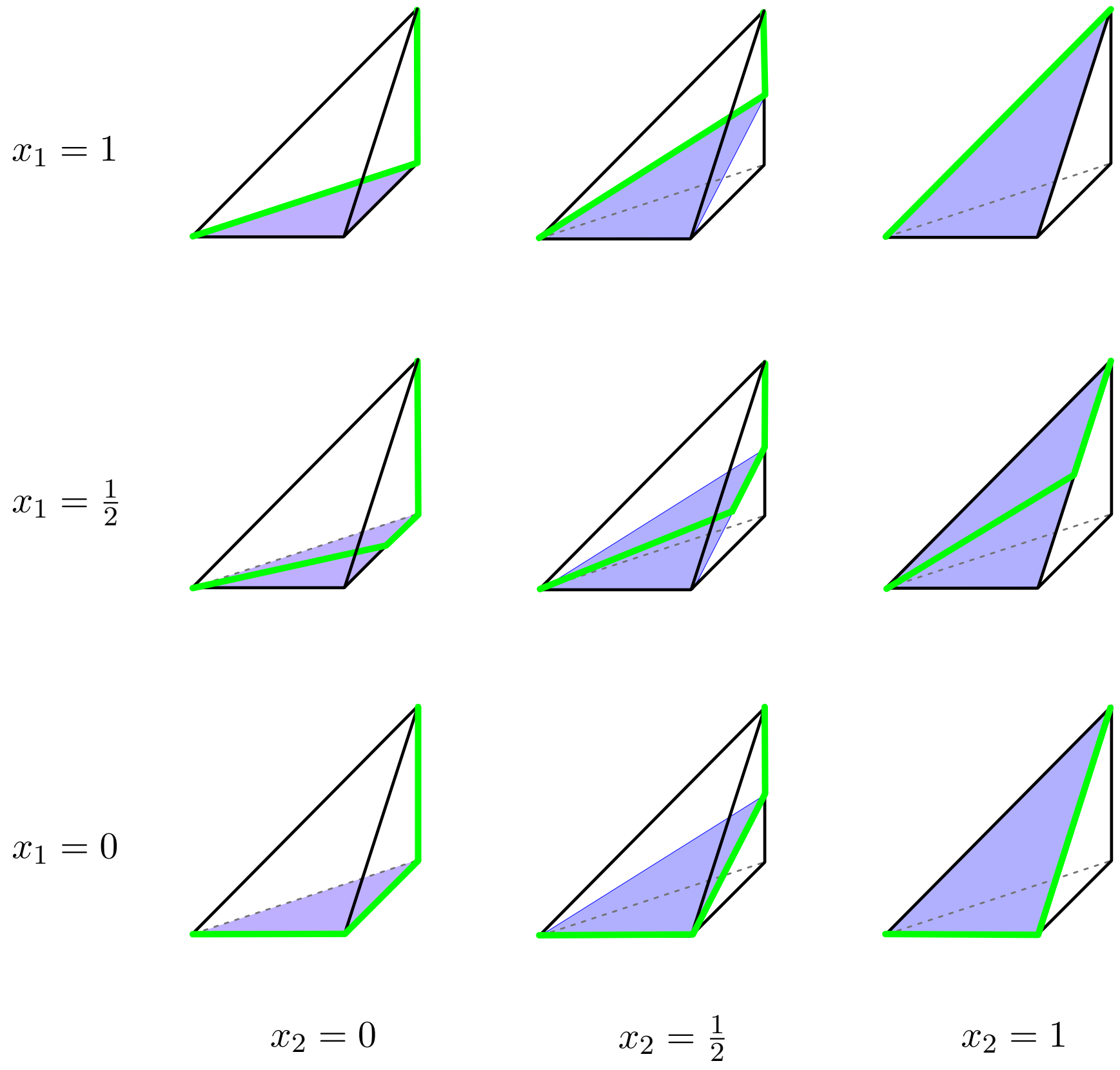}
    \caption{Illustration of the paths $t \mapsto \theta_3(t, x_1, x_2)$ from Definition \ref{def:cube-to-simplex-family-example} for the specified values of $x_1$ and $x_2$. Here, the blue surface is the space of constant $x_2$.}
    \label{fig:path_families}
\end{figure}

Figure \ref{fig:path_families} shows a schematic depiction of these cube-to-simplex families, which one can compare with the families in \cite{igusa} depicted between lemma 4.7 and lemma 4.8.

\begin{lemma}
The cube-to-simplex family $\theta$ from Definition \ref{def:cube-to-simplex-family-example} makes the diagram \eqref{diag:good-cube-to-simplex:easy} commute up to reparametrization.
\end{lemma}

\begin{proof}
Notice first that
\begin{align*}
p_k\delta_i^1(x_1,\dots,x_{k-1}) &=p_k(x_1,\dots,x_i,1,x_{i+1},\dots, x_{k-1}) \\
&=(x_1,x_1x_2,\dots, x_1\dots x_i,x_1\dots x_i,\dots, x_1\dots x_{k-1}) \\
&=\delta_ip_{k-1}(x_1,\dots, x_{k-1}).
\end{align*}
Hence, it is enough to show that $\tilde\theta_k\delta_i^1=\delta_i^1\tilde\theta_{k-1}$ (up to reparametrization).

On the one hand, $\tilde\theta_k\delta_i^1(t,x_1,\dots, x_{k-2})$ is the evaluation at $t$ of  the piecewise affine path that, in $k$ steps, passes through the points
\begin{align*}
(0,x_1,\dots,x_{i-1},1,x_i,\dots, x_{k-2}) &\to (1,x_1,\dots, x_{i-1},1,x_i,\dots, x_{k-2}) \\
&\to (1,1,x_2,\dots, x_{i-1},1,x_i,\dots, x_{k-2})\to\dots \\
&\to (1,\dots, 1,x_{i-1},1,x_i,\dots, x_{k-2}) \\
&\to (1,\dots, 1,1,1,x_i,\dots, x_{k-2}) \\
&\to(1,\dots, 1,1,1,x_i,\dots, x_{k-2}) \\
&\to (1,\dots, 1,x_{i+1},\dots, x_{k-2})\to\dots \\
&\to (1,\dots, 1),
\end{align*}
whereas $\delta_i^1\tilde\theta_{k-1}(t,x_1,\dots, x_{k-2})$ is the evaluation at $t$ of the piecewise affine path that, in $k-1$ steps, passes through the points
\begin{align*}
(0,x_1,\dots, x_{i-1},1,x_i,\dots, x_{k-2}) &\to (1,x_1,\dots, x_{i-1},1,x_i,\dots, x_{k-2}) \\
&\to (1,1,x_2,\dots, x_{i-1},1,x_i,\dots, x_{k-2})\to\dots \\
&\to (1,\dots, 1,x_{i-1},1,x_i,\dots, x_{k-2}) \\
&\to (1,\dots, 1,1,1,x_i,\dots, x_{k-2}) \\
&\to (1,\dots,1, x_{i+1},\dots, x_{k-2}) \to\dots \\
&\to(1,\dots, 1).
\end{align*}
Is is clear now that the reparametrization
\[
\chi_i^1(t):=\left\{
\begin{array}{ll}
   \frac{kt}{k-1},  & t\in[0,\frac{i}{k}]  \\
   \frac{i}{k-1},  & t\in[\frac{i}{k},\frac{i+1}{k}] \\ 
   \frac{kt-1}{k-1}, & t\in[\frac{i+1}{k},1]
\end{array}
\right.
\]
will do the job.
\end{proof}

\begin{lemma}
The cube-to-simplex family $\theta$ from Definition \ref{def:cube-to-simplex-family-example} makes the diagram \eqref{diag:good-cube-to-simplex:hard} commute up to reparametrization, with respect to the simplicial face-composition family $\lambda$ from Definition \ref{def:face-composition-family-example}.
\end{lemma}

\begin{proof}
Throughout the proof we use coordinates $(t,x,y)$ on $I^{k-1}$, as in Definition \ref{def:cubical-face-composition}.
Let 
\begin{align*}
I^{k-1}_{\leq s} &:=\{(t,x,y)\in I^{k-1}:t\leq s\}, \\
I^{k-1}_{\geq s} &:=\{(t,x,y)\in I^{k-1}:t\geq s\}.
\end{align*}

First of all, notice that
\[
(\tilde\theta_i\times\tilde\theta_{k-i})\mu_i(t,x,y) = \left\{
\begin{array}{ll}
    \alpha_i(\tilde\theta_i(2t,x)), & 0\leq t\leq \frac{1}{2} \\
    \beta_{k-i}(\tilde\theta_{k-i}(2t-1,y), & \frac{1}{2}<t\leq 1
\end{array}
\right.
\]
so $(\tilde\theta_i\times\tilde\theta_{k-i})\mu_i(I^{k-1}_{\leq 1/2})\subseteq \alpha_i(I^i)$ and $(\tilde\theta_i\times\tilde\theta_{k-i})\mu_i(I^{k-1}_{\geq 1/2})\subseteq \beta_{k-i}(I^{k-i})$.

Observe as well that
\begin{align*}
\lambda_i(p_i\times p_{k-i})\alpha_i(t,x) &= \lambda_i(p_i(t,x),v_0) =\alpha_i p_i(t,x) \\
&=(p_i(t,x),0,\dots,0)=p_k\alpha_i(t,x), \\
\lambda_i(p_i\times p_{k-i})\beta_{k-i}(t,y) &= \lambda_i(v_i,p_{k-i}(t,y)) =\beta_{k-i} p_{k-i}(t,y) \\
&=(1,\dots,1,p_{k-i}(t,y))=p_k\beta_{k-i}(t,y),
\end{align*}
which means that, as maps $\alpha_i(I^i) \cup \beta_{k-i}(I^{k-i})\to \Delta^k$, we have that $\lambda_i(p_i\times p_{k-i})=p_k$.

On the other hand, the map $\tilde\theta_k\delta_i^0$ can be described as follows: $\tilde\theta_k\delta_i^0(t,x,y)$ is the evaluation at $t$ of  the piecewise affine path that, in $k$ steps, passes through the points
\begin{align}
(0,x,0,y) &\to (1,x,0,y) \label{eq:thetak-deltai0} \\
&\to (1,1,x_2,\dots, x_{i-1},0,y)\to\dots \nonumber\\
&\to (1,\dots, 1,x_{i-1},0,y) \nonumber\\
&\to (1,\dots, 1,1,0,y) \nonumber\\
&\to(1,\dots, 1,1,1,y) \nonumber\\
&\to (1,\dots, 1,y_2,\dots, y_{k-i-1})\to\dots \nonumber\\
&\to (1,\dots, 1). \nonumber
\end{align}
Hence, we see that $\tilde\theta_k\delta_i^0(I^{k-1}_{\leq i/k})\subseteq I^i\times \{0\}\times I^{k-i-1}$, and $\tilde\theta_k\delta_i^0(I^{k-1}_{\geq i/k})\subseteq \beta_{k-i}(I^{k-i})$.

Consider the reparametrization $\varphi_i^0=\chi_i^0\times\id:I^{k-1}\to I^{k-1}$ given by
\[
\chi_i^0(t)=\left\{
\begin{array}{ll}
    \frac{2it}{k} & t\in[0,\frac{1}{2}] \\
    \frac{2(k-i)t}{k}+\frac{2i-k}{k}, & t\in[\frac{1}{2},1].
\end{array}
\right.
\]
Then $\varphi_i^0(I^{k-1}_{\leq 1/2})=I^{k-1}_{\leq i/k}$ and $\varphi_i^0(I^{k-1}_{\geq 1/2})= I^{k-1}_{\geq i/k}$. We will prove the commutativity of \eqref{diag:good-cube-to-simplex:hard} in $I^{k-1}_{\leq 1/2}$ and $I^{k-1}_{\geq 1/2}$ separately, up to the reparametrization $\varphi_i^0$.

For $I^{k-1}_{\leq 1/2}$, consider the map $m:I^i\times\{0\}\times I^{k-i-1}\to \alpha_i(I^i)$ given by $m(t,x,0,y):=(t,x,0,0)$, and notice that $p_km(t,x,0,y)=p_k(t,x,0,y)$. We can describe the situation by the following diagram:
\[\begin{tikzcd}
	{I^{k-1}_{\leq 1/2}} & {\alpha_i(I^i)} & {\alpha_i(I^i)} & {\Delta^i\times\Delta^k} \\
	{I^{k-1}_{\leq i/k}} & {I^k} & {I^i\times\{0\}\times I^{k-i-1}} & {\Delta^k.}
	\arrow["{\mu_i}", from=1-1, to=1-2]
	\arrow["{\varphi_i^0}"', from=1-1, to=2-1]
	\arrow["{\tilde\theta_i\times\tilde\theta_{k-i}}", from=1-2, to=1-3]
	\arrow["{p_i\times p_{k-i}}", from=1-3, to=1-4]
	\arrow["{p_k}"', from=1-3, to=2-4]
	\arrow["{\lambda_i}", from=1-4, to=2-4]
    \arrow["{\delta_i^0}"', from=2-1, to=2-2]
	\arrow["{\tilde\theta_k}"', from=2-2, to=2-3]
	\arrow["{m}"', from=2-3, to=1-3]
    \arrow["{p_k}"',from=2-3, to=2-4]
\end{tikzcd}\]
The two right-most triangles commute, so to prove that the outer diagram commutes it suffices to prove that the left rectangle commutes.

On the one hand, $(\tilde\theta_i\times\tilde\theta_{k-i})\mu_i(t,x,y)$, for $t\leq \frac{1}{2}$, is the evaluation at $t$ of the piecewise affine path that, in $i$ steps and in total time $\frac{1}{2}$, passes through the points
\begin{align*}
(0,x,0,\dots, 0) &\to (1,x,0,\dots, 0) \\
&\to (1,1,x_2,\dots, x_{i-1},0,\dots, 0)\to\dots \\
&\to (1,\dots, 1,x_{i-1},0,\dots, 0) \\
&\to (1,\dots, 1,1,0,\dots, 0),
\end{align*}
whereas $m\tilde\theta_k\delta_i^0(t,x,y)$, for $t\leq \frac{i}{k}$, is the evaluation at $t$ of the piecewise affine path that, in $i$ steps and in total time $\frac{i}{k}$, passes through the points (see \eqref{eq:thetak-deltai0})
\begin{align*}
(0,x,0,\dots, 0) &\to (1,x_1,\dots, x_{i-1},0,\dots, 0) \\
&\to (1,1,x_2,\dots, x_{i-1},0,\dots, 0)\to\dots \\
&\to (1,\dots, 1,x_{i-1},0,\dots, 0) \\
&\to (1,\dots, 1,1,0,\dots, 0).
\end{align*}
It is clear now that the reparametrization $\varphi_i^0$ indeed makes the diagram commute.

Turning now to $I^{k-1}_{\geq 1/2}$, the situation is easier, since we have the diagram
\[\begin{tikzcd}
	{I^{k-1}_{\geq 1/2}} & {\beta_{k-i}(I^{k-i})} & {\beta_{k-i}(I^{k-i})} & {\Delta^i\times\Delta^k} \\
	{I^{k-1}_{\geq i/k}} & {I^k} & & {\Delta^k.}
	\arrow["{\mu_i}", from=1-1, to=1-2]
	\arrow["{\varphi_i^0}"', from=1-1, to=2-1]
	\arrow["{\tilde\theta_i\times\tilde\theta_{k-i}}", from=1-2, to=1-3]
	\arrow["{p_i\times p_{k-i}}", from=1-3, to=1-4]
	\arrow["{p_k}"', from=1-3, to=2-4]
	\arrow["{\lambda_i}", from=1-4, to=2-4]
    \arrow["{\delta_i^0}"', from=2-1, to=2-2]
	\arrow["{\tilde\theta_k}"', from=2-2, to=1-3]
\end{tikzcd}\]
The right triangle commutes, so it is enough to show that the left-hand side commutes as well.

On the one hand, $(\tilde\theta_i\times\tilde\theta_{k-i})\mu_i(t,x,y)$, for $t\geq \frac{1}{2}$, is the evaluation at $t$ of the piecewise affine path that, in $k-i$ steps and in total time $\frac{1}{2}$, passes through the points
\begin{align*}
(1,\dots 1,0, y) &\to (1,\dots 1,1, y) \\
&\to (1,\dots 1,1, 1, y_2,\dots, y_{k-i-1})\to\dots \\
&\to (1,\dots, 1),
\end{align*}
whereas $\tilde\theta_k\delta_i^0(t,x,y)$, for $t\geq \frac{i}{k}$, is the evaluation at $t$ of the piecewise affine path that, in $k-i$ steps and in total time $\frac{k-i}{k}$, passes through the points (see \eqref{eq:thetak-deltai0} again)
\begin{align*}
(1,\dots, 1,0,y) &\to(1,\dots, 1,1,y)\\
&\to (1,\dots, 1,y_2,\dots, y_{k-i-1})\to\dots \\
&\to (1,\dots, 1). 
\end{align*}
It is clear again that the reparametrization $\varphi_i^0$ indeed makes the diagram commute.
\end{proof}

These two lemmas put together finally prove the following.

\begin{proposition} \label{prop:good-cube-to-simplex-family-example}
The  cube-to-simplex family $\theta$ from Definition \ref{def:cube-to-simplex-family-example} is a good cube-to-simplex family with respect to the simplicial face-composition family $\lambda$ from Definition \ref{def:face-composition-family-example}.
\end{proposition}

\bigskip
\section{The Riemann--Hilbert functor} \label{sect:riemann-hilbert-functor}

\subsection{Curved differential graded categories}
The following is a standard categorification of the material of section \ref{sect:curved-dg-algebras}.
\begin{definition}
A \textbf{curved differential graded category} (curved dg-category for short) is a category $C$ enriched over real graded vector spaces together with the following data:
\begin{enumerate}[label=(\roman*)]
    \item for all $E_0,E_1 \in C$, a linear degree 1 map $D : \Hom(E_0, E_1) \to \Hom(E_0, E_1)$, referred to as the \textbf{connection},
    \item for all $E \in C$, a degree 2 endomorphism $K_E \in \End^2(E)$, referred to as the \textbf{curvature},
\end{enumerate}
such that
\begin{enumerate}[label=(\roman*)]
    \item for all $E\in C$, the identity $1_E$ and the curvature $K_E$ are $D$-closed, i.e., $D1_E=0$ and $DK_E=0$,
    \item the composition 
    \[
    \Hom(E_1,E_0)\times \Hom(E_2,E_1)\to \Hom(E_2,E_0),
    \]
    for $E_i\in C$, is graded bilinear and satisfies the Leibniz rule: if $\Phi_i\in\Hom(E_i,E_{i-1})$, then
    \[
    D(\Phi_1\Phi_2) = (D\Phi_1)\Phi_2+(-1)^{\abs{\Phi_1}}\Phi_1(D\Phi_2),
    \]
    \item for all $\Phi \in \Hom(E_1, E_0)$, we have that
    \[
    D^2\Phi = K_{E_0} \Phi - \Phi K_{E_1}.
    \]
\end{enumerate}

A \textbf{dg-category} is a curved dg-category all of whose curvature morphisms are zero. In such case, the connections on the sets of morphisms we call \textbf{differentials} and we write them as $d$. For a curved dg-category $C$, its \textbf{dg-category of flat objects} $\Flat(C)$ is the full subcategory of $C$ whose objects are $E\in C$ with $K_E=0$.
\end{definition}

The theory of curved dg-categories is part of the wider theory of $A_\infty$-categories, the categorification of $A_\infty$-algebras.

\begin{definition}
An \textbf{$A_\infty$-category} $A$ consists of a collection $\textrm{Obj}(A)$ called objects and for each pair of objects $E_1, E_2 \in \textrm{Obj}(A)$ a graded vector space denoted $\Hom(E_2, E_1)$ together with a collection of degree 1 linear maps
\[
b_\ell: \bigotimes_{i=1}^\ell \Hom(E_i,E_{i-1})[1]\to \Hom(E_\ell,E_0)[1]
\]
for $\ell\geq 0$ and $E_i\in A$, such that
\[
\sum_{i+j+k=\ell} b_{i+k+1}\circ (\id^{\otimes i}\otimes b_j\otimes \id^{\otimes k})=0.
\]

A functor $\psi:A_0\to A_1$ between $A_\infty$-categories is a map $\psi$ at the level of objects and a collection of degree 0 maps
\[
\psi_\ell:\bigotimes_{i=1}^\ell \Hom(E_i,E_{i-1})[1]\to \Hom(\psi(E_\ell),\psi(E_0))[1]
\]
for $\ell\geq 0$ and $E_i\in A_0$, such that
\[
\sum_{k\geq 0}\sum_{i_1+\dots+i_k=\ell} b^1_k \circ (\psi_{i_1}\otimes\dots\otimes \psi_{i_k}) =\sum_{i+j+k=\ell} \psi_{i+k+1}\circ (\id^{\otimes i}\otimes b^0_j\otimes \id^{\otimes k}).
\]
\end{definition}

From now on the collection $\mathrm{Obj}(A)$ we will just denote by $A$ as well.

\begin{example}
An $A_\infty$-category with a single object is an $A_\infty$-algebra, and an $A_\infty$-category with $b_\ell=0$ for $\ell\geq 3$ is a curved dg-category.
\end{example}

As in Section \ref{sect:curved-dg-algebras}, for an $A_\infty$-functor $\psi:A_0\to A_1$ we will use $\hat\psi_\ell$ to denote the induced degree $1-\ell$ map $\bigotimes_{i=1}^\ell \Hom(E_i,E_{i-1})\to \Hom(\psi(E_\ell),\psi(E_0))$. For the case of curved dg-categories, note that the formulas for an $A_\infty$-functor are formally the same as those of an $A_\infty$-map between curved dgas (cf. Corollary \ref{cor:A-infty-map-for-curved-dgas}).

We now introduce the two curved dg-categories that we will be working with.

\begin{definition}
Let $(\Ascr\to M,Q)$ be a higher Lie algebroid. We define the \textbf{category of connections $\HConn(\Ascr)$ on $\Ascr$} as follows:
\begin{enumerate}
    \item objects are pairs $(E,D)$, where $E\to M$ is a graded vector bundle and $D$ is a higher $\Ascr$-connection on $E$,
    \item $\Hom((E_1,D_1),(E_0,D_0)):=\Hom(\Escr_1,\Escr_0)\cong C(\Ascr,\Hom(E_1,E_0))$, with the induced $\Ascr$-connection from Lemma \ref{lem:higher-connection-on-homomorphisms},
    \item if $(E,D)\in\HConn(\Ascr)$, then $K_D:=D^2$.
\end{enumerate}
We also define the \textbf{category of ruths of $\Ascr$} by $\HRep(\Ascr):=\Flat(\HConn(\Ascr))$.

Let now $G$ be a higher Lie groupoid. We define the \textbf{category of connections $\HConn(G)$ on $G$} as follows:
\begin{enumerate}
    \item objects are pairs $(E,D)$, where $E\to M$ is a graded vector bundle and $D$ is a $G$-connection on $E$,
    \item $\Hom((E_1,D_1),(E_0,D_0)):=\Hom_{C(G)}(C(G,E_1),C(G,E_0))\cong C_G(\Hom(E_1,E_0))$, with the induced $G$-connection, analogous to Lemma \ref{lem:higher-connection-on-homomorphisms},
    \item if $(E,D)\in\HConn(G)$, then $K_D:=D^2$.
\end{enumerate}
We also define the \textbf{category of ruths of $G$} by $\HRep(G):=\Flat(\HConn(G))$.
\end{definition}

\subsection{The Riemann--Hilbert functor for Lie algebroids}

In this section, we fix a good cube-to-simplex family $\theta$, for instance the one from Proposition \ref{prop:good-cube-to-simplex-family-example}.

This section will hinge on the following crucial gauge equivariance of $\Sigma^\nabla$.

\begin{proposition} \label{prop:gauge-equivariance}
Let $\alpha\in\Omega^1(A,\End^0(E))$. Then
\[
\hat S_\ell^{\nabla+\alpha}(\eta_1\otimes\dots\otimes\eta_\ell)=\sum_{m_i\geq 0}\hat S^\nabla_{m_0+\dots+m_\ell+\ell}(\alpha^{m_0}\otimes\eta_1\otimes \alpha^{m_1}\otimes\dots\otimes \alpha^{\otimes m_{\ell-1}}\otimes\eta_\ell\otimes\alpha^{m_\ell}).
\]
In particular, the same holds for $\hat\Sigma^\nabla$.
\end{proposition}

\begin{proof}
By Example \ref{ex:parallel-transport-as-iterated-integral}, we have that $\tau^{\nabla+\alpha}_{a^{-1}}=\tau^\nabla_{a^{-1}}+\sum_{\ell\geq 1}\int_a\alpha^{\otimes\ell}$ for any $A$-path $a$. Then, if $\phi\in I_\infty(A)$ and $\eta_i\in \phi^*\Omega(A,\End(E))$ are homogeneous, one can directly see from \eqref{eq:iterated-integral-explicit} that
\[
\int^{\nabla+\alpha}_\phi (\eta_1\otimes\dots\otimes\eta_\ell) =\sum_{m_i\geq 0}\int^\nabla_\phi(\alpha^{m_0}\otimes\eta_1\otimes \alpha^{m_1}\otimes\dots\otimes \alpha^{\otimes m_\ell}\otimes\eta_\ell\otimes\alpha^{m_\ell}).
\]
The claim for $\hat\Sigma^\nabla$ follows from that fact that
\[
\nu(\eta_1,\dots,\eta_i,\alpha,\eta_{i+1},\dots, \eta_a)=\nu(\eta_1,\dots,\eta_a). \qedhere
\]
\end{proof}

It is now that we make use of Remark \ref{rem:integral-different-bundles}.

\begin{definition}\label{def:algebroid_RH_functor}
We define the \textbf{Riemann--Hilbert functor} as the $A_\infty$-functor
\[
\RH:\HConn(A)\to \HConn(\Pi_\infty(A))
\]
defined as follows:
\begin{enumerate}
    \item If $(E,D)\in \HConn(A)$, we set
    \[
    \RH(D):=\hat\Sigma^\nabla(\omega)=\sum_{m\geq 0}\hat \Sigma^\nabla_m(\omega^{\otimes m}),
    \]
    where $\nabla$ is any degree-preserving $A$-connection on $E$ and $\omega:=D-d^\nabla\in\Omega(A,\End(E))^1$.

    \item If $(E_i,D_i)\in\HConn(A)$, for $i=0,\dots,\ell$, and $\eta_i\in\Omega(A,\Hom(E_i,E_{i-1}))$, then
    \begin{align*}
    \widehat{\RH}_\ell(\eta_1 & \otimes\dots\otimes\eta_\ell):= \hat\Sigma^{\nabla^0,\dots,\nabla^\ell}(\omega_0\otimes\eta_1\otimes\omega_1\otimes\dots\otimes\omega_{\ell-1}\otimes\eta_\ell\otimes\omega_\ell) \\
    &:=\sum_{m_i\geq 0} \hat\Sigma^{\nabla^0,\dots,\nabla^\ell}_{m_0+\dots+m_\ell+\ell} (\omega_0^{\otimes m_0}\otimes\eta_1\otimes \omega_1^{\otimes m_1}\otimes\dots\otimes\omega_{\ell-1}^{\otimes m_{\ell-1}}\otimes\eta_\ell\otimes\omega_\ell^{\otimes m_\ell}),
    \end{align*}
    where $\nabla^i$ is any degree-preserving $A$-connection on $E_i$ and $\omega_i:=D_i-d^{\nabla^i}\in\Omega(A,\End(E_i))^1$. In particular, $\widehat{\RH}_0(1)=0$.
\end{enumerate}
\end{definition}

When we set $\RH(D)=\hat\Sigma^\nabla(\omega)\in C_{\Pi_\infty(A)}(\End(E))^1$, we regard this as a $G$-connection on $E$ through Proposition \ref{prop:connections-simplicial-description}, so that the $G$-connection properly is $D_0+\hat\Sigma^\nabla(\omega)$.

\begin{theorem} \label{thm:riemann-hilbert-lie-algebroids}
The Riemann--Hilbert functor for Lie algebroids is indeed a well-defined $A_\infty$-functor. In particular, curvature commutes with integration
\begin{equation} \label{eq:curvature-commutes-with-integration}
K_{\RH(D)}=\widehat{\RH}_1(K_D)
\end{equation}
and hence $\RH$ descends to an $A_\infty$-functor
\[
\RH:\HRep(A)\to\HRep(\Pi_\infty(A)).
\]
\end{theorem}

\begin{proof}
Let us first see that it is well defined. Let $(E,D)\in\HConn(A)$, let $\nabla$ be any degree-preserving $A$-connection on $E$ and $\omega:=D-d^\nabla$. First of all, $\hat\Sigma(\omega)$ converges, since the only part of $\omega$ that can lead to infinite sums is its part in $\Omega^1(A,\End^0(E))$, and the sums concerning such part converge by a similar argument to Example \ref{ex:parallel-transport-as-iterated-integral}. Secondly, to see that $\RH(D)$ does not depend on $\nabla$, notice that for any $\alpha\in\Omega^1(A,\End^0(E))$ we have that, by Proposition \ref{prop:gauge-equivariance},
\begin{align*}
\hat\Sigma^{\nabla+\alpha}(\omega)&=\sum_{p\geq 0}\hat\Sigma^{\nabla+\alpha}_p(\omega^{\otimes p})\\
&=\sum_{p\geq 0} \sum_{m_i\geq 0} \hat\Sigma^\nabla_{m_0+\dots+m_p+p}(\alpha^{m_0}\otimes\omega\otimes \alpha^{m_1}\otimes\dots\otimes \alpha^{\otimes m_{p-1}}\otimes\omega\otimes\alpha^{m_p}) \\
&=\sum_{p\geq 0}\hat\Sigma^\nabla_p((\omega+\alpha)^{\otimes p}) \\
&=\hat\Sigma^\nabla(\omega+\alpha).
\end{align*}
Since $D-d^{\nabla+\alpha}=\omega-\alpha$, we then see that
\[
\hat\Sigma^{\nabla+\alpha}(\omega-\alpha)=\hat\Sigma^\nabla(\omega),
\]
which proves that $\RH(D)$ is well defined.

Let now $(E_i,D_i)\in\HConn(A)$, for $i=0,\dots,\ell$, and $\eta_i\in\Omega(A,\Hom(E_i,E_{i-1}))$. Let $\nabla^i$ be any degree-preserving $A$-connection on $E_i$ and $\omega_i:=D_i-d^{\nabla^i}$. Let $E:=\bigoplus_{i=0}^\ell E_i$, and consider $\omega_i$ and $\eta_i$ as elements of $\Omega(A,\End(E))$. Let $\omega:=\sum_{i=0}^\ell\omega_i$ and $\eta:=\sum_{i=1}^\ell\eta_i$, and $\nabla:=\bigoplus_{i=0}^\ell \nabla^i$, which is a degree-preserving $A$-connection on $E$. If $\Phi\in\End(E)$, we denote by $\Phi_{ij}$ the component $E_j\to E_i$, so that
\[
\omega_{ij}=\begin{cases}
    \omega_j, & \text{if $i=j$}, \\
    0, & \text{otherwise.}
\end{cases}
\quad\text{and}\quad 
\eta_{ij}=\begin{cases}
    \eta_j, & \text{if $i=j-1$}, \\
    0, & \text{otherwise.}
\end{cases}
\]
With this notation, we have that
\[
\widehat{\RH}_\ell(\eta_1\otimes\dots\otimes\eta_\ell)=\hat\Sigma^\nabla(\omega+\eta)_{0\ell}.
\]
In particular, if $\alpha_i\in\Omega^1(A,\End^0(E^i))$ and $\alpha:=\sum_{i=0}^\ell\alpha_i\in\Omega^1(A,\End^0(E))$, then, by the previous computation,
\[
\hat\Sigma^{\nabla+\alpha}(\omega-\alpha+\eta)=\hat\Sigma^\nabla(\omega+\eta),
\]
so $\widehat{\RH}_\ell(\eta_1\otimes\dots\otimes\eta_\ell)$ is also well defined.

Using the notation from Section \ref{sect:curved-dg-algebras} and Theorem \ref{thm:good-cube-to-simplex-family-first-version}, we can now write
\begin{align*}
d&(\widehat{\RH}_\ell(\eta_1\otimes\dots\otimes\eta_\ell))+\RH(D_0)\widehat{\RH}_\ell(\eta_1\otimes\dots\otimes\eta_\ell) \\
&\qquad -(-1)^{P_1^\ell+Q_1^\ell+1}\widehat{\RH}_\ell(\eta_1\otimes\dots\otimes\eta_\ell)\RH(D_\ell) \\
&\qquad +\sum_{i=1}^{\ell-1}(-1)^{P_1^i+Q_1^i}\widehat{\RH}_i(\eta_1\otimes\dots\otimes\eta_i)\widehat{\RH}_{\ell-i}(\eta_{i+1}\otimes\dots\otimes\eta_\ell) \\
&=d(\hat\Sigma^\nabla(\omega+\eta)_{0\ell})-\epsilon(\hat\Sigma^\nabla(\omega+\eta)_{00})\hat\Sigma^\nabla(\omega+\eta)_{0\ell} \\
&\qquad -\epsilon(\hat\Sigma^\nabla(\omega+\eta)_{0\ell})\hat\Sigma^\nabla(\omega+\eta)_{\ell\ell} \\
&\qquad -\sum_{i=1}^{\ell-1} \epsilon(\hat\Sigma^\nabla(\omega+\eta)_{0i} )\hat\Sigma^\nabla(\omega+\eta)_{i\ell} \\
&=\left( d(\hat\Sigma^\nabla(\omega+\eta))-\epsilon(\hat\Sigma^\nabla(\omega+\eta))\hat\Sigma^\nabla(\omega+\eta) \right)_{0\ell} \\
&=\left( \hat K_{\hat\Sigma^\nabla(\omega+\eta)} \right)_{0\ell}.
\end{align*}
Applying now Proposition \ref{prop:curvature-in-A-infty-algebras} (or, rather, Equation \eqref{eq:curvature-in-A-infty-algebras-hat-notation}), we have that
\[
\hat K_{\hat\Sigma^\nabla(\omega+\eta)}= \sum_{i,j\geq 0} (-1)^i\hat\Sigma^\nabla\left( \epsilon(\omega+\eta)^{\otimes i} \otimes \hat K_{\omega+\eta}\otimes (\omega+\eta)^{\otimes j} \right).
\]
But, if we set $D:=d^\nabla+\omega$, which is a higher $A$-connection on $E$, then
\begin{align*}
\hat K_{\omega+\eta} &= K^\nabla+d^\nabla(\omega+\eta)-\epsilon(\omega+\eta)\wedge(\omega+\eta) \\
&=K^\nabla+d^\nabla\omega+\omega\wedge\omega+d^\nabla\eta+\omega\wedge\eta-\epsilon(\eta)\wedge\omega -\epsilon(\eta)\wedge\eta \\
&=D^2+D\eta-\epsilon(\eta)\wedge\eta.
\end{align*}
Hence,
\begin{align*}
\hat K_{\hat\Sigma^\nabla(\omega+\eta)} &=\sum_{i,j\geq 0} \hat\Sigma^\nabla((\omega-\epsilon(\eta))^{\otimes i} \otimes D\eta \otimes (\omega+\eta)^{\otimes j}) \\
&\qquad +\sum_{i,j\geq 0} \hat\Sigma^\nabla((\omega-\epsilon(\eta))^{\otimes i} \otimes (-\epsilon(\eta)\wedge\eta) \otimes (\omega+\eta)^{\otimes j}) \\
&\qquad +\sum_{i,j\geq 0} \hat\Sigma^\nabla((\omega-\epsilon(\eta))^{\otimes i} \otimes D^2 \otimes (\omega+\eta)^{\otimes j}).
\end{align*}
Looking at the $0\ell$ component gives that $\RH$ is indeed an $A_\infty$-functor.
\end{proof}

It is instructive to write what equation \eqref{eq:curvature-commutes-with-integration} states in detail. If $\sigma\in \Pi_\infty(A)_k$, then
\begin{equation} 
\begin{aligned}\label{eq:curvature-commutes-with-integration-explicit}
-\sum_{j=1}^{k-1} (-1)^j \sum_{m\geq 0}\int_{d_j\sigma} \omega^{\otimes m}+\sum_{j= 0}^k (-1)^j \sum_{m_0,m_1 \geq 0} \int_{a_j\sigma} \omega^{\otimes m_0} \int_{b_{k-j}\sigma} \omega^{\otimes m_1}= \\
= \sum_{m_0,m_1\geq 0} \int_\sigma \omega^{\otimes m_0}\otimes K_D\otimes \omega^{\otimes m_1}.
\end{aligned}
\end{equation}
Here the integrals over simplices should be interpreted as the integrals along pinched cubes once everything is pulled back by $\theta$, e.g.,
\[
\int_{d_j\sigma}\omega^{\otimes m}:=\int_{d_j\sigma\circ \theta_{k-1}} (\theta_{k-1}^*\sigma^*\omega)^{\otimes m}.
\]

For a usual $TM$-connection $\nabla$ on a vector bundle, equation \eqref{eq:curvature-commutes-with-integration} means that for a contractible path $\gamma$ with null-homotopy $H$, we have 
\begin{equation*}
    \tau^\nabla_\gamma = \1+\int_H R^\nabla.
\end{equation*}

\begin{remark}
One can define the Riemann--Hilbert functor slightly differently. If $(E,D)\in\HConn(A)$ and $\sigma\in\Pi_\infty(A)_k$, one can fix a trivialization $\sigma_0^*E\cong \Delta^k\times V$, for a graded vector space $V$, so that under such trivialization one can write $\sigma^*D=d+\omega$, for $d$ the de Rham differential on $\Omega(\Delta^k,V)=\Omega(\Delta^k)\otimes V$ and $\omega\in\Omega(\Delta^k, \End(V))^1$, and set
\[
\RH(D)(\sigma):=\hat\Sigma(\omega)(\id)=\sum_{m\geq 0} \hat\Sigma_m(\omega^{\otimes m})(\id),
\]
where now $\Sigma:\Omega(\Delta^k)\otimes \End(V)\to C(\Pi_\infty(\Delta^k))\otimes \End(V)$ is the $\End(V)$-valued de Rham $A_\infty$-map \cite{gugenheim}, which hinges on the classical (abelian) Stokes's theorem. One then has to prove that this does not depend on the trivialization \cite{arias-abad-schatz}, i.e., that this definition is gauge invariant. The proof that the Riemann--Hilbert functor thus defined is an $A_\infty$-functor is implicit in \cite[Prop. 3.3]{block-smith} or \cite[Prop 1.5.2]{chen}, although they do not consider the curved case. Similarly, equation \eqref{eq:curvature-commutes-with-integration-explicit} is equation (3.4) in \cite{block-smith}. The difference between our definition and this definition is that we work globally by choosing a global $A$-connection $\nabla$ and use its parallel transport to trivialize $E$ along $\sigma$ (or, to be more precise, along $\sigma\circ\theta_k$). The independence on $\nabla$ (cf. Proposition \ref{prop:gauge-equivariance}) is the analog of the local gauge invariance of \cite{arias-abad-schatz}.
\end{remark}

\subsection{The Riemann--Hilbert functor for higher Lie algebroids}

In this section $(\Ascr\to M,Q)$ is a higher Lie algebroid, and we fix a cube-to-simplex family $\theta$. 

The definition of $\Pi_\infty(\Ascr)$ is the obvious one.

\begin{definition}
The \textbf{homotopy groupoid $\Pi_\infty(\Ascr)$ of $\Ascr$} is defined by setting $\Pi_\infty(\Ascr)_k$ to be the space of higher Lie algebroid maps $\sigma:T\Delta^k\to \Ascr$ and by taking the face and degeneracy maps to be
\begin{align*}
    d_i(\sigma^*,\sigma_0)&:=( \delta_i^*\sigma^*,\sigma_0\circ\delta_i), \\
    s_i(\sigma^*,\sigma_0)&:=(\sigma_i^*\sigma^*,\sigma_0\circ\sigma_i).
\end{align*}
\end{definition}
The below diffeology corresponds to the diffeology used in \cite{zeng2025higher} in case $\Ascr$ is the Lie algebroid associated to a regular foliation.
\begin{lemma} \label{lem:pi-infty-higher-lie-algebroids}
$\Pi_\infty(\Ascr)$ is a diffeological higher groupoid, where the set of plots for $\Pi_\infty(\Ascr)_k$ is the set of maps $p:U\to\Pi_\infty(\Ascr)_k$ for some open $U\subseteq\R^n$ such that the induced map 
\[
\begin{array}{cccc}
\hat p_0: & U\times\Delta^k & \longrightarrow & M \\
& (x,t) & \longmapsto & p(x)_0(t)
\end{array}
\]
is smooth and for every $\xi\in C^\ell(\Ascr)$ the map $p^*_\xi:U\to\Omega^\ell(\Delta^k)$ given by $x\mapsto p(x)^*\xi$ is smooth. In the above, $p_0$ denotes the map $U \to \Pi_\infty(M)_k$ induced by $p$.
\end{lemma}

\begin{remark}
If we denote by $\pi_1:U\times\Delta^k\to U$ and $\pi_2:U\times\Delta^k\to\Delta^k$ the projections, then $\ker \pi_{1*}\to U\times\Delta^k$ is a Lie algebroid and, for $\xi\in C^\ell(\Ascr)$, the expression
\[
(\hat p^*\xi)_{(x,t)}:=(\pi_2^* p(x)^*\xi)_{(x,t)}
\]
defines a set-theoretic differential form on $\ker\pi_{1*}$. The conditions on $p$ from Lemma \ref{lem:pi-infty-higher-lie-algebroids} are equivalent to $(\hat p^*,\hat p_0)$ being a higher Lie algebroid morphism $\ker \pi_{1*}\to\Ascr$.
\end{remark}

\begin{proof}
The fact that $\Pi_\infty(\Ascr)$ is a higher groupoid follows easily from the fact that the standard horn spaces 
\[
\Lambda_{k,\ell}=\bigcup_{\substack{i=0 \\ i\neq \ell}}^k \delta_i(\Delta^{k-1})\subseteq \Delta^k
\]
are strong deformation retracts of $\Delta^k$, as can be seen from taking the barycenter of $\delta_\ell(\Delta^{k-1})\subseteq\Delta^k$ and ``pulling'' it to the vertex $v_\ell$.

Let us now check that the set of plots on $\Pi_\infty(\Ascr)_k$ indeed defines a diffeology. Firstly, if $p:\R^n\to \Pi_\infty(\Ascr)_k$ is constant, say $p(x)=\sigma$ for some $\sigma\in\Pi_\infty(\Ascr)_k$, then $\hat p_0(x,t)=\sigma_0(t)$ and $p_\xi^*(x)=\sigma^*\xi$, both of which are smooth, so $p$ is a plot. Secondly, if $p:U\to\Pi_\infty(\Ascr)_k$ is a map such that for every $x\in U$ there is a neighborhood $V$ containing $x$ such that $p|_V$ is a plot, then by the glueing property for smooth maps $U\times\Delta^k\to M$ and for smooth functions on $U\times\Delta^k$ we have that $p$ is a plot. Lastly, if $p:U\to\Pi_\infty(\Ascr)_k$ is a plot and $F:V\to U$ is smooth, for some open $V\subseteq \R^m$, then
\[
\widehat{(p\circ F)}_0=\hat p_0\circ (F\times\id)
\]
and
\[
(p\circ F)^*_\xi=p_\xi^*\circ F,
\]
both of which are smooth, so $p\circ F$ is a plot.

For any plot $p:U\to\Pi_\infty(\Ascr)_k$ we have that
\[
\widehat{(d_i\circ p)}_0=\hat p_0\circ(\id\times\delta_i)
\]
and
\[
(d_i\circ p)_\xi^*=\delta_i^*\circ p_\xi^*,
\]
both of which are smooth, so that $d_i\circ p$ is a plot. Therefore, $d_i$ is smooth. An analogous argument shows that $s_i$ is also smooth.
\end{proof}

In our construction of the Riemann--Hilbert functor for higher algebroids, we exploit that $T\Delta^k$ is a Lie algebroid. Given a simplex $\sigma : T\Delta^k \to \Ascr$, we first pullback the $\Ascr$-connection on $E$ to $T\Delta^k$ and then apply Riemann--Hilbert functor for algebroids to $T\Delta^k$.

\begin{definition}
We define the \textbf{Riemann--Hilbert functor} as the $A_\infty$-functor
\[
\RH:\HConn(\Ascr)\to \HConn(\Pi_\infty(\Ascr))
\]
constructed as follows:
\begin{enumerate}
    \item If $(E,D)\in \HConn(\Ascr)$, we define $\RH(D)\in C_{\Pi_\infty(\Ascr)}(\End(E))^1$ as follows: for $\sigma\in\Pi_\infty(\Ascr)_k$ and any degree-preserving $T\Delta^k$-connection $\nabla$ on $\sigma_0^*E\to \Delta^k$, let $\omega:=\sigma^*D-d^\nabla$ and set
    \[
    \RH(D)(\sigma):=\hat\Sigma^\nabla(\omega)(\id),
    \]
    where we view $\id:\Delta^k\to\Delta^k$ as a $k$-simplex in $\Delta^k$.

    \item If $(E_i,D_i)\in\HConn(\Ascr)$, for $i=0,\dots,\ell$, and $\eta_i\in C(\Ascr,\Hom(E_i,E_{i-1}))$, then we define $\widehat{\RH}_\ell(\eta_\otimes\dots\otimes\eta_\ell)\in C_{\Pi_\infty(\Ascr)}(\End(E))$ as follows: for $\sigma\in\Pi_\infty(\Ascr)_k$ and any degree-preserving $T\Delta^k$-connections $\nabla^i$ on $\sigma_0^*E_i\to\Delta^k$, let $\omega_i:=\sigma^*D-d^{\nabla^i}$ and set
    \begin{align*}
    \widehat{\RH}_\ell(\eta_1 \otimes & \dots  \otimes\eta_\ell)(\sigma):= \\
    & := \hat\Sigma^{\nabla^0,\dots,\nabla^\ell}(\omega_0\otimes\sigma^*\eta_1\otimes\omega_1\otimes\dots\otimes\omega_{\ell-1}\otimes\sigma^*\eta_\ell\otimes\omega_\ell)(\id).
    \end{align*}
    In particular, $\widehat{\RH}_0(1)=0$.
\end{enumerate}
In the above, $\hat\Sigma^\nabla$ denotes the $A_\infty$-map from Definition \ref{def:Sigma} for $A = T\Delta^k$.
\end{definition}

\begin{theorem}\label{thm:riemann-hilbert-higher-algebroids}
The Riemann--Hilbert functor for higher Lie algebroids is indeed a well-defined $A_\infty$-functor. In particular, curvature commutes with integration
\begin{equation*} 
K_{\RH(D)}=\widehat{\RH}_1(K_D)
\end{equation*}
and hence $\RH$ descends to an $A_\infty$-functor
\[
\RH:\HRep(\Ascr)\to\HRep(\Pi_\infty(\Ascr)).
\]
\end{theorem}

\begin{proof}
This follows directly from Theorem \ref{thm:riemann-hilbert-lie-algebroids}.
\end{proof}

\bigskip
\section{The higher general linear groupoid} \label{sect:higher-general-linear-groupoid}

Throughout this section, $E \to M$ denotes a graded vector bundle of length $n$. We can associate to it a curved dg-category $C_E$, whose objects are $(x,\delta_x)$, for $x\in M$ and $\delta_x\in\End^1(E_x)$, and whose space of morphisms are
\[
\Hom((x_1,\delta_{x_1}),(x_0,\delta_{x_0})):=\Hom(E_{x_1},E_{x_0}),
\]
with the connection
\[
D\Phi:=\delta_{x_0}\Phi-(-1)^{\abs{\Phi}}\Phi\delta_{x_1},
\]
and with curvature morphisms $K_{(x,\delta_x)}:=\delta_x^2\in\End^2(E_x)$. Note that the objects assemble into a bundle $\End^1(E)\to M$ and similarly the morphisms assemble into a bundle $\Hom(\pi^*_1E, \pi_0^*E)\to M \times M$, where $\pi_i$ are the projections. This inspires the central definitions of this section (Definitions \ref{def:globular-collection} and \ref{def:simplicialization}).

\subsection{Simplicialization}

Simplicialization is a process through which we construct a simplicial set from some raw data, that we call a globular collection. The simplicialization has many redundancies, so as to organize the globular collection into a simplicial set.

\begin{definition} \label{def:globular-collection}
    Let $M$ be a manifold. A \textbf{globular collection} on $M$ consists of
    \begin{enumerate}
        \item a collection of manifolds $\{X_k\}_{k \geq 0}$;
        \item a surjective submersion $q_0 : X_0 \to M$;
        \item a collection of surjective submersions $q_k : X_k \to M \times M$ for $k \geq 1$;
        \item a map $u_1 : M \to X_1$ such that 
        \begin{center}
            \begin{tikzcd}
                & X_1 \ar[d, "q_1"]\\
                M \ar[ru, "u_1"] \ar[r, "\Delta"] & M \times M
            \end{tikzcd}
        \end{center}
        commutes;
        \item a collection of sections $u_k$ of $q_k$ for $k \geq 2$.
    \end{enumerate}
    We say that a globular collection $\{X_k\}$ has length $n$ if $q_n$ is not a diffeomorphism and $q_k$ is a diffeomorphism for $k > n$. We call the maps $u_k$ the unit maps of $X$. 
\end{definition}

The examples of globular collections that we care about, are all associated to a graded vector bundle $E \to M$ of finite length $n$.

\begin{definition} \label{def:globular-collections-on-E}
The \textbf{degree 1 globular collection} for $E$ is given by the bundles
    \begin{align*}
        X_0 & = \End^1(E), \\
        X_k & = \Hom^{1 - k}(\pi_1^*E, \pi_0^*E),
    \end{align*}
where $k \geq 1$ and $\pi_j$ denotes the $(j +1)$th projection $M^2 \to M$. The units maps are given by 
\begin{align*}
    u_1(x) & = \id_{E_x}\in\Hom^0(E_x,E_x), \\
    u_k(x,y) & = 0 \in \Hom^{1 - k}(E_y, E_x).
\end{align*}

The \textbf{groupoidal globular collection} for $E$ has bundle data
\begin{align*}
    X_0 & = \End^1(E), \\
    X_1 & = \GL^0(\pi_1^*E, \pi_0^*E), \\
    X_k & = \Hom^{1 - k}(\pi_1^*E, \pi_0^*E),
\end{align*}
with $k \geq 2$, and where
\[
\GL^0(\pi_1^*E,\pi_0^*E):=\{ A\in\Hom^0(\pi_1^*E,\pi_0^*E): \text{$A$ is an isomorphism}\}.
\]
The unit maps are given by 
\begin{align*}
    u_1(x) & = \id_{E_x}\in\GL^0(E_x,E_x), \\
    u_k(x,y) & = 0 \in \Hom^{1 - k}(E_y, E_x).
\end{align*}

The \textbf{degree 2 globular collection} for $E$ is given by the bundles
    \begin{align*}
        X_0 & = \End^2(E), \\
        X_k & = \Hom^{2 - k}(\pi_1^*E, \pi_0^*E),
    \end{align*}
    with $k \geq 1$. The unit maps are given by
    \begin{align*}
        u_1(x)& =0\in \Hom^1(E_x,E_x), \\
        u_k(x,y)& =0\in\Hom^{2-k}(E_y,E_x).
    \end{align*}
\end{definition}

\begin{notation} Let $q_i:X_i\to M$ be surjective submersions with $i=0,\dots, k$ and let $\phi : N \to M$ be a smooth map, then we denote by:
\begin{enumerate}
    \item $\phi^*X_i$ the fibered product $N\times_MX_i$.

    \item $\bigoplus_{i=0}^k X_i$ the fibered product $X_0\times_MX_1\times_M\dots\times_MX_k$.
\end{enumerate}
\end{notation}

Let us now fix a globular collection $X = \{X_k\}_{k \in \mathbb{N}}$ of length $n$.

\begin{definition}\label{def:simplicialization}
    The \textbf{simplicialization} of a globular collection $X$ is the simplicial set $S(X)$ defined as follows:
    \begin{enumerate}
        \item $S(X)_k$ is the total space of the surjective submersion
    \begin{equation*}
        \left(\bigoplus_{\alpha : [0] \hookrightarrow [k]} \pi_{\alpha(0)}^*X_0 \right) \oplus \bigoplus_{\substack{\alpha : [i] \hookrightarrow [k]\\ i > 0}}(\pi_{\alpha(0)}, \pi_{\alpha(i)})^*X_i  \longrightarrow M^{k+1},
    \end{equation*}
    where $\pi_j$ denotes the $(j +1)$th projection $M^{k +1} \to M$. The surjective submersion $S(X)_k\to M^{k+1}$ is denoted by $p_k$. 
    For $\alpha : [i] \hookrightarrow [k]$ and $A \in S(X)_k$, we write $A_\alpha$ for the $\alpha$ component of $A$.
    \item For $\gamma : [k] \to [j]$ order-preserving, the induced map $\gamma^*:S(X)_j\to S(X)_k$ makes
    \[\begin{tikzcd} \label{cd:simplicialization_proj}
	{S(X)_j} & {S(X)_k} \\
	{M^{j+1}} & {M^{k+1}}
	\arrow["{\gamma^*}", from=1-1, to=1-2]
	\arrow["{p_j}"', from=1-1, to=2-1]
	\arrow["{p_k}", from=1-2, to=2-2]
	\arrow["{\gamma^*}"', from=2-1, to=2-2]
    \end{tikzcd}\]
    commute, where $\gamma^*:M^{j+1}\to M^{k+1}$ is
    \[
    \gamma^*(x_0,\dots, x_j):=(x_{\gamma(0)},\dots, x_{\gamma(k)}),
    \]
    and it is defined as follows: if $A\in S(X)_j$ and $\alpha:[i]\hookrightarrow [k]$, then
    \begin{equation}\label{eq:structure_maps_simplicialization}
    (\gamma^*A)_\alpha=\begin{cases}
        A_{\gamma\circ\alpha}, & \text{if $\gamma\circ\alpha$ is injective,} \\
        u_1(\pi_{\gamma\circ\alpha(0)}\circ p_j(A)), & \text{if $i=1$ and $\gamma \circ \alpha$ is constant,} \\
        u_i(\pi_{\gamma\circ\alpha(0)}\circ p_j(A),\pi_{\gamma\circ\alpha(i)}\circ p_j(A)), & \text{otherwise.}
    \end{cases}
    \end{equation}
    \end{enumerate}
\end{definition}

\begin{figure}
    \centering
    \includesvg[width=0.3\textwidth]{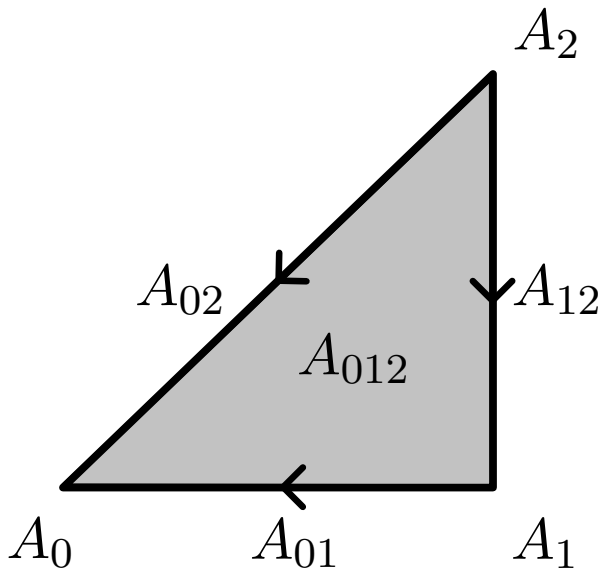}
    \caption{Illustration of a 2-simplex in $S(X)_2$.}
    \label{fig:2-simplex-simplialization}
\end{figure}

\begin{theorem} \label{thm:simplicialization}
The simplicialization $S(X)$ of a globular collection $X$ is a contractible higher Lie groupoid. Moreover, if $X$ has length $n$, then $S(X)$ is a Lie $(n + 1)$-groupoid.
\end{theorem}

\begin{proof}
    We defer this proof to Appendix \ref{app:simplicialization}. Lemma \ref{lemma:simplicialization_simp_manifold} shows that $S(X)$ is a simplicial manifold, Proposition \ref{prop:simplicialization_groupoid} shows that $S(X)$ is a higher Lie groupoid and proposition \ref{prop:simplicialization_contractible} shows that $S(X)$ is contractible.
\end{proof}

 Thus, we see that simplicialization organizes the globular data into a contractible groupoid. Figure \ref{fig:2-simplex-simplialization} shows an illustration of how it does so for the 2-simplices. Even if $S(X)$ is contractible, there can be interesting relations within simplicializations of globular data that give non-contractible groupoids.

\subsection{Higher endomorphism and general linear groupoids}

\begin{definition} \label{def:higher_end_groupoid}
    The \textbf{higher endomorphism groupoid} $\End_\infty(E)$ of $E$ is the simplicialization of its degree 1 globular data. Explicitly, we have that
    \begin{equation*}
        \End_\infty(E)_k = \bigoplus_{\alpha : [i] \hookrightarrow[k]} \Hom^{1 - i}(\pi^*_{\alpha(i)}E, \pi^*_{\alpha(0)}E)\to M^{k+1}.
    \end{equation*}
    The \textbf{curvature bundle} $\Curv_\infty(E)$ of $E$ is the simplicialization of its degree 2 globular data. Explicitly, we have that
    \begin{equation*}
        \Curvi(E)_k = \bigoplus_{\alpha : [i] \hookrightarrow[k]} \Hom^{2 - i}(\pi^*_{\alpha(i)}E, \pi^*_{\alpha(0)}E)\to M^{k+1}.
    \end{equation*}
\end{definition}

\begin{remark}
    By proposition \ref{prop:simplicialization_vector_bundle}, $\Curvi(E)\to M^{\bullet +1}$ is a simplicial vector bundle. On the contrary, $\Endi(E)$ is not, since
    \begin{equation*}
        s_0 : \End^1(E) \to \End^1(\pi^*_0E) \oplus \End^0(\pi^*_1E, \pi^*_0E) \oplus \End^1(\pi^*_1E)
    \end{equation*}
    is given by
    \begin{equation*}
        \delta_x \mapsto (\delta_x, \id_x, \delta_x),
    \end{equation*}
    which is not linear.
\end{remark}

$\Endi(E)$ carries a canonical connection on $E$, which is universal in the sense that any higher groupoid connection on $E$ is a pullback of this universal connection, as we will see.

\begin{definition} \label{def:universal-connection-and-curvature-section}
    The \textbf{universal connection} of $\Endi(E)$ is the $\Endi(E)$-connection $D_E$ on $E$ defined by the degree 1 cochain $F \in C_{\Endi(E)}(\End(E))^1$ with components
    \begin{equation*}
        F_i(A) := A_{\id:[i]\hookrightarrow[i]} \in\Hom^{1-i}(E_{v_k(A)}, E_{v_0(A)}),
    \end{equation*}
    where $A \in \End_\infty(E)_k$. We call the curvature $K_{D_E}$ of $D_E$ the \textbf{universal curvature} of $\Endi(E)$.

    In addition, the \textbf{curvature map} of $\Endi(E)$ is the map $K_E : \Endi(E) \to \Curvi(E)$ defined by
    \begin{equation} \label{eq:curvature_end}
        K_E(A)_{\alpha:[i]\hookrightarrow [k]} := -\sum_{j=1}^{i-1}(-1)^jA_{\alpha \circ \delta^j}+\sum_{j=0}^i(-1)^j A_{\alpha \circ \alpha^j} A_{ \alpha\circ \beta^{i-j}}. 
I    \end{equation}
Furthermore, let $G$ be a higher Lie groupoid with $G_0 = M$ and let $F : G \to \Endi(E)$ be a morphism of simplicial manifolds such that $F_0 : M \to \Endi(E)_0$ is a section of $\End^1(E)$. Then, we call $F$ \textbf{flat} if
\begin{equation*}
    K_E \circ F = 0 \circ p,
\end{equation*}
where $p : G \to M^{\bullet + 1}$ is the map given on $G_k$ by 
\begin{equation}\label{eq:canonical_proj_p}
    p_k(x) = (v_0(x), \dots, v_k(x)).
\end{equation}
\end{definition}

\begin{lemma} \label{lem:universal-curvature-is-the-curvature-of-universal-connection}
    The curvature map $K_E : \End_\infty(E) \to \Curv_\infty(E)$ is a map of simplicial manifolds such that $p\circ K_E=p$. In particular, it induces a simplicial section $K_E\in\Gamma(\Endi(E),p^*\Curvi(E))$, which represents the universal curvature of $E$, i.e., for all $A\in\Endi(E)_k$ we have that
    \[
    K_E(A)_{\alpha:[i]\hookrightarrow [k]}=K_{D_E}(\alpha^*A).
    \]
\end{lemma}
\begin{proof}
    First of all, notice that if $A\in\Endi(E)_k$ and $\gamma:[k]\to[\ell]$ and $\alpha:[i]\hookrightarrow [k]$ are such that $\gamma\circ\alpha$ is injective, then
    \begin{align*}
        (\gamma^*\circ K_E)(A)_\alpha & = -\sum_{j=1}^{i-1}(-1)^jA_{\gamma\circ \alpha \circ \delta^j}+\sum_{j=0}^i(-1)^j A_{\gamma \circ \alpha \circ \alpha^j} A_{ \gamma\circ\alpha\circ \beta^{i-j}} \\
        & =  -\sum_{j=1}^{i-1}(-1)^j(\gamma^*A)_{ \alpha \circ \delta^j}+\sum_{j=0}^i(-1)^j (\gamma^*A)_{\alpha \circ \alpha^j} (\gamma^*A)_{ \alpha\circ \beta^{i-j}} \\
        & = (K_E\circ\gamma^*)(A)_\alpha.
    \end{align*}

    To prove that $K_E$ is a map of simplicial manifolds we need to show that $K_E$ commutes with the face and degeneracy maps of $\End_\infty(E)$ and $\Curvi(E)$. Firstly, remark that, as $\delta_i : [k] \to [k + 1]$ is injective, $\delta_i \circ \alpha$ is injective for any injective $\alpha : [i] \hookrightarrow [k]$. Thus, by the computation above, we see that $\delta_i^*$ commutes with $K_E$.

    Now consider an injective $\alpha\colon [i] \hookrightarrow [k]$ and $\sigma^\ell \colon [k] \to [k -1]$. Then, if $\sigma^\ell \circ \alpha$ is injective, we are done by the previous computation. If it is not injective, then we see that $((\sigma^\ell)^*\circ K_E)(A)_\alpha = 0$ and 
    \begin{equation*}
        (K_E\circ(\sigma^\ell)^*)(A)_\alpha = -\sum_{j=1}^{i-1}(-1)^j((\sigma^\ell)^*A)_{ \alpha \circ \delta^j}+\sum_{j=0}^i(-1)^j ((\sigma^\ell)^*A)_{\alpha \circ \alpha^j} ((\sigma^\ell)^*A)_{ \alpha\circ \beta^{i-j}}.
    \end{equation*}
    Note that as $\sigma^\ell \circ \alpha$ is not injective, we have $i > 0$. Suppose first that $i = 1$. Since $\sigma^\ell \circ \alpha$ is not injective, we see that $\alpha(0) = \ell$, $\alpha(1) = \ell + 1$ and $\ell = \sigma^\ell(\alpha(0)) = \sigma^\ell(\alpha(1))$. Then, using that $((\sigma^\ell)^*A)_\alpha = \id$, we have that
    \begin{align*}
        (K_E\circ(\sigma^\ell)^*)(A)_\alpha & = ((\sigma^\ell)^*A)_{\alpha \circ \alpha^0} ((\sigma^\ell)^*A)_{ \alpha\circ \beta^{1}} - ((\sigma^\ell)^*A)_{\alpha \circ \alpha^1} ((\sigma^\ell)^*A)_{ \alpha\circ \beta^{0}} \\
        & = ((\sigma^\ell)^*A)_{\alpha \circ \delta^1} ((\sigma^\ell)^*A)_{ \alpha} - ((\sigma^\ell)^*A)_{\alpha} ((\sigma^\ell)^*A)_{ \alpha\circ \delta^{0}} \\
        & = A_{\upsilon^\ell} - A_{\upsilon^\ell} = 0.
    \end{align*}
    
    On the other hand, suppose now that $i > 1$. Note that as $\alpha$ is injective and order preserving, it is fully characterized by its image. Then, since $\sigma^\ell \circ \alpha$ is not injective, we see that $\ell, \ell +1 \in \im(\alpha)$ and we have a unique element $\ell_0$ such that $\alpha(\ell_0) = \ell$ and $\alpha(\ell_0 + 1) = \ell + 1$. 
    
    Now, suppose that $\ell_0 \neq 0$ or $\ell_0 \neq i-1$. We first consider the terms involving
    \begin{equation*}
        ((\sigma^\ell)^*A)_{\alpha \circ \delta^j},
    \end{equation*}
    with $0 < j < i$.  Note that $\sigma^\ell \circ \alpha \circ \delta^j$ is injective if and only if $j = \ell_0$ or $j = \ell_0 + 1$.
    \begin{equation*}
        -\sum_{j=1}^{i-1}(-1)^j((\sigma^\ell)^*A)_{ \alpha \circ \delta^j} = -(-1)^{\ell_0}A_{\sigma^\ell \circ \alpha \circ \delta^{\ell_0}} - (-1)^{\ell_0 + 1}A_{\sigma^\ell\alpha \circ \delta^{\ell_0 + 1}} =0,
    \end{equation*}
    since $\sigma^\ell \circ \alpha \circ\delta^{\ell_0} = \sigma^\ell \circ \alpha \circ\delta^{\ell_0 + 1}$.
    Secondly, let us consider the terms involving
    \begin{equation}\label{eq:prod_curvature_term}
        ((\sigma^\ell)^*A)_{\alpha \circ \alpha^j} ((\sigma^\ell)^*A)_{ \alpha\circ \beta^{i-j}}.
    \end{equation}
    Note that 
    \begin{align*}
        \im(\alpha \circ \alpha^j) & = \{\alpha(0), \dots, \alpha(j)\}, \\
        \im(\alpha \circ \beta^{i - j}) & = \{\alpha(j), \dots, \alpha(i)\}.
    \end{align*}
    Thus, we see that $\{\ell, \ell + 1\}$ is contained in either $\im(\alpha \circ \alpha^j)$ or $\im(\alpha \circ \beta^{i - j})$. Thus, either $\sigma^\ell \circ \alpha \circ \alpha^i$ or $\sigma^\ell \circ \alpha \circ \beta^{i -j }$ is not injective. As $\ell_0 \neq 0$ and $\ell_0 \neq i-1$, we see that the domain of the function that is not injective has at least 3 elements, which implies that \eqref{eq:prod_curvature_term} vanishes. Therefore, we conclude that $(K_E \circ (\sigma^\ell)^*)(A)_\alpha$ vanishes.
    
    If $\ell_0=0$, then we have that 
    \begin{equation*}
         -\sum_{j=1}^{i-1}(-1)^j((\sigma^\ell)^*A)_{ \alpha \circ \delta^j} = A_{\sigma^\ell \circ\alpha \circ \delta^1},
    \end{equation*}
    and since $((\sigma^\ell)^*A)_{\alpha \circ \alpha^1} = \id$, we have that
    \begin{equation*}
        \sum_{j=0}^i(-1)^j ((\sigma^\ell)^*A)_{\alpha \circ \alpha^j} ((\sigma^\ell)^*A)_{ \alpha\circ \beta^{i-j}} = - A_{\sigma^\ell \circ \alpha \circ \beta^{i -1}}.
    \end{equation*}
    Then, since $\beta^{i - 1} = \delta^0$ and $\sigma^\ell \circ \alpha \circ \delta^0 = \sigma^\ell \circ \alpha \circ \delta^1$, we conclude that $(K_E \circ(\sigma^\ell)^*)(A)_\alpha = 0$. An analogous argument can be made in the case $\ell_0 = i - 1$.
\end{proof}
\begin{remark} \label{rm:quasi-smooth-groupoids}
To be complete, let us also mention that one may form the simplicial subset
\begin{equation*}
    \Neri(E):=K_E^{-1}(0)\subseteq \Endi(E).
\end{equation*}
Set theoretically, this is equivalent to the dg-nerve \cite[\href{https://kerodon.net/tag/00PK}{Subsect. 00PK}]{kerodon} of $\Flat(C_E)$, where $C_E$ is the curved dg-category described at the beginning of this section. This is an $n$-category. Inside any $n$-category, one has an $n$-groupoid called the core \cite[\href{https://kerodon.net/tag/01D0}{Subsect. 01D0}]{kerodon}. The core of $\Neri(E)$ is the simplicial subset $\qGaui(E)\subseteq \Neri(E)$ given by
\begin{equation*}
    \qGaui(E)_k:=\{A\in\Neri(E)_k\sep \text{$A_\alpha$ is a quasi-isomorphism for all $\alpha:[1]\hookrightarrow[k]$}\}.
\end{equation*}
So far, these have been set theoretic descriptions. However, there is geometry present in $\Neri(E)$. For example, in the case of a bounded complex of vector spaces $V$ with fixed differential $\delta$ is studied in \cite[Sect. 8]{behrend-getzler}. There, it is shown that $\Neri(V)|_\delta$ is a regular $k$-category inside the descent category of finite type schemes. Furthermore, it is shown that the core $\qGaui(V)|_\delta$ becomes a $k$-groupoid. 

Now, set theoretically $\Neri(E)_k$ is realized for each $k \in \mathbb{N}$ as the pullback:
\begin{center}
     \begin{tikzcd}
         \Neri(E)_k \ar[r] \ar[d] & M^{k + 1} \ar[d, "0"] \\
         \Endi(E)_k \ar[r, "K_E"] & \Curvi(E)_k.
     \end{tikzcd}
\end{center}
Then using the theory of \cite[Sect. 1]{behrend-liao-xu}, $\Neri(E)_k$ becomes a quasi-smooth manifold $(\Endi(E)_k, p^*\Curvi(E)_k, p^*K_E)$, i.e., $\Neri(E)$ is cut out inside $\Endi(E)_k$ by the section $p^*\Curvi(E)_k$ of $p^*K_E$. These quasi-smooth manifolds together with the face and degeneracy maps of $\Endi(E)$ and $p^*\Curvi(E)$ define an simplicial object inside quasi-smooth manifolds.

So, there is a more geometric description of $\Neri(E)$. A natural question is if there exists such a more geometric description of $\qGaui(E)$? Does this define a quasi-smooth simplicial manifold? Furthermore, there are several choices of pretopologies on derived manifolds that can be used to define higher groupoids inside derived manifolds \cite[Sect. 4 and App. B]{cueca2026shifted}. With respect to which choices of pretopologies does $\qGaui(E)$ become a higher groupoid? 
\end{remark}

\begin{theorem} \label{thm:groupoid-connections-as-functors}
For any $G$-connection $D$ on $E$ there is a unique smooth functor $F:G\to \Endi(E)$ such that $p\circ F=p$ and $F^*D_E = D$. Then, $K_D$ is given by $F^*K_{D_E}$. Furthermore, if $K_D = 0$, then $F$ becomes a functor $G \to \qGaui(E)$.
\end{theorem}

\begin{proof}
Let $\tilde F\in C_G(\End(E))^1$ be the cochain representing the $G$-connection (cf. Proposition \ref{prop:connections-simplicial-description}). Define $F:G\to\Endi(E)$ as follows: for $x\in G_k$ set $F(x)_{\alpha:[i]\hookrightarrow[k]}:=\tilde F_i(\alpha^*x)$. Clearly, $F$ is smooth and $p\circ F=p$. Finally, to  
see that $F$ is simplicial, we firstly remark that for $x \in G_k$ and $\alpha : [i] \hookrightarrow [k-1]$, the face maps give
\begin{align*}
    (d_jF(x))_{\alpha} & = F(x)_{\delta_j \circ \alpha} \\
    & = \tilde F_i(\alpha^*\delta_j^*x) = F(d_jx)_\alpha.
\end{align*}
For the degeneracy maps, we take $\alpha : [i] \hookrightarrow [k + 1]$. Then, if $\sigma_j \circ \alpha$ is injective, the same computation as before goes through. On the other hand, if $\sigma_j \circ \alpha$ is not injective, then $\sigma_j \circ \alpha = \beta \circ \sigma_r$ with $\beta$ injective. Now, assuming without loss of generality that $i > 1$, we see using the normalization condition on $G$-connections that
\begin{align*}
    F(s_jx)_\alpha & = \tilde F_i(\alpha^*\sigma_j^*x) \\
    & =  \tilde F_i(\sigma_r^*\beta^*x) \\
    & = 0 = (\sigma_jF(x))_\alpha. 
\end{align*}
Conversely, if $F:G\to\Endi(E)$ is a smooth functor, then $\tilde F$ is defined through the pullback $F^*D_E$. Concretely, we obtain $\tilde F_i(x) = F(x)_{\id:[i]\to[i]}$. The constructions above are inverse to each other.

To see that $K_D = F^*K_{D_E}$, we remark that for $x \in G_i$, we have
\begin{align*}
(K_D)_i(x) & =-\sum_{j=1}^{i-1}(-1)^j\tilde F_{i-1}(d_jx) +\sum_{j= 0}^i (-1)^j \tilde F_j(a_jx) \tilde F_{i-j}(b_{i-j}x) \\
&=-\sum_{j=1}^{i-1}(-1)^jF(x)_{\delta^j} +\sum_{j= 0}^i (-1)^jF(x)_{\alpha^j}F(x)_{\beta^{i-j}} \\
&=K_E(F(x))_{\id} = F^*K_{D_E}(x).
\end{align*}

Lastly, assume that $K_D = F^*K_{D_E}=0$. Note as $K_E$ is defined in terms of $K_{D_E}$, we have $F:G\to\Neri(E)$. To see that actually $F:G\to\qGaui(E)$, we need to check that for all $x\in G_1$ we have that $g:=F(x)_{\id:[1]\to[1]}=\tilde F_1(x)$ is a quasi-isomorphism between $\delta_0:=F(x)_{\upsilon^0:[0]\to[1]}=\tilde F_0(v_0x)$ and $\delta_1:=F(x)_{\upsilon^1:[0]\to[1]}=\tilde F_0(v_1x)$. This is a standard argument, but let us make it explicit for completeness. First of all, we have that
\[
0=\tilde K_1(x)=\tilde F_0(v_0x)\tilde F_1(x)-\tilde F_1(x)\tilde F_0(v_1x)=\delta_0g-g\delta_1,
\]
so $g$ is a cochain map. Consider now $(x,s_0v_1x)\in G_{2,2}$ and let $y\in G_2$ be such that $d_{2,2}(y)=(x,s_0v_1x)$. Then
\begin{align*}
0 &=\tilde K_2(y) = \tilde F_1(d_1y)+\tilde F_0(v_0y)\tilde F_2(y)-\tilde F_1(a_1y)\tilde F_1(b_1y)+\tilde F_2(y)\tilde F_0(v_2y)\\
&= \tilde F_1(s_0v_1x)+\tilde F_0(v_1x)\tilde F_2(y)-\tilde F_1(d_2y)\tilde F_1(x)+\tilde F_2(y)\tilde F_0(v_1x)\\
&= \id +\delta_1 \tilde F_2(y)-\tilde F_1(d_2y)g+\tilde F_2(y)\delta_1,
\end{align*}
so $g$ admits a left inverse up to homotopy. Similarly, using $(s_0v_0x,x)\in G_{2,0}$, we see that the $g$ also admits a right inverse up to homotopy, so $g$ is indeed a quasi-isomorphism.
\end{proof}

As we will see in Section \ref{sect:higher-holonomy}, and is already hinted at in Section \ref{sect:riemann-hilbert-functor}, the higher parallel transport $\tau$ of an $\Ascr$-connection $D$ on $E$, for $\Ascr$ a higher Lie algebroid, satisfies a special property: if $\sigma\in\Pi_\infty(\Ascr)_1$, then $\tau_\sigma:E_{\sigma_0(1)}\to E_{\sigma_0(0)}$ is an isomorphism.

\begin{definition} \label{def:higher-general-linear-and-gauge-groupoid}
    The {\bf higher general linear groupoid} $\GL_\infty(E) \subseteq \End_\infty(E)$ of $E$ is the open simplicial submanifold of $\End_\infty(E)$ given by
    \begin{equation*}
        \GL_\infty(E)_k :=\left\{A \in \Endi(E)_k \sep  \text{$A_\alpha$ is an isomorphism for all $\alpha:[1]\hookrightarrow [k]$}\right\}.
    \end{equation*}
    Furthermore, we define the \textbf{higher gauge groupoid of $E$} as $\Gaui(E) := \GLi(E) \cap \Neri(E)$.
\end{definition}

Notice that $\GLi(E)$ is the simplicialization of the groupoidal globular data for $E$ from Definition \ref{def:globular-collections-on-E}. 

\subsection{The transitive case}
Let us now restrict our attention to ruths of transitive higher algebroids. As we will explain in this section, the image of the parallel transport functor in this case is severely restricted and will actually land in a higher Lie groupoid without any quasi-smoothness considerations.

In this section we always assume that $M$ are connected.

\begin{definition}
Let $\Ascr\to M$ be a higher Lie algebroid. Let $A\to M$ be the vector bundle such that $C^1(\Ascr)=\Gamma(A)$. Then the cohomological vector field $Q:C^0(\Ascr)\to C^1(\Ascr)$ defines a vector bundle map $\rho:A\to TM$, which we call the \textbf{anchor map} of $\Ascr$. We say that $\Ascr$ is \textbf{transitive} if $\rho$ is fiberwise surjective.

If $G$ is a higher Lie groupoid over $M$, by which we mean that $G_0=M$, we say that it is \textbf{transitive} if for all $x,y\in M$ there is $g\in G_1$ such that $v_0(g)=x$ and $v_1(g)=y$.
\end{definition}

Let $\Ascr\to M$ be a transitive higher Lie algebroid and $D$ a ruth of $\Ascr$ on a graded vector bundle $E=\bigoplus_{k=0}^{n-1}E^k$. Let $\delta\in\End^1(E)$ be the induced differential on $E$. For any path $\sigma\in\Pi_\infty(\Ascr)_1$, parallel transport (cf. Theorem \ref{thm:riemann-hilbert-higher-algebroids} and Definition \ref{def:Sigma}) gives an isomorphism $\tau(\sigma):E_{\sigma_0(1)}\stackrel{\cong}{\to} E_{\sigma_0(0)}$ commuting with the differentials. Hence,
\[
\delta_{\sigma_0(0)}=\tau(\sigma)\delta_{\sigma_0(1)}\tau(\sigma)^{-1},
\]
so the isomorphism type of $\delta$ is fixed.

Now, fix $x\in M$ and $\delta_x \in \End^1(E_x)$, and write
\begin{equation*}
    \End^1(E,\delta_x) := \left\{\delta_y\in\End^1(E_y) \sep \,\begin{matrix}\text{with $y\in M$ and $g \in \GL^0(E_x, E_y)$} \\
    \text{such that $g\delta_x = \delta_yg$}\end{matrix}\right\},
\end{equation*} 
for the space of ``differentials'' (they might not square to zero) with the same isomorphism type as $\delta_x$. Furthermore, we write
\begin{equation*}
    \GLi(E, \delta_x)_k := \{ A \in \GLi(E)_k : A_\alpha \in \End^1(E, \delta_x) \text{ for all } \alpha : [0] \hookrightarrow [k]\}
\end{equation*}
and $\Gaui(E, \delta_x) = \GLi(E, \delta_x) \cap \Neri(E)$.

\begin{lemma}
    Let $V$ be a graded vector space with a degree $1$ map $\delta$. Then, $\End^1(V, \delta) \subseteq \End^1(V)$ is an embedded submanifold.
\end{lemma}
\begin{proof}
    Consider the $\GL^0(V)$ action on $\End^1(V)$ given by
    \begin{equation*}
        (g, \delta) \mapsto g\cdot \delta := g \delta g^{-1},
    \end{equation*}
    The statement of the lemma is equivalent to saying that all the orbits of this action are embedded submanifolds. Recall that orbits of an action are always immersed submanifolds. Thus, it suffices to show that the orbits of this action are locally closed. This follows by applying \cite[Prop. 1.65]{milne2017algebraic} to the action map once $\delta$ is fixed. Also see the text at the start of Section C, Chapter 7.
\end{proof}

\begin{lemma}
    The set $\End^1(E,\delta_x)$ is a embedded submanifold of $\End^1(E)$.
\end{lemma}
\begin{proof}
    Write $V := E_x$ and $\delta := \delta_x$. Upon choosing a degree preserving trivialization, we obtain a commutative diagram  
    \begin{center}
        \begin{tikzcd}
            \End^1(E|_U, \delta_x) \ar[d, "\cong"'] \ar[r, hook] & \End^1(E|_U) \ar[d, "\cong"] \\ 
            U \times \End^1(V, \delta) \ar[r, hook ] & U \times \End^1(V).
        \end{tikzcd}
    \end{center}
    By the previous lemma, the space $\End^1(V, \delta)$ is an embedded submanifold of $\End^1(V)$, so $\End^1(E, \delta_x)$ is an embedded submanifold of $\End^1(E)$.
\end{proof}

\begin{corollary}
    The manifold $\GLi(E, \delta_x)_k$ is an embedded submanifold of $\GLi(E)_k$.
\end{corollary}

\begin{lemma} \label{lem:gau(E,delta_x)-is-embedded}
     Let $\delta_x^2=0$. Then for each $k \geq 0$, $\Gaui(E, \delta_x)_k$ is an embedded submanifold of $\GL_\infty(E, \delta_x)_k$. 
\end{lemma}
\begin{proof}
    Let us first construct an embedded submanifold $M(E, \delta_x)_k \subseteq \GLi(E, \delta_x)_k$. We will then show that this embedded submanifold coincides with $\Gaui(E, \delta_x)_k$. Write 
    \[
    T_k := \pi_0^*\End^1(E, \delta_x) \oplus\bigoplus_{\substack{\alpha : [i] \hookrightarrow [k] \\ i > 0, \, \alpha(0) = 0}} \Hom^{1-i}(\pi^*_{\alpha(i)}E, \pi^*_{0}E)\to M^{k+1}.
    \]
    For $j, \ell >0$ and $A\in T_k$, define $\iota(A)_{j\ell}\in \Hom^0(E_{v_\ell(A)},E_{v_j(A)})$ by
    \[
        \iota(A)_{j\ell} := A_{0j}^{-1}(A_{0\ell} + A_0 A_{0j\ell} + A_{0j\ell}A_{0\ell}^{-1}A_0A_{0\ell}),
    \]
    where $A_{0j\ell}$ refers to $\alpha:[2]\hookrightarrow [k]$ with $\alpha(0)=0$, $\alpha(1)=j$ and $\alpha(2)=\ell$, and similarly for $A_{0j}$. We define
    \[
    M(E, \delta_x)_k := \{ A \in T_k :  \iota(A)_{j\ell} \text{ is an isomorphism for all $j, \ell > 0$}\}
    \]
    Note that being an isomorphism is an open condition, so $M(E, \delta_x)_k$ is an open subset of $T_k$. Let us write 
    \[
    S = \{\alpha : [i] \hookrightarrow [k] : \alpha(0) = 0\},
    \]
    where $\alpha : [0] \hookrightarrow [k]$ with $\alpha(0) = 0$ is used to index the $\pi_0^*\End^1(E, \delta_x)$-elements inside $M(E, \delta_x)_k$. If $\alpha:[i]\hookrightarrow[k]$ is such that $\alpha(0) > 0$, then we write $\alpha^+ : [i + 1] \hookrightarrow [k]$ for $\alpha^+(0) = 0$ and $\alpha^+(m + 1) = \alpha(m)$ for $m \in [i]$.
    
    We will now define
    \begin{align*}
        \iota : M(E, \delta_x)_k & \longrightarrow  \GLi(E, \delta_x)_k 
    \end{align*}
    inductively. For $\alpha \in S$, we define
    \begin{equation*}
        \iota(A)_\alpha = A_\alpha.
    \end{equation*}
    For $\alpha : [i] \hookrightarrow [k]$ with $\alpha(0) > 0$, we proceed by induction on $i$. For $i=0$ and $\alpha(0) = j>0$, we set
    \begin{equation*}
         \iota(A)_{\alpha} =  A_{{0j}}^{-1}A_{{0}}A_{{0j}}.
    \end{equation*}
    For $i=1$ and $\alpha = j\ell$, with $ j, \ell  >0$, we define
    \begin{equation*}
        \iota(A)_\alpha = A_{{0j}}^{-1}(A_{{0\ell}} + A_{{0}}A_{{0j\ell}} + A_{{0j\ell}} \iota(A)_{{\ell}}).
    \end{equation*}
    Note that this coincides with the definition of $\iota(A)_{j\ell}$ we gave before. Moreover, the condition on $\iota(A)_{j\ell}$ inside the definition of $M(E, \delta_x)_k$ is precisely what guarantees that the map $\iota$ lands inside $\GLi(E, \delta_x)_k$. Finally, for all $\alpha$ with $i>1$ and $\alpha(0) > 0$, we define
    \begin{equation*}
        \iota(A)_{\alpha} = A_{\alpha^+ \circ \alpha^1}^{-1}\left(-\sum_{j = 1}^{i}(-1)^jA_{\alpha^+ \circ \delta^j} + A_{{0}} A_{\alpha^+} + \sum_{j > 1}^{i+1} (-1)^j A_{\alpha^+ \circ \alpha^j} \iota(A)_{\alpha^+ \circ \beta^{i -j }}\right).
    \end{equation*}
    Notice that the projection back onto the $S$-components of $\GLi(E, \delta_x)_k$ defines a retract of $\iota$, so $\iota$ is indeed an embedding. 
    
    Let us now show that the image of $\iota$ is $\Gaui(E, \delta_x)_k$. Let us first show that $K(\iota(A))_\alpha = 0$ for all $\alpha : [i] \hookrightarrow [k]$. 
    Note that for $\alpha = 0$, the curvature equation holds, as for all $\delta \in \End^1(E, \delta_x)$, we have $\delta^2 = 0$. Now, if $\alpha(0) = 0$ and $i > 0$, then we note that $\alpha = \beta^+$ for some $\beta : [i-1] \hookrightarrow [k]$. Then, as $\iota(A)_\beta$ is defined by solving the curvature equation \eqref{eq:curvature_end} of $\beta^+ = \alpha$ for $\beta$, it follows that $K(\iota(A))_\alpha = 0$.
    
    For $i=0$ and $\alpha = j>0$, we have
    \begin{equation*}
        \iota(A)_{{j}} \iota(A)_{{j}} = A_{{0j}}^{-1}A_{{0}}^2A_{{0j}} = 0.
    \end{equation*}
    
    For $i=1$ $\alpha = j\ell$ with $j, \ell > 0$, we have
    \begin{align*}
        \iota(A)_{j}\iota(A)_{j\ell} - \iota(A)_{j\ell} \iota(A)_{\ell} & = A_{{0j}}^{-1}A_{{0}}(A_{{0\ell}} + A_{{0}}A_{{0j\ell}} + A_{{0j\ell}} \iota(A)_{{\ell}}) + \\
        & \qquad - A_{{0j}}^{-1}(A_{{0\ell}} + A_{{0}}A_{{0j\ell}} ) A_{{0\ell}}^{-1}A_{{0}}A_{{0\ell}} \\
        & =  A_{{0j}}^{-1}A_{{0}}A_{{0j\ell}} \iota(A)_{{\ell}} - A_{{0j}}^{-1}(A_{{0}}A_{{0j\ell}} ) A_{{0\ell}}^{-1}A_{{0}}A_{{0\ell}} \\
        &=0.
    \end{align*}
    
    Finally, let $i>1$ and $\alpha = r_0 \dots r_i$ with $r_0, \dots, r_i > 0$. Firstly, we note that
    \begin{align*}
        \sum^{i-1}_{j = 1}(-1)^j & A_{0r_0}^{-1}\sum_{\ell = 1}^j (-1)^\ell A_{0r_0\dots \hat{r}_{\ell - 1} \dots \hat{r}_j \dots r_i} +  \sum^{i-1}_{j = 1}(-1)^j A_{0r_0}^{-1}\sum_{\ell = {j + 1}}^{i-1} (-1)^\ell A_{0r_0\dots \hat{r}_j \dots \hat{r}_{\ell} \dots  r_i} =\\
& =  A_{0r_0}^{-1} \sum^{i-1}_{j = 1}\left(\sum_{\ell = 1}^j (-1)^{j + \ell} A_{0r_0\dots \hat{r}_{\ell - 1} \dots \hat{r}_j \dots r_i} +  \sum_{\ell = {j + 1}}^{i-1} (-1)^{\ell + j} A_{0r_0\dots \hat{r}_j \dots \hat{r}_{\ell} \dots  r_i}\right)\\
& =  A_{0r_0}^{-1} \sum^{i-1}_{j = 1}\left(\sum_{\ell = 0}^{j-1} (-1)^{j + \ell + 1} A_{0r_0\dots \hat{r}_{\ell} \dots \hat{r}_j \dots r_i} +  \sum_{\ell = {j + 1}}^{i-1} (-1)^{\ell + j} A_{0r_0\dots \hat{r}_j \dots \hat{r}_{\ell} \dots  r_i}\right)\\
& =  A_{0r_0}^{-1} \left(- \sum_{0 \leq \ell < j \leq {i -1}} (-1)^{\ell + j} A_{0r_0\dots \hat{r}_{\ell} \dots \hat{r}_j \dots r_i} +  \sum_{1 \leq j < \ell \leq i -1} (-1)^{\ell + j} A_{0r_0\dots \hat{r}_j \dots \hat{r}_{\ell} \dots  r_i}\right)\\
& =  -A_{0r_0}^{-1}\sum_{j = 1}^{i-1} (-1)^{ j} A_{0r_1\dots \hat{r}_j \dots r_i} 
    \end{align*}
Then, using this and the curvature equation $K_E(\iota(A))_\beta = 0$ for $\beta:[j]\hookrightarrow[k]$ with $j < i$ and $\beta(0) = 0$, we obtain
    \begin{align*}
        K_E(& \iota(A))_{r_0 \dots r_i} \\
        & = - \sum^{i-1}_{j = 1}(-1)^j\iota(A)_{r_0 \dots \hat{r}_j\dots r_i} + \sum_{j=0}^i(-1)^j\iota(A)_{r_0 \dots r_j}\iota(A)_{r_j\dots r_i} \\
& = \sum^{i-1}_{j = 1}(-1)^j A_{0r_0}^{-1}\sum_{\ell = 1}^j (-1)^\ell A_{0r_0\dots \hat{r}_{\ell - 1} \dots \hat{r}_j \dots r_i} +  \sum^{i-1}_{j = 1}(-1)^j A_{0r_0}^{-1}\sum_{\ell = {j + 1}}^{i-1} (-1)^\ell A_{0r_0\dots \hat{r}_j \dots \hat{r}_{\ell} \dots  r_i}\\
& \qquad - A_{0r_0}^{-1} A_0 \sum_{j=1}^{i-1} (-1)^j A_{0r_0 \dots \hat{r}_j \dots r_i}  + A_{0r_0}^{-1}\sum_{j = 1}^{i-1} \sum_{\ell > j}^i(-1)^{j+\ell}A_{0r_0 \dots \hat{r}_j \dots r_\ell}\iota(A)_{r_\ell \dots r_i}\\
& \qquad - A_{0r_0}^{-1}\sum_{j = 1}^{i-1} \sum_{\ell = 1}^{j-1}(-1)^{j+\ell}A_{0r_0 \dots r_\ell}\iota(A)_{r_\ell \dots \hat{r}_j \dots r_i} \\
&\qquad + \sum_{j = 0}^i (-1)^jA^{-1}_{0r_0}\left(- \sum^{j-1}_{\ell = 0} (-1)^\ell A_{0r_0 \dots \hat{r}_\ell \dots r_j}+ A_0 A_{0r_0 \dots r_j}\right. \\
& \qquad \left. + \sum_{\ell = 1}^j (-1)^\ell A_{0r_0 \dots r_\ell}\iota(A)_{r_\ell \dots r_j}\right) \iota(A)_{r_j\dots r_i} \\
& = A_{0r_0}^{-1}\left( -\sum_{j = 1}^{i-1} (-1)^{j} A_{0r_1\dots \hat{r}_j \dots r_i} - \color{red} A_0 \sum_{j=1}^{i-1} (-1)^j A_{0r_0 \dots \hat{r}_j \dots r_i}\right. \\
& \qquad + \sum_{1 \leq j < \ell \leq i}(-1)^{j+\ell}A_{0r_0 \dots \hat{r}_j \dots r_\ell}\iota(A)_{r_\ell \dots r_i} - \color{cyan} \sum_{1 \leq \ell < j \leq i-1}(-1)^{j+\ell}A_{0r_0 \dots r_\ell}\iota(A)_{r_\ell \dots \hat{r}_j \dots r_i} \\
& \qquad - \sum_{0 \leq \ell < j \leq i} (-1)^{\ell +j } A_{0r_0 \dots \hat{r}_\ell \dots r_j} \iota(A)_{r_j\dots r_i} + \color{red} \sum_{j = 0}^i (-1)^{j }  A_0 A_{0r_0 \dots r_j} \iota(A)_{r_j\dots r_i}\\
& \qquad \left. + {\color{cyan}{\sum_{1 \leq \ell \leq j \leq i } (-1)^{\ell +j} A_{0r_0 \dots r_\ell}\iota(A)_{r_\ell \dots r_j}\iota(A)_{r_j\dots r_i}}}\right) \\
& = A_{0r_0}^{-1}\left( -\sum_{j = 1}^{i-1} (-1)^{j} A_{0r_1\dots \hat{r}_j \dots r_i} - \sum_{1 \leq j \leq i} (-1)^{j } A_{0r_1 \dots r_j} \iota(A)_{r_j\dots r_i}\right. \\
& \qquad \color{red} - A_0 \cancelto{0}{K_E(\iota(A))_{0 r_1 \dots r_i}} + \cancelto{0}{A_0^2}A_{0r_0\dots r_i} + A_0A_{0r_1\dots r_i}\\
& \qquad \left.+ {\color{cyan}{\sum_{\ell = 1}^{i}(-1)^\ell A_{0r_0 \dots r_\ell}\cancelto{0}{K_E(\iota(A))_{r_\ell \dots r_i}}}}\right) \\
& = A_{0r_0}^{-1}\cancelto{0}{K_E(\iota(A))_{0r_1\dots r_i}} \\
&= 0.
    \end{align*}
Here, the colours indicate what terms combine to give which other terms.

    This shows that $\iota(M(E, \delta_x)_k) \subseteq \Gaui(E, \delta_x)_k$. On the other hand, we note that due to the curvature equations all components of $A \in \Gaui(E, \delta_x)$ are unique fixed by the components $A_\alpha$ with $\alpha(0) = 0$. So, we indeed obtain an equality. 
\end{proof}

\begin{proposition}
    Let $\delta_x \in \End^1(E)$ such that $\delta_x^2=0$. Then the manifolds $\Gaui(E, \delta_x)_k \subseteq \GL_\infty(E)_k$ define a higher Lie groupoid. 
\end{proposition}
\begin{proof}
    Note that $\Gaui(E, \delta_x)$ is constructed using the following steps. We first restrict the simplicial set $\Endi(E)$ to the 0-simplices isomorphic to $\delta_x$, which gives the simplicial subset $\Endi(E, \delta_x)\subseteq \Endi(E)$. Then we further restrict $\Endi(E, \delta_x)$ to the 1-simplices such that the $\id : [1] \to [1]$ component is an isomorphism, giving the simplicial subset $\GLi(E, \delta_x)\subseteq\Endi(E)$. We then intersect this with the simplicial subset $\Neri(E)\subseteq\Endi(E)$, i.e.
    \[
    \Gaui(E,\delta_x)=\GLi(E,\delta_x)\cap\Neri(E),
    \]
    which is again a simplicial subset of $\Endi(E)$.
    
    In the proof of Lemma \ref{lem:gau(E,delta_x)-is-embedded} we constructed $M(E, \delta_x)$ based on $\alpha : [i] \hookrightarrow [k]$ with $0 \in \im(\alpha)$. However, this choice was arbitrary: we could have also constructed $M(E, \delta_x)$ around any other $j \in [k]$. For the rest of this proof, we will focus on the case $j = 0$, since the other cases are done similarly.
    
    Note that $d_i : M(E, \delta_x)_k \to M(E, \delta_x)_{k-1}$ for $i > 0$ is given by 
    \begin{equation*}
        A \mapsto (A_{\delta^i\circ \alpha})_{\alpha \in S},
    \end{equation*}
    i.e., we project $A \in M(E, \delta_x)_k$ onto all its components that factor through $\delta^i$. Now, let us write 
    \begin{equation*}
        T_{k,0} := \pi_0^*\End^1(E, \delta_x) \oplus\bigoplus_{\substack{\alpha : [i] \hookrightarrow [k] \\ i > 0, \, \alpha(0) = 0 \\ \alpha \neq \id}} \Hom^{1-i}(\pi^*_{\alpha(i)}E, \pi^*_{0}E)\to M^{k+1}.
    \end{equation*}
    Then we define
    \begin{equation*}
        M(E, \delta_x)_{k, 0} := \{ A \in T_{k,0} :  \iota(A)_{j\ell} \text{ is an isomorphism for all $j, \ell > 0$}\},
    \end{equation*}
    where $\iota(A)_{j\ell}$ is as in the proof of Lemma \ref{lem:gau(E,delta_x)-is-embedded}. This embeds into $\prod_{i > 0}^k M(E, \delta_x)_{k -1}$ through the map 
    \begin{equation*}
        A \mapsto (i, A_{\delta^i\circ \alpha}).
    \end{equation*}
    A proof analogous to the one for Lemma \ref{lemma:horn_of_end} shows that the image of this map coincides with the $(k, 0)$-horn of $\Gaui(E, \delta_x)$. Then, the $(k, 0)$-horn map $M(E, \delta_x)_k \to M(E, \delta_x)_{k, 0}$ is the surjective submersion that forgets the component of $\id:[k]\hookrightarrow[k]$. An analogous story shows that $(k, \ell)$-horn map is a surjective submersion for $\ell > 0$. Thus, $\Gaui(E, \delta_x)$ defines a higher Lie groupoid. 
\end{proof}

\begin{corollary}
    Let $G$ be a transitive higher groupoid over $M$ and let $\sum_{i\geq 0}F_i$ be a ruth on a bounded graded vector bundle $E$ such that $F_1$ acts by isomorphisms, i.e., $F_1(g): E_{v_1(g)}\to E_{v_0(g)}$ is an isomorphism for all $g\in G_1$. Then, the functor $F: G \to \GLi(E)$ from Theorem \ref{thm:groupoid-connections-as-functors} factors through a smooth functor $F:G\to \Gaui(E, \delta_x)$, where $\delta_x:=F_0(x)$ for any $x\in M$:
    \[\begin{tikzcd}
	G & {\GLi(E)} \\
	& {\Gaui(E,\delta_x).}
	\arrow["F", from=1-1, to=1-2]
	\arrow["F"', from=1-1, to=2-2]
	\arrow[hook, from=2-2, to=1-2]
\end{tikzcd}\]

    In particular, this holds for the parallel transport functor $\tau :\Pi_\infty(\Ascr) \to \GL_\infty(E)$ coming from a ruth of a transitive higher Lie algebroid $\Ascr$.
\end{corollary}
\begin{proof}
    The fact that $(F_1)_g$ defines an isomorphism means that we may restrict the image to $\GL_\infty(E)$. Then, since $D$ is a ruth and $G$ is transitive, the isomorphism type of the $C^0_G(\End^1(E))$ part of $D$ is fixed by the equation
    \begin{equation*}
        \delta_{t(g)} (F_1)_g - (F_1)_g \delta_{s(g)} = 0.
    \end{equation*}
    So, $F$ factors through $\GLi(E, \delta_x)$. Then, as $D$ is a ruth, $K_E \circ F = 0 \circ p$. Thus, $F$ factors through $\Gaui(E)$.

    For a transitive algebroid $\Ascr$, the higher groupoid $\Pi_\infty(A)$ is transitive. Furthermore, as  argued at the beginning of this section, $(F_1)_{\sigma_1}$ for parallel transport is an isomorphism for every $\sigma_1 : T\Delta^1 \to \Ascr$, since it is locally given by the usual parallel transport of vector bundles along a connection. Finally, by theorem \ref{thm:riemann-hilbert-higher-algebroids}, the integration of a ruth again defines a ruth. So, we can apply the first part of this corollary to obtain the desired result. 
\end{proof}

\bigskip
\section{Higher holonomy} \label{sect:higher-holonomy}

\subsection{Higher parallel transport}

Throughout this section $(\Ascr\to M,Q)$ is a higher Lie algebroid and $E\to M$ is a graded vector bundle.

\begin{definition}
Let $D$ be an $\Ascr$-connection on $E$. We define the \textbf{higher parallel transport} of $D$ as $\tau:=\RH(D)\in C_{\Pi_\infty(\Ascr)}(\End(E))^1$. 
\end{definition}

As an immediate consequence of Theorem \ref{thm:groupoid-connections-as-functors} we get the following.

\begin{corollary}[{of Theorem \ref{thm:groupoid-connections-as-functors}}] \label{cor:higher-parallel-transport}
Let $D$ be an $\Ascr$-connection on $E$. Then the higher parallel transport of $D$ defines a smooth functor
\[
\tau: \Pi_\infty(\Ascr)\to \GLi(E).
\]
Moreover, if $D$ is flat then the higher parallel transport defines a quasi-smooth functor
\[
\tau: \Pi_\infty(\Ascr)\to \Gaui(E).
\]
\end{corollary}

\begin{proof}
The set-theoretic existence follows from Theorem \ref{thm:groupoid-connections-as-functors}. To see that it is smooth, let $p:U\to\Pi_\infty(\Ascr)_k$ be a plot. Then, write $\pi_1 : U \times \Delta^k \to U$ for the projection onto the first term. Then, $\ker\pi_{1*} = 0_U \boxplus T\Delta^k$ defines a Lie algebroid over $U \times \Delta^k$. Let $\nabla$ be a degree-preserving $0_U \boxplus T\Delta^k$-connection on $\hat p_0^*E\to U\times\Delta^k$, and define
\[
\omega:=\hat p_0^*D-d^\nabla\in\Omega(0_U\boxplus T\Delta^k,\End(\hat p_0^*E))^1.
\]
Then $\tau(p(x))$ is fully determined by $\RH(D)(p(x))=\hat\Sigma^{\iota_x^*\nabla}(\iota_x^*\omega)$, where $\iota_x:\Delta^k\to U\times\Delta^k$ is the inclusion $\iota_x(t):=(x,t)$. For any path $\gamma$ in $\Delta^k$, we have that $\tau^{\iota_x^*\nabla}_\gamma = \tau^\nabla_{(x,\gamma)}$ depends smoothly on $x$. Furthermore, the form $\omega$ depends smoothly on $x$. We conclude that $\hat\Sigma^{\iota_x^*\nabla}(\iota_x^*\omega)$ is an integral over a compact set with an integrand that depends smoothly on $x$. 
\end{proof}

According to Theorem \ref{thm:simplicialization}, if $E$ has length $n$ then $\GLi(E)$ and $\Curvi(E)$ are Lie $(n+1)$-groupoids, and $\Gaui(E)$ is a quasi-smooth $n$-groupoid, in which case we denote them by $\GL_{n+1}(E)$, $\Curv_{n+1}(E)$ and $\Gau_n(E)$. Hence, higher parallel transport is a functor $\tau:\Pi_\infty(\Ascr)\to\GL_{n+1}(E)$ or $\tau:\Pi_\infty(\Ascr)\to\Gau_n(E)$. We would like to go even further: in this case also the source of parallel transport descends from $\Pi_\infty(\Ascr)$ to $\Pi_{n+1}(\Ascr)$ or $\Pi_n(\Ascr)$. These homotopy groupoids can be defined abstractly using higher homotopy categories \cite[\href{https://kerodon.net/tag/0596}{Subsect. 0596}]{kerodon}, but here we opt for a more down-to-earth approach.

\begin{definition}
The \textbf{$n$-skeleton} of the standard $k$-simplex $\sk_n(\Delta^k)$ is the space of points $(t_1,\dots,t_k)\in\Delta^k$ such that there are $0\leq i_1<\dots< i_{k-n}\leq k$ with $t_{i_j}=t_{i_j+1}$ for all $j=1,\dots,k-n$, where $t_0:=1$ and $t_{k+1}:=0$.
\end{definition}

The $n$-skeleton of $\Delta^k$ is the union of all its $n$-dimensional faces. So, we have that $\sk_n(\Delta^k)=\Delta^k$ for all $n\geq k$, and $\sk_{k-1}(\Delta^k)=\del\Delta^k$. In particular, we see that $\sk_{n}(\Delta^k)$ is not an embedded submanifold for $k > n$. However, each face is an $n$-dimensional manifold and the intersection of two distinct $n$-faces will either be empty or an $m$-dimensional face with $0 \leq m \leq n-1$. 
Finally, observe that $\sk_{n-1}(\Delta^k)\subseteq \sk_n(\Delta^k)$.

\begin{definition}
    Let $U$ be a manifold. A \textbf{higher Lie algebroid map} $\sigma : T\sk_n(\Delta^k) \to \mathscr{A}$ is defined as a collection of higher Lie algebroid maps, namely for each $m \leq n$ and $m$-dimensional face $f \subseteq \Delta^k$, we have a higher Lie algebroid map
    \begin{equation*}
        \sigma_{f}: Tf \to \mathscr{A},
    \end{equation*}
    such that if $f_1 \subseteq f_2$, we have $\sigma_{f_1} = \sigma_{f_2}|_{Tf_1}$.
\end{definition}

\begin{remark}
     Given a Lie algebroid $A$, we can also form the products $T\sk_n(\Delta^k) \times A$. Then, a higher Lie algebroid map $\sigma : T\sk_n(\Delta^k) \times A \to \mathscr{A}$ is a collection of higher algebroid maps
     \begin{equation*}
         \sigma_f : Tf \times A \to \mathscr{A}
     \end{equation*}
     for each $m$-dimensional face $f$ of $\Delta^k$ with $m \leq n$.
\end{remark}

\begin{definition} \label{def:homotopy-n-groupoid}
The \textbf{homotopy $n$-groupoid $\Pi_n(\Ascr)$ of $\Ascr$} is defined by setting $\Pi_n(\Ascr)_k$ to be the space of higher Lie algebroid maps $\sigma:T\sk_n(\Delta^k)\to\Ascr$ that can be extended to a map from $T\sk_{n+1}(\Delta^k)$, modulo the relation $\sim_k$ defined as follows. Two such maps $\sigma_0$ and $\sigma_1$ are equivalent if they are homotopic relative $\sk_{n-1}(\Delta^k)$, meaning that there is a higher Lie algebroid map $H:T(\sk_n(\Delta^k)\times\Delta^1)\to\Ascr$ such that 
\begin{enumerate}[label=(\roman*)]
    \item $\iota_i^*H=\sigma_i$, where $\iota_i:\sk_n(\Delta^k)\times\{i\}\hookrightarrow \sk_n(\Delta^k)\times\Delta^1$ are the inclusions,
    \item and $j^*H$ factors through the projection $\sk_{n-1}(\Delta^k)\times\Delta^1\to \sk_{n-1}(\Delta^k)$, where $j:\sk_{n-1}(\Delta^k)\times\Delta^1\hookrightarrow \sk_n(\Delta^k)\times\Delta^1$ is the inclusion.
\end{enumerate}
\end{definition}

Observe that, as expected:
\[
\Pi_n(\Ascr)_k =\begin{cases}
    \Pi_\infty(\Ascr)_k, & \text{if $k<n$}, \\
    \Pi_\infty(A)_n/\sim_n, & \text{if $k=n$},
\end{cases}
\]
where $\sigma_0,\sigma_1\in\Pi_\infty(\Ascr)_n$ are equivalent if they are homotopic relative $\del\Delta^n$. Each $\Pi_n(\Ascr)_k$ is naturally a diffeological space by first endowing the space of Lie algebroid maps $\sigma : T \sk_n(\Delta^k) \to \Ascr$ that can be extended to a map from $T\sk_{n+1}(\Delta^k)$ with the pullback diffeology (Example \ref{ex:pullback}) and then endowing their quotient by the relation $\sim_k$ with the pushforward diffeology (Example \ref{ex:pushforward}).

Very explicitly, this means that the plots of $\Pi_n(\Ascr)_k$ are given by maps $p : U \to \Pi_n(\mathscr{A})_k$ such that for every $x\in U$ there is an open neighborhood $V\subseteq U$ of $x$ and a higher Lie algebroid map $q:0_V\boxplus T\sk_n(\Delta^k)\to\Ascr$ that can be extended to a higher Lie algebroid map $0_V\boxplus T\sk_{n+1}(\Delta^k)\to\Ascr$ such that $p|_V=[q]$, where $[\cdot]$ represents the equivalence class with respect to $\sim_k$.

\begin{lemma} \label{lem:pi-n-diffeological-groupoid}
The simplicial structure of $\Pi_\infty(\Ascr)$ descends to $\Pi_n(\Ascr)$ making it into a diffeological $n$-groupoid.
\end{lemma}

\begin{proof}
To see that the simplicial structure of $\Pi_\infty(\Ascr)$ descends to $\Pi_n(\Ascr)$, first note that if $\gamma:[k]\to[j]$ is order-preserving and we consider the induced map $\gamma_*:T\Delta^k\to T\Delta^j$, then the following diagram commutes:
\[
\begin{tikzcd}
    T\sk_n(\Delta^k) & T\sk_n(\Delta^j) \\
    T\sk_{n+1}(\Delta^k) & T\sk_{n+1}(\Delta^j). 
    \ar["\gamma_*", from=1-1, to=1-2]
    \ar[hook, from=1-1, to=2-1]
    \ar["\gamma_*", from=2-1, to=2-2]
    \ar[hook, from=1-2, to=2-2]
\end{tikzcd}
\]
Hence, if $\sigma:T\sk_n(\Delta^j)\to\Ascr$ extends to a map $T\sk_{n+1}(\Delta^j)\to\Ascr$, then $\gamma^*\sigma:T\sk_n(\Delta^k)\to\Ascr$ extends to a map $T\sk_{n+1}(\Delta^k)\to\Ascr$. Moreover, if two such elements $\sigma_0,\sigma_1:T\sk_n(\Delta^j)\to\Ascr$ are homotopic relative $\sk_{n-1}(\Delta^j)$ through a homotopy $H:T(\sk_n(\Delta^j)\times\Delta^1)\to\Ascr$, then $H\circ (\gamma_*\times\id)$ defines a homotopy between $\gamma^*\sigma_0$ and $\gamma^*\sigma_1$ relative $\sk_{n-1}(\Delta^k)$.

Consider now the skeletal horn spaces
\[
\Lambda^n_{k,\ell}:=\Lambda_{k,\ell}\cap \sk_n(\Delta^k)\subseteq\Delta^k.
\]
Then a $(k,\ell)$-horn in $\Pi_n(\Ascr)$ is the same as a higher Lie algebroid map $T\Lambda^n_{k,\ell}\to \Ascr$ that can be extended to a higher Lie algebroid map $T\Lambda^{n+1}_{k,\ell}\to\Ascr$ modulo homotopy relative $\Lambda^{n-1}_{k,\ell}$. To see that $\Pi_n(\Ascr)$ is a higher groupoid, we have to show that the horn maps are surjective. Note that for $k < n$, we have $\Pi_n(\Ascr)_k = \Pi_\infty(\Ascr)_k$. Then, as $\Pi_\infty(\Ascr)_k$ is a higher groupoid, the $(k, \ell)$-horn maps for $\Pi_n(\Ascr)$ are surjective for $k < n$. 
Moreover, if $k>n+1$ then $\Lambda^n_{k,\ell}=\sk_n(\Delta^k)$, so $\Pi_n(\Ascr)$ is a $(n+1)$-groupoid. 

For $k = n$, we have $\Pi_n(\Ascr) = \Pi_\infty(\Ascr)/\sim_n$ and $\Pi_n(\Ascr)_{n, \ell} = \Pi_\infty(\Ascr)_{n, \ell}$. Since $\Lambda^{n, \ell} \subseteq \partial\Delta^n = \sk_{n-1}(\Delta^n)$ and the homotopy defining $\sim_n$ is relative $\sk_{n-1}(\Delta^n)$, we see that the $(n, \ell)$-horn map factors through $\Pi_n(\Ascr)_n$:
\begin{center}
    \begin{tikzcd}
        \Pi_\infty(\Ascr)_n \ar[rd, "d_{n,\ell}"] \ar[d] \\
        \Pi_\infty(\Ascr)_n/\sim_n \ar[r] & \Pi_\infty(\Ascr)_{n, \ell}.
    \end{tikzcd}
\end{center}
This factorization coincides with the $(n, \ell)$-horn map of $\Pi_n(\Ascr)$. Then, as $d_{n, \ell}$ in the above diagram is surjective, we see that its factorization is also surjective. 

Finally, in the case that $k=n+1$, we have that $\Lambda^n_{n+1,\ell}= \Lambda_{n+1,\ell}$. Given a higher Lie algebroid map $\lambda:T\Lambda_{n+1,\ell}\to\Ascr$, a horn filling is an extension of $\lambda$ by a higher Lie algebroid map $T\sk_n(\Delta^{n+1})=T\del\Delta^{n+1}\to\Ascr$ that can be extended to a map $T\Delta^{n+1}\to\Ascr$, modulo homotopy relative $\sk_{n-1}(\Delta^{n+1})$. Since $\Lambda_{n+1,\ell}$ is a strong deformation retract of $\Delta^{n+1}$, any two extensions of $\lambda$ to $\Delta^{n+1}$ are homotopic relative $\Lambda_{n+1,\ell}$, and since $\sk_{n-1}(\Delta^{n+1})\subseteq \Lambda_{n+1,\ell}$, we have that any two such extensions represent the same element in $\Pi_n(\Ascr)_{n+1}$, which finally implies that $\Pi_n(\Ascr)$ is actually an $n$-groupoid.
\end{proof}

\begin{example}
If $A$ is a Lie algebroid, then $\Pi_1(A)$ is its fundamental or Weinstein groupoid. It is well known that it is not always a Lie groupoid \cite{crainic-fernandes}, but Lemma \ref{lem:pi-n-diffeological-groupoid} shows that it is always a diffeological groupoid.
\end{example}

\begin{proposition} \label{prop:parallel-transport-functor-descends-to-n+1-homotopy-groupoid}
Let $D$ be an $\Ascr$-connection on $E$. If $E$ has length $n$, then the higher parallel transport descends to a smooth functor
\[
\tau:\Pi_{n+1}(\Ascr)\to \GL_{n+1}(E).
\]
Moreover, if $D$ is flat then the higher parallel transport defines a functor
\[
\tau: \Pi_n(\Ascr)\to \Gau_n(E).
\]
\end{proposition}

\begin{proof}
Since $E$ has length $n$, then $C^p_{\Pi_\infty(\Ascr)}(\End^q(E))\neq 0$ only for $-n<q<n$. Hence, if $\sigma\in\Pi_\infty(\Ascr)_k$ then $\tau(\sigma)$ depends only on the $m$-dimensional data of $\sigma$, for $0\leq m\leq n$, i.e., on the restriction of $\sigma$ to $\sk_n(\Delta^k)$. In particular, any two $\sigma_0,\sigma_1:T\sk_{n+1}(\Delta^k)\to\Ascr$ which are homotopic relative $\sk_n(\Delta^k)$ have the same image through $\tau$.

If now $D$ is flat, then $\tau:\Pi_{n+1}(\Ascr)\to \Gau_n(E)$. Let now $\sigma_0,\sigma_1:T\sk_n(\Delta^k)\to\Ascr$ be homotopic relative $\sk_{n-1}(\Delta^k)$ through a homotopy $H$. The only discrepancy of $\tau(\sigma_0)$ and $\tau(\sigma_1)$ can be along the $n$-faces of $\sk_n(\Delta^k)$, since $\sigma_0$ and $\sigma_1$ agree on $\sk_{n-1}(\Delta^k)$. Let $C$ be one such face, which we can write as the image of a map $\alpha:\Delta^n\hookrightarrow \sk_n(\Delta^k)$. Then $(\alpha\times\id)^*H$ provides a homotopy from $\alpha^*\sigma_0$ to $\alpha^*\sigma_1$ relative $\sk_{n-1}(\Delta^n)=\del\Delta^n$. This is the same as a simplex $\sigma:T\Delta^{n+1}\to\Ascr$ such that
\[
d_i\sigma = \begin{cases}
    \pi^*(\delta_i\times\id)\iota^*\sigma, & \text{if $i=0,\dots, n-1$}, \\
    \alpha^*\sigma_1, & \text{if $i=n$}, \\
    \alpha^*\sigma_0, & \text{if $i=n+1$},
\end{cases}
\]
where $\pi$ and $\iota$ refer to the maps
\[\begin{tikzcd}
	{\Delta^n} & {\Delta^{n+1}} \\
	{\Delta^{n-1}\times\{0\}} & {\Delta^n\times\{0\}.}
	\arrow["{\delta_i}", from=1-1, to=1-2]
	\arrow["\pi"', from=1-1, to=2-1]
	\arrow["{\delta_i\times\id}"', from=2-1, to=2-2]
	\arrow["{\iota}"', from=2-2, to=1-2]
\end{tikzcd}\]
Hence, any integral along any lower-dimensional face of $\sigma$ that contains some $d_i$ for some $i=0,\dots,n-1$ vanishes. Since $b_i\sigma=d_0\dots d_0\sigma$, so that most terms vanish, we finally have that
\[
0=K_E(\tau(\sigma))=-(-1)^n\tau(\alpha^*\sigma_1)-(-1)^{n+1}\tau(\alpha^*\sigma_0),
\]
so $\tau(\alpha^*\sigma_0)=\tau(\alpha^*\sigma_1)$.
\end{proof}

\subsection{Higher general linear groupoids as simplicial groupoids} \label{sect:1-gpd-structures-on-higher-general-linear-gpds}

In this section $E\to M$ will always be a graded vector bundle. As we will show now, the higher groupoid $\GLi(E)$ has more structure than just that of a higher Lie groupoid: it is an augmented simplicial Lie groupoid. The reason we are interested in viewing $\GLi(E)$ as a simplicial Lie groupoid comes from the fact that for classical holonomy groups, the group operation is composition of parallel transports, but the higher groupoid structure of $\GLi(E)$ does not reflect this composition. We will see that there is a way of composing elements in $\GLi(E)$ which will become crucial for our definition of higher holonomy in Section \ref{sect:higher-holonomy:simplicial-holonomy}.

\begin{definition}
The \textbf{augmented simplex category $\Delta_+$} has as objects the (possibly empty) linearly ordered sets and as morphisms the order-preserving maps. An \textbf{augmented simplicial set} is a functor $X:\Delta_+^{\mathrm{op}}\to\mathrm{Set}$. Equivalently, it is a simplicial set $X$ together with a set $X_{-1}$ and a map $d_0:X_0\to X_{-1}$ such that $d_0d_1=d_0d_0:X_1\to X_{-1}$.

A(n) \textbf{(augmented) simplicial Lie group(oid)} is a functor $G: \Delta^{\mathrm{op}}_{(+)}\to \mathrm{LieGrp(d)}$, where $\mathrm{LieGrp(d)}$ is the category of Lie group(oid)s with Lie group(oid) morphisms. A(n) \textbf{(augmented) simplicial (graded) vector bundle} is a functor $V:\Delta^{\mathrm{op}}_{(+)}\to \mathrm{(Gr)Vect}$, where $\mathrm{(Gr)Vect}$ is the category of (graded) vector bundles with (graded) vector bundle morphisms.
\end{definition}

If $G$ is a simplicial Lie groupoid, its base $N$ is a simplicial manifold. Notice that a simplicial vector bundle $V\to N$ is a simplicial Lie groupoid, where the source and target maps coincide. 

We will use the following simplified notation: if $\alpha:[i]\to[k]$ and $\beta:[j]\to[k]$, then we write $\alpha\subseteq\beta$ if $\alpha([i])\subseteq \beta([j])$.

\begin{definition} \label{def:1-gpd-structure-on-GLinfty}
For each $k>0$ we define a groupoid structure $\GLi(E)_k\rightrightarrows \GLi(E)_{k-1}$ as follows:
\begin{enumerate}
    \item the source map is $d_{k-1}$,
    \item the target map is $d_k$,
    \item composition is given as follows: if $A^0,A^1\in \GLi(E)_k$ are such that $d_kA^1=d_{k-1}A^0$, then we set $A^0A^1\in\GLi(E)_k$ to be
    \[
    (A^0A^1)_\alpha:=\begin{cases}
        A^0_\alpha, & \text{if $\alpha\subseteq\delta^k$,} \\
        A^1_\alpha, & \text{if $\alpha\subseteq \delta^{k-1}$,} \\
        A^0_\varepsilon A^1_\varepsilon, & \text{if $\alpha=\varepsilon$}, \\
        A^0_\alpha A^1_\varepsilon+A^1_\alpha, & \text{otherwise,}
    \end{cases}
    \]
    where $\varepsilon:=\beta^1:[1]\hookrightarrow[k]$ is given by $\varepsilon(0)=k-1$ and $\varepsilon(1) = k$,
    \item the unit map is $s_{k-1}$,
    \item the inverse is given as follows: if $A\in\GLi(E)_k$, then
    \begin{equation} \label{eq:1-groupoid-structures-on-gliE}
    (A^{-1})_{\alpha:[i]\hookrightarrow[k]}=\begin{cases}
        A_\alpha, & \text{if $\alpha\subseteq \delta^{k-1}\cap \delta^k$,} \\
        A_\varepsilon^{-1}, & \text{if $\alpha=\varepsilon$,} \\
        -A_\alpha A_\varepsilon^{-1}, & \text{if $\varepsilon\subsetneq \alpha$,} \\
        A_{\hat\alpha}, & \text{otherwise.}
    \end{cases}
    \end{equation}
    where $\hat\alpha:[i]\hookrightarrow[k]$ is
    \[
    \hat\alpha(j):=\begin{cases}
        \alpha(j), & \text{if $\alpha(j)\neq k-1,k$,} \\
        k, & \text{if $\alpha(j)=k-1$,} \\
        k-1, & \text{if $\alpha(j)=k$.}
    \end{cases}
    \]
\end{enumerate}
\end{definition}
Notice that $\hat\alpha$ is well defined whenever $\alpha$ satisfies the ``otherwise'' condition from \eqref{eq:1-groupoid-structures-on-gliE}, which means that $\alpha\subseteq \delta^{k-1}$ or $\alpha\subseteq \delta^k$, but $\alpha\nsubseteq \delta^{k-1}\cap \delta^k$. 

\begin{figure}
    \centering
    \includesvg[width=0.9\linewidth]{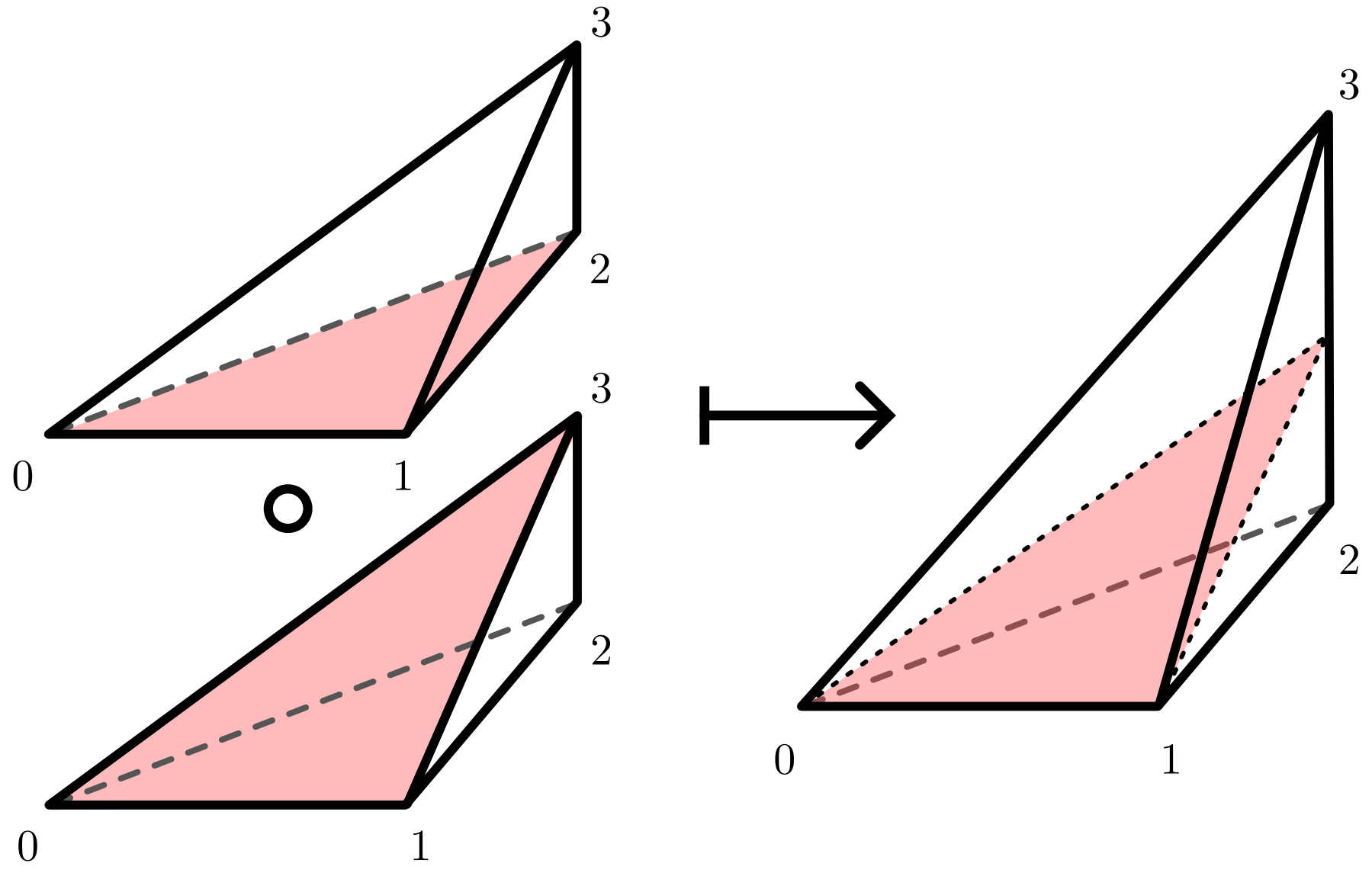}
    \caption{A geometric illustration of the composition inside $\GLi(E)_3 \rightrightarrows \GLi(E)_2$, where we compose two 3 simplices along the $\delta^3$ and $\delta^2$ face.}
    \label{fig:composition_glk}
\end{figure}

The geometric intuition for the above definition is that we glue the $n$-simplex along its top and bottom $n-1$-simplices. Figure \ref{fig:composition_glk} illustrates this for the 3-simplex.

\begin{remark}
    Note that the above construction is not particular to $\delta^{k-1}$ and $\delta^k$. It works for any $\delta^i$ and $\delta^j$. However, as we will see in the next section, this particular choice is compatible with our choice of cube-to-simplex families. 
\end{remark}

\begin{proposition} \label{prop:1-gpd-structure-on-GL-inf}
For each $k>0$, the structure $\GLi(E)_k\rightrightarrows \GLi(E)_{k-1}$ is a Lie groupoid. Moreover, for all $\gamma:[k]\to [j]$ order-preserving, the induced map
\[
\begin{tikzcd}
	{\GLi(E)_{j+2}} & {\GLi(E)_{k+2}} \\
	{\GLi(E)_{j+1}} & {\GLi(E)_{k+1}}
	\arrow["{\gamma[2]^*}", from=1-1, to=1-2]
	\arrow[shift right, from=1-1, to=2-1]
    \arrow[shift left, from=1-1, to=2-1]
	\arrow["{\gamma[1]^*}"', from=2-1, to=2-2]
	\arrow[shift right, from=1-2, to=2-2]
    \arrow[shift left, from=1-2, to=2-2]
\end{tikzcd}
\]
is a Lie groupoid morphism, where $\gamma[i]:[k+i]\to [j+i]$ is given by
\[
\gamma[i](p):=\begin{cases}
    \gamma(p), & \text{if $p\leq k$,} \\
    j+p-k, & \text{if $p>k$.}
\end{cases}
\]
\end{proposition}

\begin{proof}
The operations are clearly smooth, so let us check that $\GLi(E)_k\rightrightarrows\GLi(E)_{k-1}$ is indeed a groupoid.

Let $A^i\in\GLi(E)_k$, with $i=0,1,2$, such that $d_kA^2=d_{k-1}A^1$ and $d_kA^1=d_{k-1}A^0$, and let $\alpha:[i]\hookrightarrow [k]$. Then:
\begin{enumerate}
    \item if $\alpha\subseteq\delta^k$ we have that
    \[
    ((A^0A^1)A^2)_\alpha=(A^0A^1)_\alpha=A^0_\alpha=(A^0(A^1A^2))_\alpha,
    \]
    
    \item if $\alpha\subseteq \delta^{k-1}$, then
    \[
    ((A^0A^1)A^2)_\alpha=A^2_\alpha=(A^1A^2)_\alpha=(A^0(A^1A^2))_\alpha,
    \]

    \item if $\alpha=\varepsilon$, then
    \[
    ((A^0A^1)A^2)_\alpha=(A^0_\varepsilon A^1_\varepsilon)A^2_\varepsilon= A^0_\varepsilon (A^1_\varepsilon A^2_\varepsilon)=(A^0(A^1A^2))_\varepsilon,
    \]

    \item if $\alpha$ satisfies none of the above, then
    \begin{align*}
    ((A^0A^1)A^2)_\alpha &=(A^0A^1)_\alpha A^2_\varepsilon+A^2_\alpha=(A^0_\alpha A^1_\varepsilon+A^1_\alpha)A^2_\varepsilon+A^2_\alpha \\
    &=A^0_\alpha A^1_\varepsilon A^2_\varepsilon + A^1_\alpha A^2_\varepsilon+A^2_\alpha=A^0_\alpha (A^1A^2)_\varepsilon+(A^1A^2)_\alpha \\
    &=(A^0(A^1A^2))_\alpha.
    \end{align*}
\end{enumerate}
We conclude that the product is associative.

If now $A\in\GLi(E)_k$, let us study $AA^{-1}$. First of all, we have to check that $A$ and $A^{-1}$ are composable. For any $\beta:[i]\hookrightarrow [k-1]$ we have that $(d_{k-1}A)_\beta=A_{\delta^{k-1}\circ\beta}$ and that
\[
(d_kA^{-1})_\beta=(A^{-1})_{\delta^k\circ\beta}=\begin{cases}
    A_{\delta^k\circ\beta}, & \text{if $\beta\subseteq \delta^{k-1}$,} \\
    A_{\widehat{\delta^k\circ\beta}}, & \text{otherwise.}
\end{cases}
\]
If $\beta\subseteq \delta^{k-1}$, then $\delta^k\circ\beta=\delta^{k-1}\circ\beta$, and if $\beta\nsubseteq \delta^{k-1}$, then there is $j\in[i]$ such that $\beta(j)=k-1$, and hence
\[
\widehat{\delta^k\circ\beta}(p)=\begin{cases}
    \delta^k\circ\beta(p), & \text{if $p\neq j$} \\
    k, & \text{if $p=j$}
\end{cases}
=\delta^{k-1}\circ\beta(p).
\]
This settles that $A$ and $A^{-1}$ are composable.

We now compare $AA^{-1}$ and $s_{k-1}d_kA$. For $AA^{-1}$ we have that
\begin{enumerate}
    \item if $\alpha\subseteq \delta^{k-1}\cap \delta^k$, then
    \[
    (AA^{-1})_\alpha=A_\alpha=(A^{-1})_\alpha,
    \]

    \item if $\alpha\subseteq\delta^k$ but $\alpha \nsubseteq \delta^{k-1}\cap\delta^k$, then
    \[
    (AA^{-1})_\alpha=A_\alpha.
    \]

    \item if $\alpha\subseteq\delta^{k-1}$ but $\alpha \nsubseteq \delta^{k-1}\cap\delta^k$, then
    \[
    (AA^{-1})_\alpha=(A^{-1})_\alpha= A_{\hat\alpha},
    \]

    \item if $\alpha=\varepsilon$, then
    \[
    (AA^{-1})_\varepsilon=A_\varepsilon(A^{-1})_\varepsilon=A_\varepsilon A_\varepsilon^{-1}=\id,
    \]

    \item if $\varepsilon\subsetneq \alpha$, then
    \[
    (AA^{-1})_\alpha=A_\alpha(A^{-1})_\varepsilon+(A^{-1})_\alpha=A_\alpha A^{-1}_\varepsilon-A_\alpha A_\varepsilon^{-1}=0.
    \]
\end{enumerate}
As for $s_{k-1}d_kA$, we have that
\begin{enumerate}
    \item if $\alpha\subseteq \delta^{k-1}\cap \delta^k$, then $\delta^k\circ\sigma^{k-1}\circ\alpha=\alpha$, so
    \[
    (s_{k-1}d_kA)_\alpha=A_\alpha,
    \]

    \item if $\alpha\subseteq\delta^k$ but $\alpha \nsubseteq \delta^{k-1}\cap\delta^k$, then $\delta^k\circ\sigma^{k-1}\circ\alpha=\alpha$, so
    \[
    (s_{k-1}d_kA)_\alpha=A_\alpha,
    \]

    \item if $\alpha\subseteq\delta^{k-1}$ but $\alpha \nsubseteq \delta^{k-1}\cap\delta^k$, then $\delta^k\circ\sigma^{k-1}\circ\alpha=\hat\alpha$, so
    \[
    (s_{k-1}d_kA)_\alpha=A_{\hat\alpha},
    \]

    \item if $\alpha=\varepsilon$, then $\sigma^{k-1}\circ\varepsilon$ is constant, so
    \[
    (s_{k-1}d_kA)_\varepsilon=\id,
    \]

    \item if $\varepsilon\subsetneq \alpha$, then $\sigma^{k-1}\circ\alpha$ is not injective but $i>1$, so
    \[
    (s_{k-1}d_kA)_\alpha=0.
    \]
\end{enumerate}
This proves at last that $AA^{-1}=s_{k-1}d_kA$. An analogous argument shows that $A^{-1}A=s_{k-1}d_{k-1}A$.

Finally, let $\gamma:[k]\to [j]$ be order-preserving. It is an easy check that $\gamma[2]^*$ and $\gamma[1]^*$ commute with source and target maps. Let $\alpha:[i]\hookrightarrow [k+2]$ and suppose first that $\gamma[2]\circ\alpha$ is injective. Then
\begin{enumerate}
    \item $\gamma[2]\circ\alpha\subseteq \delta^{k+2}$ if and only if $\alpha\subseteq \delta^{j+2}$, so in such case
    \[
    (\gamma[2]^*(A^0A^1))_\alpha=A^0_{\gamma[2]\circ\alpha}=(\gamma[2]^*A^0)_\alpha=((\gamma[2]^*A^0)(\gamma[2]^*A^1))_\alpha,
    \]

    \item $\gamma[2]\circ\alpha\subseteq \delta^{k+1}$ if and only if $\alpha\subseteq \delta^{j+1}$, so in such case
    \[
    (\gamma[2]^*(A^0A^1))_\alpha=A^1_{\gamma[2]\circ\alpha}=(\gamma[2]^*A^0)_\alpha=((\gamma[2]^*A^0)(\gamma[2]^*A^1))_\alpha,
    \]

    \item $\gamma[2]\circ\alpha=\varepsilon:[1]\hookrightarrow [j+2]$ if and only if $\alpha=\varepsilon:[1]\hookrightarrow [k+2]$, so in such case
    \begin{align*}
    (\gamma[2]^*(A^0A^1))_\alpha & =(A^0A^1)_\varepsilon=A^0_\varepsilon A^1_\varepsilon=(\gamma[2]^*A^0)_\varepsilon (\gamma[2]^*A^1)_\varepsilon \\
    &=((\gamma[2]^*A^0)(\gamma[2]^*A^1))_\alpha,
    \end{align*}

    \item if $\alpha$ satisfies none of the above conditions, then
    \begin{align*}
    (\gamma[2]^*(A^0A^1))_\alpha & =A^0_{\gamma[2]\circ\alpha} A^1_\varepsilon+A^1_{\gamma[2]\circ\alpha}=(\gamma[2]^*A^0)_\alpha (\gamma[2]^*A^1)_\varepsilon+(\gamma[2]^*A^1)_\alpha \\
    &=((\gamma[2]^*A^0)(\gamma[2]^*A^1))_\alpha.
    \end{align*}
\end{enumerate}
If $i=1$ and $\gamma[2]\circ\alpha$ is constant, then $\alpha\neq\varepsilon$, which means that either $\alpha\subseteq\delta^{k+1}$ or $\alpha\subseteq\delta^{k+2}$. If for instance $\alpha\subseteq\delta^{k+1}$, then
\[
(\gamma[2]^*(A^0A^1))_\alpha=\id=(\gamma[2]^*A^1)_\alpha=((\gamma[2]^*A^0)(\gamma[2]^*A^1))_\alpha,
\]
and similarly if $\alpha\subseteq \delta^{k+2}$. Lastly, if $\gamma[2]\circ\alpha$ is not injective but $i>1$, then a quick inspection of the definition of the product leads to
\[
(\gamma[2]^*(A^0A^1))_\alpha=0=((\gamma[2]^*A^0)(\gamma[2]^*A^1))_\alpha. \qedhere
\]
\end{proof}

\begin{corollary} \label{cor:GLinf(E)-as-simplicial-grpd}
$\GLi(E)[2]\rightrightarrows \GLi(E)[1]$ is an augmented simplicial Lie groupoid, where $\GLi(E)[i]$ is the simplicial set given by $\GLi(E)[i]_k:=\GLi(E)_{k+i}$ and structure maps $\gamma[i]^*:\GLi(E)[i]_j\to \GLi(E)[i]_k$, for $\gamma:[k]\to [j]$ order-preserving.
\end{corollary}

It can be helpful to have the following diagram of $\GLi(E)[2]\rightrightarrows \GLi(E)[1]$ in mind:
\[\begin{tikzcd}
	\dots & {\GLi(E)_3} & {\GLi(E)_2} & {\GLi(E)_1} \\
	\dots & {\GLi(E)_2} & {\GLi(E)_1} & {\GLi(E)_0.}
	\arrow[from=1-1, to=1-2]
	\arrow[shift left=2, from=1-1, to=1-2]
	\arrow[shift right=2, from=1-1, to=1-2]
	\arrow[shift right, from=1-2, to=1-3]
	\arrow[shift left, from=1-2, to=1-3]
	\arrow[shift left, from=1-2, to=2-2]
	\arrow[shift right, from=1-2, to=2-2]
	\arrow[from=1-3, to=1-4]
	\arrow[shift left, from=1-3, to=2-3]
	\arrow[shift right, from=1-3, to=2-3]
	\arrow[shift left, from=1-4, to=2-4]
	\arrow[shift right, from=1-4, to=2-4]
	\arrow[from=2-1, to=2-2]
	\arrow[shift left=2, from=2-1, to=2-2]
	\arrow[shift right=2, from=2-1, to=2-2]
	\arrow[shift left, from=2-2, to=2-3]
	\arrow[shift right, from=2-2, to=2-3]
	\arrow[from=2-3, to=2-4]
\end{tikzcd}\]

The augmented simplicial Lie groupoid $\GLi(E)[2]$ comes with a canonical representation.

\begin{definition}
Let $G\rightrightarrows N$ be an augmented simplicial Lie groupoid. A \textbf{representation} of $G$ is a simplicial vector bundle $V\to N$ with a representation of $G_k$ on $V_k\to N_k$ for each $k$ such that the structure maps of $V$ are equivariant, in the sense that for each $\gamma:[k]\to[j]$ order-preserving and $g\in G_j$, the following diagram commutes
\[
\begin{tikzcd}
	{V_{j,s(g)}} & {V_{j,t(g)}} \\
	{V_{k, s(\gamma^*g)}} & {V_{k,t(\gamma^*g)}.}
	\arrow["{g}", from=1-1, to=1-2]
	\arrow["{\gamma^*}"', from=1-1, to=2-1]
	\arrow["{\gamma^*g}"', from=2-1, to=2-2]
	\arrow["\gamma^*"', from=1-2, to=2-2]
\end{tikzcd}
\]
\end{definition}

\begin{proposition}
If we set
\[
\hat E_k:=\left( \bigoplus_{\alpha:[i]\hookrightarrow [k]} \pi_{\alpha(0)}^*E[-1-i]\right) \oplus \pi_{k+1}^*E\to M^{k+2},
\]
then $\hat E\to M^{\bullet+1}[1]$ becomes an augmented simplicial vector bundle as follows: if $\gamma:[k]\to[j]$ is order-preserving, then
\[
\begin{tikzcd}
	{\hat E_j} & {\hat E_k} \\
	{M^{j+2}} & {M^{k+2}}
	\arrow["{\gamma^*}",from=1-1, to=1-2]
	\arrow[from=1-1, to=2-1]
	\arrow["{\gamma[1]^*}"', from=2-1, to=2-2]
	\arrow[from=1-2, to=2-2]
\end{tikzcd}
\]
is the vector bundle map defined by $\gamma^*((v_\beta)_\beta,u)=(((\gamma^*v)_\alpha)_\alpha,u)$, where
\[
(\gamma^*v)_\alpha :=\begin{cases}
    v_{\gamma\circ\alpha}, & \text{if $\gamma\circ\alpha$ is injective,} \\
    0, & \text{otherwise.}
\end{cases}
\]
\end{proposition}

\begin{proof}
Clearly $\hat E_k\to M^{k+1}[1]$ is a vector bundle and the map $\gamma^*:\hat E_j\to \hat E_k$ is fiberwise linear and covers $\gamma[1]^*:M^{j+2}\to M^{k+2}$. The proof that $\hat E$ is simplicial follows along the same lines as the proof of Lemma \ref{lemma:simplicialization_simp_manifold}.
\end{proof}

\begin{definition} \label{def:action-of-GLinf}
We define the action of $\GLi(E)[2]$ on $p[1]^*\hat E\to\GLi(E)[1]$, where $p:\GLi(E)\to M^{\bullet+1}$ is the map defined in equation \eqref{eq:canonical_proj_p}, as follows: if $A\in \GLi(E)[2]_k$ and we write $(x_0,\dots, x_{k+2}):=p(A)$, then
\[
\begin{array}{cccc}
     A: & \displaystyle \left( \bigoplus_{\alpha:[i]\hookrightarrow [k]} E_{x_{\alpha(0)}}[-1-i] \right) \oplus E_{x_{k+2}} & \longrightarrow & \displaystyle \left( \bigoplus_{\alpha:[i]\hookrightarrow [k]} E_{x_{\alpha(0)}}[-1-i] \right) \oplus E_{x_{k+1}} \\
     & ((v_\alpha)_\alpha,u) & \longmapsto & ((v_\alpha+A_{\alpha[2]} u)_\alpha, A_\varepsilon u).
\end{array}
\]
\end{definition}

\begin{proposition} \label{prop:canonical-rep-of-GLinf(E)[2]}
$p[1]^*\hat E\to \GLi(E)[1]$ is a representation of $\GLi(E)[2]\rightrightarrows \GLi(E)[1]$.
\end{proposition}

\begin{proof}
Let $A^0,A^1\in\GLi(E)[2]_k$ be composable. Then
\begin{align*}
A^0\cdot (A^1\cdot ((v_\alpha)_\alpha,u)) &= A^0\cdot ((v_\alpha +A^1_{\alpha[2]} u)_\alpha, A^1_\varepsilon u) \\
&=((v_\alpha+A^1_{\alpha[2]} u+A^0_{\alpha[2]} A^1_\varepsilon u)_\alpha, A^0_\varepsilon A^1_\varepsilon u) \\
&=(A^0A^1)\cdot ((v_\alpha)_\alpha,u).
\end{align*}

Let now $\gamma:[k]\to[j]$ be order-preserving and $A\in\GLi(E)_{j+2}$ and $((v_\beta)_\beta,u)\in \hat E_{j, p(d_{k+1}A)}$. Then
\[
\gamma^*(A\cdot ((v_\beta)_\beta, u))=\gamma^*((v_\beta+A_{\beta[2]}u)_\beta, A_\varepsilon u)=((\gamma^*(v+A_{\cdot [2]}u))_\alpha)_\alpha, A_\varepsilon u)
\]
and
\[
\gamma[2]^*A\cdot \gamma^*((v_\beta)_\beta,u)=(((\gamma^*v)_\alpha+(\gamma[2]^*A)_{\alpha[2]} u)_\alpha, (\gamma[2]^*A)_\varepsilon u)=(((\gamma^*v)_\alpha+(\gamma[2]^*A)_{\alpha[2]} u)_\alpha, A_\varepsilon u).
\]
Consider $\alpha:[i]\hookrightarrow [k]$. If $\gamma\circ\alpha$ is not injective, then
\[
(\gamma^*(v+A_{\cdot [2]}u))_\alpha =0=(\gamma^*v)_\alpha+(\gamma[2]^*A)_{\alpha[2]}u,
\]
since $(\gamma[2]^*A)_{\alpha[2]}=0$. If, on the other hand, $\gamma\circ\alpha$ is injective, then
\begin{align*}
(\gamma^*(v+A_{\cdot [2]}u))_\alpha & = v_{\gamma\circ\alpha} +A_{(\gamma\circ\alpha)[2]}u= v_{\gamma\circ\alpha}+A_{\gamma[2]\circ\alpha[2]}u \\
&=(\gamma^*v)_\alpha+(\gamma[2]^*A)_{\alpha[2]}u. \qedhere
\end{align*}
\end{proof}

\begin{definition}
Let $G\rightrightarrows N$ be an augmented simplicial Lie groupoid. Its \textbf{base isotropy at $x\in N_0$} is the augmented simplicial Lie group $G_x$ given by
\[
(G_x)_k:= \begin{cases}
    (G_{-1})_{d_0x}, & \text{if $k=-1$}, \\
    (G_0)_x, & \text{if $k=0$}, \\
    (G_k)_{s_0^k(x)}, & \text{if $k>0$.}
\end{cases}
\]
If $V\to N$ is a representation of $G$, then the \textbf{base isotropy representation} is the induced representation of $G_x$ on $V_x$.
\end{definition}

Notice that in the case of $\GLi(E)$, for $x\in M$ and $\delta_x\in\End^1(E_x)$ we can consider the base isotropy $\GLi(E)[2]_{s_0\delta_x}$, which acts on the simplicial vector space $(p[1]^*\hat E)_{s_0\delta_x}$, explicitly given by
\[
((p[1]^*\hat E)_{s_0\delta_x})_k:=(p[1]^*\hat E)_{k, s_0^{k+1}\delta_x}=\hat E_{k,s_0^{k+1}x}=\left( \bigoplus_{\alpha:[i]\hookrightarrow [k]} E_x[-1-i]\right) \oplus E_x.
\]
We have the following immediate characterization, which is straightforward to prove.

\begin{proposition}
The base isotropy group $\GLi(E)[2]_{s_0\delta_x}$ is isomorphic to the augmented simplicial Lie group $\End(E_x)\rtimes\GL^0(E_x)$ defined as follows:
\[
(\End(E_x)\rtimes\GL^0(E_x))_k:=\left(\prod_{\alpha:[i]\hookrightarrow [k]}\End^{-1-i}(E_x)\right)\rtimes \GL^0(E_x),
\]
where $\GL^0(E_x)$ acts on $\End^{-1-i}(E_x)$ by $(H,g)\mapsto Hg$, and if $\gamma:[k]\to[j]$ is order-preserving, then $\gamma^*(H, g)=(\gamma^*H,g)$, where
\[
(\gamma^*H)_\alpha=\begin{cases}
    H_{\gamma\circ\alpha}, & \text{if $\gamma\circ\alpha$ is injective,} \\
    0, & \text{otherwise.}
\end{cases}
\]

In particular, the base isotropy representation of $\GLi(E)[2]_{s_0\delta_x}$ on $(p[1]^*\hat E)_{s_0\delta_x}$ is isomorphic to the representation of $\End(E_x)\rtimes\GL^0(E_x)$ on the simplicial vector space whose $k$-simplices are
\[
\left( \bigoplus_{\alpha:[i]\hookrightarrow [k]} E_x[-1-i] \right) \oplus E_x
\]
acting by $(H,g)\cdot (v,u)=((v_\alpha+H_\alpha u)_\alpha, gu)$.
\end{proposition}

\subsection{Simplicial holonomy groups} \label{sect:higher-holonomy:simplicial-holonomy}

We can now introduce the higher holonomy groups of a higher connection. In this section, we fix an $\Ascr$-connection $D$ on the graded vector bundle $E\to M$ and we denote by $\tau$ its higher parallel transport, whose base map is given by a section $\delta\in\Gamma(\End^1(E))$. We also fix a point $x\in M$.

A priori, we would like to say that the higher holonomy group of $D$ at $x$ is $\tau(\Pi_\infty(\Ascr)_x)$, that is, the parallel transport of all simplices all of whose vertices are $x$, but this has a clear problem: it is not in general a higher group. The reason is that a horn in $\tau(\Pi_\infty(\Ascr)_x)$ is not necessarily the image through $\tau$ of a horn in $\Pi_\infty(\Ascr)_x$. More concretely, given two algebraic simplices in $\tau(\Pi_\infty(\Ascr)_x)$ that can be glued along a face, it is not necessarily possible to glue the underlying geometric simplices in $\Pi_\infty(\Ascr)_x$. 

A way to overcome this issue is by restricting the simplices over which we do parallel transport, in such a way that the holonomy group becomes a subgroup of the simplicial isotropy $\GLi(E)[2]_{s_0\delta_x}$. Before we can do that, though, we need to understand how parallel transport behaves with respect to the composition of simplices.

There is a clear analogue of Definition \ref{def:1-gpd-structure-on-GLinfty} in $\Pi_\infty(\Ascr)$, i.e., a way of composing simplices, with the caveat that associativity and inverses only hold up to homotopy. To define it we will need the following subspaces of $\Delta^k$:
\begin{align*}
\Delta^k_+ &:=\left\{(t_1,\dots,t_k)\in\Delta^k\sep t_k\geq \frac{1}{2}t_{k-1}\right\}, \\
\Delta^k_- &:=\left\{(t_1,\dots,t_k)\in\Delta^k\sep t_k\leq \frac{1}{2}t_{k-1}\right\}.
\end{align*}
Notice that $\Delta^k_+\cup \Delta^k_-=\Delta^k$ and
\[
\Delta^k_+\cap\Delta^k_-=\left \{(t_1,\dots,t_k)\in\Delta^k\sep t_k=\frac{1}{2}t_{k-1}\right\}.
\]
Both of them are diffeomorphic to $\Delta^k$ itself, through diffeomorphisms $\varphi_\pm:\Delta^k_\pm\to\Delta^k$ given by
\begin{align*}
\varphi_+(t_1,\dots,t_k) &:=(t_1,\dots, t_{k-1},2t_k-t_{k-1}), \\
\varphi_-(t_1,\dots,t_k) &:=(t_1,\dots, t_{k-1},2t_k).
\end{align*}

\begin{definition}
For each $k>0$ we define a source, target and composition $\Pi_\infty(\Ascr)_k\rightrightarrows \Pi_\infty(\Ascr)_{k-1}$ as follows:
\begin{enumerate}
    \item the source map is $d_{k-1}$,
    \item the target map is $d_k$,
    \item composition is given as follows: if $\sigma_0,\sigma_1\in \Pi_\infty(\Ascr)_k$ are such that $d_k\sigma_1=d_{k-1}\sigma_0$, then we set
    \[
    \sigma_0\sigma_1:=\varphi_-^*\sigma_0\cup \varphi_+^*\sigma_1:T\Delta^k\to\Ascr.
    \]
\end{enumerate}
\end{definition}
\begin{remark}
    Technically, the above operation is not well-defined, since in directions transverse to $\Delta_+^k \cap \Delta_{-}^k$, $\sigma_0$ and $\sigma_1$ might not agree. Furthermore, even if it is well-defined, $\sigma_0\sigma_1$ need not be smooth in the transverse directions.

    To remedy this, the elements of $\Pi_\infty(\mathscr{A})_k$ will be allowed the have a finite number of singularities as described above along planes of the form 
    \begin{equation*}
        \{(t_1, \dots, t_k) \in \Delta^k : t_k = a t_{k-1} \}
    \end{equation*}
    with $a \in (0, 1)$.
    Note that the composition of any smooth path into $\Delta^k$ with an element of this extended $\Pi_\infty(\mathscr{A})$ is piecewise defined and smooth.
\end{remark}

The key point is that $\tau$ preserves this composition.

\begin{proposition} \label{prop:parallel-transport-preserves-composition}
If $\sigma_0,\sigma_1\in\Pi_\infty(\Ascr)_k$ are such that $d_k\sigma_1=d_{k-1}\sigma_0$, then $\tau(\sigma_0\sigma_1)=\tau(\sigma_0)\tau(\sigma_1)$.
\end{proposition}

\begin{proof}
Let $\nabla^i$ be any degree-preserving $T\Delta^k$-connection on $\sigma_i^*E\to\Delta^k$ such that $\delta_k^*\nabla^1=\delta_{k-1}^*\nabla^0$. Write $\omega_i:=\sigma_i^*D-d^{\nabla^i}\in \Omega(\Delta^k,\End(\sigma_i^*E))^1$, so that $\tau(\sigma_i)=\hat\Sigma(\omega_i)(\id)$. Let $\sigma:=\sigma_0\sigma_1$ and define $\nabla:=\varphi_-^*\nabla^0\cup \varphi_+^*\nabla^1$, which is a degree-preserving $T\Delta^k$-connection on $\sigma^*E\to \Delta^k$, and define as well $\omega:=\sigma^*D-d^\nabla\in\Omega(\Delta^k,\End(\sigma^*E))^1$. Observe that $\omega=\varphi_-^*\omega_0\cup \varphi_+^*\omega_1$. Then $\tau(\sigma)=\hat\Sigma(\omega)(\id)$.

By decomposing $\omega_i$ into its homogeneous parts we can write $\tau(\sigma_i)_\alpha$, for $\alpha:[j]\hookrightarrow [k]$, as a sum of integrals of the type
\[
\int_{\theta_j^*\alpha^*\sigma_i} \theta_j^*\alpha^*(\beta_i^1\otimes\dots\otimes \beta_i^\ell).
\]
We will do the argument in the case that $j=k$, $\alpha=\id$ and $\ell=1$, since the general case is readily seen to follow along similar lines.

The main observation is the following: if we denote by $\gamma_i(x)$ the path $t\mapsto \sigma_i\circ\theta_k(t,x)$, for $x\in I^{k-1}$, and similarly $\gamma(x)$ the path $t\mapsto \sigma\circ\theta_k(t,x)$, then
\[
\gamma(x)=\begin{cases}
    \gamma_1(x), & \text{if $x_{k-1}\geq \frac{1}{2}$,} \\
    \varepsilon_1\gamma_0(x), & \text{if $x_{k-1}\leq \frac{1}{2}$},
\end{cases}
\]
where $\varepsilon_1$ is the path $t\mapsto \sigma_1(1,\dots, 1, t)$. This can be easily seen from the definition of $\theta_k$. From this we conclude that
\begin{align*}
\int_{\theta_k^*\sigma} \theta_k^*\beta &= \int_{I^{k-1}}\int_0^1 \tau^{\theta_k^*\nabla}_{0\leftarrow t,x} \theta_k^*\beta(\del_1,\dots,\del_k)\tau^{\theta_k^*\nabla}_{t\leftarrow 1,x}dtdx \\
&=\int_{I^{k-1}}\int_0^1 \tau^{\theta_k^*\nabla^1}_{0\leftarrow t,x} \theta_k^*\beta_1(\del_1,\dots,\del_k) \tau^{\theta_k^*\nabla^1}_{t\leftarrow 1,x}dt dx \\
&\qquad + \int_{I^{k-1}}\int_0^1 \tau^{\theta_k^*\nabla^0}_{0\leftarrow t,x} \theta_k^*\beta_0(\del_1,\dots,\del_k)\tau^{\theta_k^*\nabla^0}_{t\leftarrow 1,x} dtdx \, \tau^{\theta_k^*\nabla^1}_{\varepsilon_1}.
\end{align*}
Hence, we conclude that
\[
\tau(\sigma)_{\id}=\tau(\sigma_1)_{\id}+\tau(\sigma_0)_{\id}\tau(\sigma_1)_\varepsilon=(\tau(\sigma_0)\tau(\sigma_1))_{\id}. \qedhere
\]
\end{proof}

\begin{definition}
The \textbf{simplicial holonomy group} of $D$ at $x$ we define to be the augmented simplicial set $\Hol^D_x$ whose $k$-simplices are
\[
(\Hol^D_x)_k:=\{\tau(\sigma)\in\GLi(E)_{k+2}: \sigma\in\Pi_\infty(\Ascr)_{k+2}, d_{k+2}\sigma=d_{k+1}\sigma=s_0^{k+1}(x)\},
\]
with structure maps inherited from the simplicial isotropy group $\GLi(E)[2]_{s_0\delta_x}$.
\end{definition}

\begin{theorem} \label{thm:holonomy-groups-are-smooth}
The simplicial holonomy group $\Hol^D_x$ is an augmented simplicial Lie subgroup of $\GLi(E)[2]_{s_0\delta_x}$.
\end{theorem}

\begin{proof}
It is enough to prove that $(\Hol^D_x)_k$ is a Lie subgroup of $(\GLi(E)[2]_{s_0\delta_x})_k$. By a result of Yamabe \cite{yamabe}, it suffices to prove that it is an arcwise connected subgroup of $(\GLi(E)[2]_{s_0\delta_x})_k$, and this implies that it is a Lie subgroup thereof.

That it is a subgroup follows directly from Proposition \ref{prop:parallel-transport-preserves-composition}. To see that it is arcwise connected, let $\sigma\in\Pi_\infty(\Ascr)_{k+2}$ with $d_{k+2}\sigma=d_{k+1}\sigma=s_0^{k+1}(x)$ and consider $\sigma_\epsilon:=\varphi_\epsilon^*\sigma$, where $\varphi_\epsilon:\Delta^{k+2}\to \Delta^{k+2}$ is $\varphi_\epsilon(t)=\epsilon t$, for $0\leq \epsilon\leq 1$. Then $\epsilon\mapsto \sigma_\epsilon$ is a continuous curve in $\Pi_\infty(\Ascr)_{k+2}$ with $d_{k+2}\sigma_\epsilon=d_{k+1}\sigma_\epsilon=s_0^{k+1}(x)$ such that $\sigma_1=\sigma$ and $\sigma_0$ is the trivial simplex over the constant simplex at $x$. Then $\epsilon\mapsto \tau(\sigma_\epsilon)$ is a continuous curve in $(\Hol^D_x)_k$ with $\tau(\sigma_1)=\tau(\sigma)$ and $\tau(\sigma_0)$ the identity of $(\GLi(E)[2]_{s_0\delta_x})_k$.
\end{proof}

\begin{definition}
The \textbf{holonomy representation} of $D$ at $x$ is the representation of $\Hol^D_x$ on $(p[1]^*\hat E)_{s_0\delta_x}=\left(\bigoplus_{\alpha:[i]\hookrightarrow [k]} E_x[1-i]\right) \oplus E_x$ inherited from $\GLi(E)[2]_{s_0\delta_x}$.
\end{definition}

\bigskip
\appendix

\section{Simplicializations of globular collections} \label{app:simplicialization}

In this appendix, we prove some basic properties of the simplicialization of a globular collection $X$ on $M$. To simplify proofs, we use the following notation. Let $\alpha :[i] \to [k]$ be an ordered map. Then, we write $|\alpha| := i$ and $\pi_i : M^{k + 1} \to M$ for the $(i + 1)$th projection from $M^{k +1}$ to $M$. Furthermore, if $i = 0$, then $\alpha^*X_0 := \pi^*_{\alpha(0)}X_0$. In case $i > 0$, then $\alpha^*X_i := (\pi_{\alpha(i)}, \pi_{\alpha(0)})^*X_i$. So, in particular, with this notation, the simplicialization of the globular collection $X$ is given by
\begin{equation*}
    S(X)_k = \bigoplus_{\alpha \in \underline{\Delta^k}}\alpha^*X_{|\alpha|}.
\end{equation*}

\begin{lemma}\label{lemma:simplicialization_simp_manifold}
    The simplicialization $S(X)$ of $X$ is a simplicial manifold and the projection map $p : S(X)\to M^{\bullet + 1}$ is a map of simplicial manifolds.
\end{lemma}
\begin{proof}
    In this proof, we will write $(X_i)_{(x, y)}$ for $q^{-1}_i(\{(x, y\})$ with $q_i : X_i \to M \times M$ as in the definition \ref{def:globular-collection}.
    
    Let $\gamma : [k] \to [j]$ be an order preserving map. Then, $\gamma^*$ is smooth, as all its components are smooth. To see that the $\alpha$-component of $\gamma^*$ actually lies in $\alpha^*X_{|\alpha|}$, consider $A\in S(X)_j$ with $x = p_j(A)$ and let $\alpha:[i]\hookrightarrow [k]$. Then, if $\gamma \circ \alpha$ is injective we have
    \[
    A_{\gamma \circ \alpha} \in (X_i)_{(x_{\gamma \circ \alpha(0)}, x_{\gamma \circ \alpha(i)})} = (X_i)_{((\gamma^*x)_{\alpha(0)}, (\gamma^*x)_{\alpha(i)})}.
    \]
    If $i = 1$ and $\gamma \circ \alpha$ is constant, then we see that 
    \[
    u_1(\pi_{\gamma \circ \alpha(0)}\circ p_j(A)) \in (X_1)_{(x_{\gamma\circ\alpha(0)}, x_{\gamma\circ\alpha(0)})} = (X_1)_{((\gamma^*x)_{\alpha(0)}, (\gamma^*x)_{\alpha(1)})}
    \]
    since $u_1$ covers the diagonal map and $\gamma \circ \alpha(0) = \gamma\circ \alpha(1)$. Finally, if $\alpha$ satisfies neither of the above criteria, we have
    \[
    u_i(\pi_{\gamma \circ \alpha(0)}\circ p_j(A), \pi_{\gamma \circ \alpha(i)}\circ p_j(A)) \in (X_i)_{(x_{\gamma\circ \alpha(0)}, x_{\gamma \circ \alpha(i)})} = (X_i)_{((\gamma^*x)_{\alpha(0)}, (\gamma^*x)_{\alpha(i)})}.
    \]
    This shows that $\gamma^*$ is indeed well-defined. 
    
    Next up, note that the diagram (\ref{cd:simplicialization_proj}) commutes by definition of $\gamma^*$.
    
    Finally, let us show that $S(X)(\cdot)$ is functorial. Consider, then,
    \begin{equation*}
        [i] \stackrel{\gamma}{\to} [j] \stackrel{\delta}{\to} [k]
    \end{equation*}
    and let $A \in S(X)_k$. Then, we have to show that 
    \begin{equation} \label{eq:sx_functorial}
        ((\delta \circ \gamma)^*A)_{\alpha} =(\gamma^*\delta^*A)_{\alpha}.
    \end{equation}
    for all $\alpha :[m] \hookrightarrow [i]$. First, note that if $\delta \circ \gamma \circ \alpha$ is injective, then $\gamma \circ \alpha$ is also injective, and therefore
    \begin{equation*}
        ((\delta \circ \gamma)^*A)_{\alpha} = A_{\delta \circ \gamma \circ \alpha} = (\delta^*A)_{\gamma \circ \alpha}  = (\gamma^*\delta^*A)_\alpha.
    \end{equation*}
    Thus, in this case \eqref{eq:sx_functorial} holds. 
    
    Now, suppose that $\delta \circ \gamma \circ \alpha$ is not injective with $m = 1$. Then, for the left hand side of \eqref{eq:sx_functorial} we find
    \begin{equation*}
        ((\delta \circ \gamma)^*A)_{\alpha} =  u_1(\pi_{\delta \circ \gamma \circ \alpha(0)}\circ p_k(A)).
    \end{equation*}
    For the right hand side, there are two options. On the one hand, if $\gamma \circ \alpha$ is not injective, then we find
    \begin{equation*}
        (\gamma^*\delta^*A)_\alpha = u_1(\pi_{\gamma \circ \alpha(0)}\circ p_j(\delta^*A)) = u_1(\pi_{\delta \circ \gamma \circ \alpha(0)}\circ p_k(A)).
    \end{equation*}
    On the other hand, if $\gamma \circ \alpha$ is injective, we instead find
    \begin{equation*}
        (\gamma^*\delta^*A)_\alpha = (\delta^*A)_{\gamma \circ \alpha} = u_1(\pi_{\delta \circ \gamma \circ \alpha(0)}\circ p_k(A)),
    \end{equation*}
    as $\delta \circ \gamma \circ \alpha$ is not injective. So, in this case, \eqref{eq:sx_functorial} holds as well. The final case of $\delta \circ \gamma \circ \alpha$ not injective and $m > 1$ is analogous to the $m = 1$ case.  
\end{proof}

Let $A \subseteq \Delta^k$ be simplicial subset of $\Delta^k$ generated by some non-degenerate elements $\alpha_i \in (\Delta^k)_{a_i}$ with $i = 1, \dots, \ell$. Then, for a simplicial set $G$, we have a canonical inclusion
\[
\begin{array}{cccc}
 i_{A,G}:&  \Hom_{\sSet}(A, G) & \longrightarrow & \displaystyle \prod_{i = 1}^\ell G_{a_i} \\
    & \phi &\longmapsto &(\phi(\alpha_i))_i.
\end{array}
\]
Note this is an inclusion as a map of simplicial sets. It is fully determined by what it does on generators.
\begin{lemma}
    For $G = M^{\bullet + 1}$, 
    \begin{equation}
        \Hom_{\sSet}(A, M^{\bullet + 1}) \subseteq \prod_{i = 1}^\ell M^{a_i + 1}
    \end{equation}
    is an embedded submanifold and upon writing $A_0 = \{i_0, \dots, i_m\} \subseteq [k]$, we have an isomorphism
    \begin{equation*}
        \Hom_{\sSet}(A, M^{\bullet + 1})\longrightarrow M^{|A_0|}
    \end{equation*}
    given by 
    \begin{equation*}
        \phi \longmapsto (\phi(i_0), \dots, \phi(i_m)).
    \end{equation*}
    Finally, under this isomorphism, we have 
    \[
    i_{A,M^{\bullet+1}}(x_{i_0}, \dots, x_{i_m}) = (x_{\alpha_1(0)},\dots, x_{\alpha_1(a_1)},\dots, x_{\alpha_\ell(0)},\dots, x_{\alpha_\ell(a_\ell)}).
    \]
\end{lemma}
\begin{proof}
    We claim that $\phi \in \Hom_{\sSet}(A, M^{\bullet + 1})$ is completely fixed by where it sends $A_0 \subseteq(\Delta^k)_0$. To see this, note that for $\alpha : [i] \to [k]$, we have $\pi_j(\phi(\alpha)) = \phi(\alpha(j))$ for $j \in [i]$ and $\alpha(j) \in A_0$. So, $\Hom_{\sSet}(A, M^{\bullet + 1}) \cong M^{m + 1}$, with the isomorphism given by $\phi \mapsto (\phi(i_0), \dots, \phi(i_m))$. Then, since $i_{A,M^{\bullet +1}}$ is given by evaluating on the non-degenerate simplices $\alpha_i$, we see that
    \[
    i_{A,M^{\bullet+1}}(x_{i_0}, \dots, x_{i_m}) = (x_{\alpha_1(0)},\dots, x_{\alpha_1(a_1)},\dots, x_{\alpha_\ell(0)},\dots, x_{\alpha_\ell(a_\ell)}),
    \]
    which is an embedding. 
\end{proof}

Note that $A_0$ determines an inclusion $A_0 : [|A_0|-1] \to [k]$. Furthermore, any $\alpha : [i] \to [k] \in A$ factors through $A_0$, as some $\alpha_A$. For $i \in A_0$, we denote by $\pi^A_{i} :M^{|A_0|} \to M$ the projection $(x_j)_{j \in A_0} \to x_i$. Then, for $\alpha : [i] \to [k]$ in $A$, we have 
\[
    \alpha_A^*X = \begin{cases}
    (\pi_{\alpha(i)}^A,\pi_{\alpha(0)}^A)^*X_i &\text{if } i > 0, \\
    (\pi_{\alpha(0)}^A)^* X_i&\text{if } i = 0.
    \end{cases}
\]
\begin{lemma} \label{lemma:horn_of_end}
    The map $i_{A,S(X)}$ covers $i_{A,M^{\bullet+1}}$. Furthermore, under the map $i_{A,S(X)}$, $\Hom_{\sSet}(A, S(X))$ is an embedded submanifold of $\prod_{i = 1}^\ell S(X)_{a_i}$, which is diffeomorphic to 
    \begin{equation*}
        \bigoplus_{\alpha \in \underline{A}}\alpha^*_AX_{|\alpha|}
    \end{equation*}
    through the embedding
    \begin{align*}
        \iota : \bigoplus_{\alpha \in \underline{A}}\alpha^*_AX_{|\alpha|} & \longrightarrow \prod_{i = 1}^\ell S(X)_{a_i} \\
        (x, A) & \longmapsto (({\alpha}_i)_A^*(x), (A_{\beta^*(\alpha_i)})_{\beta \in \underline\Delta^{a_i}})_{i = 1}^\ell.
    \end{align*}
\end{lemma}
\begin{proof}
    Let us first show that $\iota$ gives an embedding. First, remark that by construction $\bigoplus_{\alpha \in \underline{A}} \alpha^*_AX_{|\alpha|}$ sits over $M^{|A_0|}$ and that $\iota$ covers $i_{A,M^{\bullet +1}}$. So, $\iota$ factors through the pullback of
    \begin{equation*}
        \prod^\ell_{i = 1} S(X)_{a_i} \to \prod^\ell_{i=1}M^{a_i + 1}
    \end{equation*}
    along the embedding $i_{(A, M^{\bullet + 1})}$. Thus, it suffices to show that the map 
    \begin{equation*}
        \bigoplus_{\alpha \in \underline{A}} \alpha^*_AX_{|\alpha|} \to i_{(A, M^{\bullet + 1})} ^*\prod_{i=1}^\ell S(X)_{a_i}
    \end{equation*}
    is an embedding. To do this, we construct a smooth retract. Since the non-degenerate elements $\alpha_i$ generate $\underline{A}$, we may fix for each $\alpha : [m] \to [k] \in \underline{A}$, an $i_\alpha \in \{1, \dots, \ell\}$ and a map $\beta_\alpha : [m] \to [a_i]$ such that $\beta^*_\alpha(\alpha_{i_\alpha}) = \alpha$. Then, as $\alpha$ is injective, $\beta_\alpha$ is injective. So, we can define
    \begin{align*}
        r : i_{A,M^{\bullet +1}}^*\prod_{i =1 }^l S(X)_{a_i} & \longrightarrow \bigoplus_{\alpha \in \underline{A}}\alpha^*X_{|\alpha|} \\
        (x, A) &\longmapsto (x, ((A_{i_\alpha})_{\beta_\alpha})_{\alpha \in \underline{A}}),
    \end{align*}
    where we note that the base point of $A$ in $S(X)_{a_i}$ is fully determined by $x$, as $\alpha_A^*(x)$. Then, we have 
    \begin{align*}
        r\iota(x, A) & = r(((\alpha_i)^*_A(x), (A_{\beta^*(\alpha_i)})_{\beta \in \underline{\Delta}^{a_i}})^\ell_{i = 1}) \\
        & = (x, (A_{\beta_\alpha^*(\alpha_{i_\alpha})})_{\alpha \in \underline{A}}) \\
        & = (x, A).
    \end{align*}
    Finally, let us show that the image of $\iota$ coincides with the image of $i_{A,S(X)}$. For this, we will first understand the image of $i_{A,S(X)}$ better. Let us write $Q =  \bigsqcup_{i = 1}^\ell \Delta^{a_i}$. Now, since the $\alpha_i$ generate $A$, the map
    \[
    \begin{array}{cccc}
        g: & Q & \longrightarrow & A \\
        & \beta_{i} &\longmapsto & \alpha_i \circ \beta_{i}
    \end{array}
    \]
    is surjective on objects. Thus, $g$ is an epimorphism in simplicial sets, and we may form the kernel pair $R$ of $g$, given by
    \begin{equation*}
        R = \{(\beta_i, \beta_j) \in Q \times Q : \alpha_i \circ \beta_i = \alpha_j \circ \beta_j\}.
    \end{equation*}
    Then, $A$ is the coequalizer of $\pi_1 : R \to Q$ and $\pi_2 : R \to Q$. Since colimits commute with the contravariant functor $\Hom_{\sSet}(\cdot, S(X))$, we see that $i_{A,S(X)}$ in  
    \begin{center}
        \begin{tikzcd}
            \Hom_{\sSet}(A, S(X)) \ar[rr, "i_{A,S(X)}"] && \prod_{i =1}^\ell S(X)_{a_i} \ar[r, shift right = 1, "\pi_2^*"'] \ar[r, shift left = 1, "\pi_1^*"] & \Hom_{\sSet}(R, S(X))
        \end{tikzcd}
    \end{center}
    is an equalizer. Let us write $E$ for the image of $i_{A,S(X)}$ in $\prod_{i = 1}^\ell S(X)_{a_i}$. Note that for $A \in \prod_{i= 1}^\ell S(X)_{a_i}$, we have $A \in E$ if and only if $\pi_1^*(A) = \pi_2^*(A)$, and this happens if and only if $\beta_i^*(A_i) = \beta_j^*(A_j)$ for all $(\beta_i, \beta_j) \in R$. We now claim that $\iota$ defines a bijection with $E$. 
    
    Firstly, note that for $(\beta_i, \beta_j) \in R_m$, we have 
    \begin{align*}
        \beta_i^*(\iota(x, A)_i) & = \beta_i^*\left( (A_{\beta^*\alpha_i})_{\beta\in\underline{\Delta}^{a_i}} \right) =(A_{\gamma^*(\alpha_i \circ \beta_i)})_{\gamma \in \underline{\Delta}^{m}}  \\
        &= (A_{\gamma^*(\alpha_j \circ \beta_j)})_{\gamma \in \underline{\Delta}^{m}} = \beta_j^*(\iota(x, A)_j).
    \end{align*}
    So, $\iota$ factors through $E$. On the other hand, we claim that $\iota r$ restricted to $E$ is the identity. To see this, note that for $A \in E$, we have
    \begin{align*}
        \iota r(A)_i & = ((A_{i_{\gamma^*(\alpha_i)}})_{\beta_{\gamma^*(\alpha_i)}})_{\gamma \in \underline{\Delta}^{a_i}} \\
        & = ((\beta_{\gamma^*(\alpha_i)}^*A_{i_{\gamma^*(\alpha_i)}})_{\id})_{\gamma \in \underline{\Delta}^{a_i}}.
    \end{align*}
    Then, since $\alpha_{i_{\gamma^*(\alpha_i)}} \circ \beta_{\gamma^*(\alpha_i)} = \gamma^*(\alpha_i) = \alpha_i \circ \gamma$, we have that $(\beta_{\gamma^*(\alpha_i)},\gamma)\in R$, which implies that, since $A \in E$,
    \[
    \beta^*_{\gamma^*(\alpha_i)}A_{i_{\gamma^*(\alpha_i)}}=\gamma^*A_i.
    \]
    Therefore,
    \begin{align*}
        \iota r(A)_i = ((\gamma^*A_i)_{\id})_{\gamma \in \underline{\Delta}^{a_i}} = ((A_i)_{\gamma})_{\gamma \in \underline{\Delta}^{a_i}} = A_i.
    \end{align*}
    We conclude that $E$ coincides with image of $\iota$.
\end{proof}

\begin{proposition} \label{prop:simplicialization_groupoid}
    The simplicialization $S(X)$ of $X$ is a higher Lie groupoid. Furthermore, if $X$ has length $n$, then $S(X)$ is a Lie $(n + 1)$-groupoid.
\end{proposition}
\begin{proof}
    Note that, by Lemma \ref{lemma:horn_of_end}, the $(k,\ell)$-horns $S(X)_{k, \ell} \subseteq \prod_{i \neq \ell}S(X)_{k-1}$ are embedded submanifolds. Furthermore, the $(k, \ell)$-horn map is given by the projection
    \begin{equation*}
        d_{k,\ell}: \bigoplus_{\alpha \in \underline{\Delta_k}} \alpha^*X_{|\alpha|} \to \bigoplus_{\alpha \in \underline{\Lambda_{k, \ell}}}  \alpha^*X_{|\alpha|}
    \end{equation*}
    that leaves out the $\id$ and $\delta^\ell$ components, and it covers the horn map $d_{k, \ell}^M : M^{k + 1} \to (M^{\bullet +1})_{k, \ell}$ of $M^{\bullet + 1}$. Thus, $d_{k,\ell}$ is a surjective submersion, so $S(X)$ is a higher Lie groupoid.

    Let us now show that $S(X)$ is a Lie $(n + 1)$-groupoid if $X$ has length $n$. Notice that, since $X$ has length $n$, we have that $X_{|\alpha|} = M \times M \to M \times M$ for all $\alpha:[i]\hookrightarrow [k]$ with $i > n \geq 0$. Hence,
    \begin{equation*}
        S(X)_{n +1} = \bigoplus_{\alpha \in \underline{\partial \Delta^{n+1}}} \alpha^*X_{|\alpha|}.
    \end{equation*}
    On the other hand, we have
    \begin{equation*}
        S(X)_{n+1, \ell} = \bigoplus_{\alpha \in \underline{\Lambda_{n+1, \ell}}} \alpha^*X_{|\alpha|}.
    \end{equation*}
    So, these two spaces differ by a term $(\delta^\ell)^*X_{n}$. Since $q_n : X_n \to M \times M$ is not a diffeomorphism, we see that the $(n+1, \ell)$-horn map $S(X)_{n+1} \to S(X)_{n + 1, \ell}$ is not a diffeomorphism. However, for $k > n + 1$, we see that for all $0 \leq \ell \leq k$ we have that
    \begin{equation*}
        S(X)_{k} = \bigoplus_{\alpha \in (\underline{\Delta^k})_{\leq n}} \alpha^*X_{|\alpha|} \cong \bigoplus_{\alpha \in (\underline{\Lambda_{k,\ell}})_{\leq n}} \alpha^*X_{|\alpha|} = S(X)_{k, \ell}, 
    \end{equation*}
    where $\underline{A}_{\leq n}$ denotes the non-degenerate simplices $x \in A_d$ with $d\leq n$.
\end{proof}

\begin{proposition} \label{prop:simplicialization_contractible}
    The simplicialization $S(X)$ of $X$ is contractible.
\end{proposition}

\begin{proof}
    Note that $\partial_k : S(X)_k \to \partial_k S(X)$ is the projection,
    \begin{equation*}
        \bigoplus_{\alpha \in \underline{\Delta^k}} \alpha^*X_{|\alpha|} \to \bigoplus_{\alpha\in \underline{\partial\Delta^k}} \alpha^*X_{|\alpha|}
    \end{equation*}
    that leaves out the $\id$ component. This is a surjective submersion, so $S(X)$ is contractible.
\end{proof}

\begin{proposition} \label{prop:simplicialization_vector_bundle}
    Suppose that the bundles making the globular collection $X$ are vector bundles $X_0 \to M$ and $X_k \to M \times M$ and that the unit maps are the $0$-sections of these vector bundles. Then, $S(X)\to M^{\bullet +1}$ is a simplicial vector bundle. 
\end{proposition}
\begin{proof}
    Note that in case $X_i$ are vector bundles, $S(X)_k \to M^{k+1}$ is constructed as a direct sum of vector bundles. To see that $S(X)$ is a simplicial vector bundle, it suffices to check that the structure maps of $S(X)$ are linear. If the unit maps are the zero sections, then \eqref{eq:structure_maps_simplicialization} becomes
    \begin{equation*}
    (\gamma^*A)_\alpha=\begin{cases}
        A_{\gamma\circ\alpha}, & \text{if $\gamma\circ\alpha$ is injective,} \\
        0_{\pi_{\gamma\circ\alpha(0)}\circ p_j(A)}, & \text{if $i=1$ and $\gamma \circ \alpha$ is constant,} \\
        0_{(\pi_{\gamma\circ\alpha(0)}\circ p_j(A),\pi_{\gamma\circ\alpha(i)}\circ p_j(A))}, & \text{otherwise,}
    \end{cases}
    \end{equation*}
    which is linear. 
\end{proof}

\section{Diffeological spaces} \label{app:diffeology}

Diffeological spaces are one of the possible ways of weakening smooth structures so that spaces of maps become smooth spaces. The advantage of diffeological spaces is that it allows one to speak these smooth structures without having to perform heavy analysis. They were introduced originally in \cite{souriau}. See \cite{iglesias-zemmour} for a comprehensive introduction.

\begin{definition}
A \textbf{diffeology} on a set $X$ is a collection of maps, called the set of \textbf{plots}, which is subset of
\[
\{U\to X\sep \text{$n\in\N$ and $U\subseteq \R^n$ is open}\}
\]
satisfying the following properties:
\begin{enumerate}
    \item any constant map $\R^n\to X$ is a plot,
    \item if $p:U\to X$ is such that for every $x\in U$ there is an open set $V$ in $U$ containing $x$ such that $p|_V$ is a plot, then $p$ is a plot,
    \item if $p:U\to X$ is a plot and $F:V\to U$ is smooth, for some open $V\subseteq\R^m$, then $p\circ F$ is a plot.
\end{enumerate}

If $X$ and $Y$ are diffeological spaces, then a map $f:X\to Y$ is \textbf{smooth} if for every plot $p:U\to X$ the composition $f\circ p:U\to Y$ is a plot.
\end{definition}

\begin{example}
If $M$ is manifold, then setting the plots to be the smooth maps $U\to M$ presents $M$ as a diffeological space. 
\end{example}

\begin{example} \label{ex:functional-spaces-are-smooth}
If $M$ and $N$ are manifolds, then $\Cinf(M,N)$ is a diffeological space, for which $p:U\to\Cinf(M,N)$ is a plot if the induced map
\[
\begin{array}{ccc}
U\times M & \longrightarrow & N \\
(x,y) & \longmapsto & p(x)(y)
\end{array}
\]
is smooth. We can think of such plots $p$ as smooth families of smooth maps $M\to N$. From this perspective, if $P$ is another manifold, then a map $F:\Cinf(M,N)\to P$ is smooth if and only if it is ``smooth along smooth families'', in the sense that for any smooth family of smooth maps $p:U\to \Cinf(M,N)$ the composition $F\circ p:U\to P$ is smooth.
\end{example}

\begin{example}
As a particular case of Example \ref{ex:functional-spaces-are-smooth}, if $E\to M$ is a vector bundle then $\Gamma(E)$ is a diffeological space. In the case that $E=\wedge T^*M$, we have the following explicit description: a map $p:U\to \Omega^k(M)$ is smooth if and only if for every tuple $X_1, \dots, X_k\in\mathfrak{X}(M)$ the map $U\times M\to \R$ given by
\[
(x,y)\longmapsto \big(p(x)(X_1,\dots,X_k)\big) (y)
\]
is smooth.
\end{example}

\begin{example}\label{ex:pullback}
    Given a diffeological space $Y$ and a map of sets $f : X \to Y$, the \textbf{pullback diffeology} on $X$ has as plots maps $p : U \to X$ such that $f \circ p : U \to Y$ is a plot, where $U \subseteq \R^n$ is open. Specifically if $X \subseteq Y$, we see that $p : U \to X$ is a plot if and only if $p$ defines a plot when considered as a map $U \to Y$.
\end{example}

\begin{example}\label{ex:pushforward}
    Given a diffeological space $X$ and a map of sets $f : X \to Y$, the \textbf{pushforward diffeology} on $Y$ has as plots maps $p : U \to Y$ such that locally $p$ is of the form $f \circ q$ for $q : V \to X$ a plot of $X$, i.e., for each $x \in U$ there exists an open neighborhood $V \subseteq U$ of $x$ and a plot $q : V \to X$ such that $p|_V = f \circ q$.
\end{example}

\bibliographystyle{alpha}
\bibliography{biblio}

\end{document}